\documentclass[preprint]{elsarticle}

\usepackage[margin=1in]{geometry}
\usepackage{mathtools}
\usepackage{amsthm}
\usepackage{amssymb}

\usepackage{color}

\usepackage{graphicx}

\newtheorem{theorem}{Theorem}

\theoremstyle{definition}
\newtheorem{definition}{Definition}
\theoremstyle{remark}

\DeclarePairedDelimiter{\abs}{\lvert}{\rvert}
\DeclarePairedDelimiter{\norm}{\lVert}{\rVert}

\let\div\relax 
\DeclareMathOperator{\div}{div}
\DeclareMathOperator{\diam}{diam} 
\DeclareMathOperator{\spanset}{span} 
\DeclareMathOperator{\dist}{dist} 
\DeclareMathOperator*{\argmin}{arg\,min}

\begin{document}

\begin{abstract}
In this work we address the Multiscale Spectral Generalized Finite Element Method (MS-GFEM) developed in [I. Babu\v{s}ka and R. Lipton, Multiscale Modeling and Simulation 9 (2011), pp. 373--406]. We outline the numerical implementation of this method and present simulations that demonstrate contrast independent exponential convergence of MS-GFEM solutions. We introduce strategies to reduce the computational cost of generating the optimal oversampled local approximating spaces used here. These strategies retain accuracy while reducing the computational work necessary to generate local bases. Motivated by oversampling we develop a nearly optimal local basis based on a partition of unity on the boundary and the associated A-harmonic extensions.
\end{abstract}

\title{Multiscale-Spectral GFEM and Optimal Oversampling}
\tnotetext[t1]{Authors listed alphabetically.}
\tnotetext[t2]{This work is supported by grants: NSF DMS-1211066, NSF DMS-1211014, and NSF DMS-1813698.}

\cortext[cor1]{Corresponding author}
\author[1]{Ivo Babu\v{s}ka}
\ead{babuska@ices.utexas.edu}
\author[2]{Robert Lipton}
\ead{lipton@math.lsu.edu}
\author[3]{Paul Sinz}
\ead{sinzpaul@msu.edu}
\author[4]{Michael Stuebner\corref{cor1}}
\ead{mstuebner@gem-innovation.com}

\address[1]{Department of Mathematics, Department of Aerospace Engineering and Engineering Mechanics, University of Texas at Austin}
\address[2]{Department of Mathematics, Louisiana State University}
\address[3]{Department of Computational Mathematics, Science and Engineering, Michigan State University}
\address[4]{Global Engineering and Materials, Inc.}


\maketitle

\section{Introduction}
Modern science and technology demand  increasingly efficient methods for solution of  multiscale heterogeneous problems. 
These problems arise in structural analysis for aerospace and infrastructure or naturally in the study of biological or geological structures. The primary computational challenge in all of these problems is due to the extreme degrees of freedom associated with parameterizing material heterogeneity. Multiscale numerical methods provide one way to solve this problem by invoking  local and independent computations that enable a global solve involving a drastically reduced number of degrees of freedom.

A typical and important example is local stress analysis inside an epoxy  fiber reinforced composite structure. Here the fibers have a diameter on the order of 10 microns and are densely distributed within an epoxy filler. A critical numerical analysis of such a problem is only possible on massively parallel computers. Hence a special and very parallelizable numerical method is needed. This paper addresses such a method.
In this computational investigation we work with the generalized finite element method (GFEM) introduced in \cite{BC94} and expanded on  in \cite{BB03}, \cite{BM97}, and \cite{SB01}. This approach is a partition of unity method (PUM) \cite{BC94}, which utilizes the results of many independent and local computations carried out across the computational domain. Essentially this is  a domain decomposition and the computational domain is partitioned into an overlapping collection of preselected subsets $\omega_i$, $i = 1,2,...m$. Finite-dimensional approximation spaces $V_{\omega_i}$ are constructed over each subset using local information. In this paper the finite dimensional space of local solutions to the problem is first solved over a slightly larger domain $\omega^*_i$ and generates the finite dimensional space $V_{\omega^*_i}$. In the second step we restrict  the solutions to the smaller set (patch) $\omega_i\subset\subset\omega_i^*$. In this way the local space $V_{\omega_i}$ is now given by the restriction of elements in $V_{\omega^*_i}$ to $\omega_i$. This way of generating the local basis is called oversampling. As is the case in all GFEM schemes, each subspace is computed independently and the full global solution is obtained by solving a global (macro) system which can be orders of magnitude smaller than the system corresponding to a direct application of the finite element method (FEM) to the full structure. 

In this paper we carry out numerical investigations using optimal local bases introduced in \cite{BL11} and \cite{BH14}. 
These bases are proved to have the best approximation properties over all local bases constructed by oversampling A-harmonic functions \cite{BL11}, \cite{BH14}. Here an A-harmonic function $v$ satisfies 
\begin{equation}
{\rm div}(A(x)\nabla v)=0 \hbox{ on $\omega^\ast$},
\end{equation}
where $A(x)$ is a given $L^\infty$ coefficient satisfying the usual coercivity and boundedness conditions.
The N-dimensional optimal local basis is shown to be the span of the first $N$ eigenfunctions associated with the singular values of the restriction operator used in the oversampling. These local bases have been theoretically proved to have approximation error that asymptotically approaches  an exponential decay with the dimension of the approximation space (\cite{BL11} and \cite{BH14}). The convergence result applies to any elastic or conductivity tensor field with $L^\infty$ coefficients satisfying coercivity and boundedness conditions. 
The associated scheme introduced in \cite{BL11} and \cite{BH14} is the Multiscale-Spectral GFEM (MS-GFEM) and has global approximation error that converges exponentially asymptotically in the energy norm.

As discussed above the ultimate motivation behind the MS-GFEM and its exponential convergence rate is its potential for use in large parallel implementations. The method has two essential parts. A. Parallel construction of the local non-polynomial  shape functions leading to an exponential approximation of the actual solution in the energy norm. Here we note that each local basis is defined over separate patches and can be computed independently on separate processors using local memory. B. Efficient parallel solution of the global system of linear algebraic equations. This part is made efficient noting that
global basis vectors associated with each patch are only influenced by immediate neighboring patches and a Schwartz alternating method \cite{B1958}, \cite{Lions88} can be applied as a preconditioner for an exponentially convergent iterative scheme. The current paper concentrates on A while B is discussed in section \ref{iterative} and its details will be presented in a forthcoming paper.  Another very expensive aspect of large implementations is mesh generation. Here the GFEM allows each patch to be meshed independently. 

In this work we provide a numerical implementation of MS-GFEM and our aim is threefold. We focus on a simple two patch domain decomposition to illustrate 1) contrast independent convergence of the local basis and global basis, 2) new methods for inexpensively generating local basis functions, and 3) identification of alternate local bases with good approximation properties. The computational domain is composed of heat conducting  particles included within a connected second heat conducting phase. The connected phase is often referred to as the matrix phase. Computations are carried out here for a large range of contrasts between particle and matrix conductivities. For inclusions that do not touch and have smooth boundary our computations show that the MS-GFEM approximations converge exponentially and independently of the contrast between phases in the preasymptotic regime. This is consistent with the simulations given in \cite{SZ03}  for a medium with separated holes. This is also consistent with the more recent numerical investigation carried out for scalar and elastic problems using the optimal local basis functions in \cite{SW17}. The theoretical proof of contrast independent exponential convergence for MS-GFEM applied to heterogeneous elastic and conduction problems is given in a forthcoming publication.

The primary numerical work in generating the best local basis is not in the solution of the eigenvalue problem but in the numerical generation of the local A-harmonic subspace used to construct the optimal local basis.  
In section \ref{sec:aharmonic} we introduce a method to reduce the computational cost of generating optimal local basis functions. We investigate the inexpensive construction of A-harmonic subspaces for numerically generating the optimal local basis. The optimal $N$-dimensional approximation space is given by the span of the optimal local basis given by the first $N$ eigenfunctions associated with the largest $N$ singular values of the restriction operator acting on A-harmonic functions. These eigenfunctions are seen to define the Kolmogorov $n$-width \cite{BL11} and are referred to as the $n$-width eigenfunctions. 

We explore the use of discrete A-harmonic extensions of hat functions defined on the boundary of $\omega^\ast$ necessary to generate these eigenfunctions numerically. Here we are motivated by the fact that longer wavelength boundary data applied to $\partial\omega^*$ penetrates further into the subdomain $\omega$ than short wavelength ``oscillatory'' boundary data \cite{BLS08} and \cite{LS16}. In light of the exponential decay of the $n$-width eigenfunctions  it becomes clear that the space spanned by the lower $n$-width modes are generated by linear combinations of A-harmonic extensions of low frequency or slowly changing boundary data. A-harmonic extensions of more oscillatory boundary data decay faster into the domain and thus contribute only to the higher $n$-width modes. Thus we use the A-harmonic extension of a collection of hat functions with a relatively large support set that taken together form a partition of unity on the boundary.  Here the support of the boundary data given by hat functions can extend over tens or hundreds of boundary nodes associated with the discrete A-harmonic extension of the boundary data. Our numerical experiments reveal the size of support of the boundary hat function that can be used for accurate yet inexpensive computation of the optimal local basis.  It turns out that the number of independent boundary hat functions can be roughly the same as the number of discrete $n$-width basis functions needed for a good approximation.  The numerical experiments presented in Section \ref{sec:aharmonic} show that this method retains accuracy while effectively reducing the computational work necessary to generate local bases. 
Motivated by our results in section \ref{sec:aharmonic} we use the A-harmonic extensions of the independent boundary hat functions restricted to the subdomain $\omega$ as elements of a local basis. In Section \ref{sec:oversampledgfem} we find that the relative global error of solution using these as local shape functions in an {\em oversampled-GFEM} scheme developed here while not as small still remains comparable to that of our MS-GFEM.  This can be deduced from the numerical experiments of Section \ref{sec:aharmonic} that show that for a given number of local basis functions the span of this type of local basis is ``nearly the same'' as the span of the $n$-widths used in MS-GFEM. In forthcoming work we illustrate this reason in rigorous mathematical terms for scalar problems in terms of Tchebycheff approximations.

Related recent work is motivated by the theory of low rank approximations based on randomized SVD \cite{HMT11} or rSVD. The work \cite{CLLW19} outlines a reduction of the dimension of the space used to generate $n$-width eigenfunctions numerically by sampling randomly generated basis functions. Another strategy \cite{BS18} provides a method based on rSVD for choosing boundary data and using their A-harmonic extensions as a sub optimal basis and demonstrating a nearly exponential decay in the expected value of the error. In a recent independent development the numerical investigations  \cite{SW17} show that local bases constructed by A-harmonically extending traces of harmonic polynomials deliver exponential decay of the error. 


We  conclude the introduction noting that there are several other strategies to reduce the computational work for multiscale numerical implementations.
One way to address multiscale computation is to exploit the local periodicity or stochasticity of a microstructure and make use of homogenization theory \cite{BLP78}. These methods include the original Multiscale FEM \cite{HW97} and the variational multiscale method \cite{HFM98}. In the absence of local structure one can consider approaches to numerical homogenization for rough coefficients (i.e., $L^\infty$ coefficients); however, the lack of local structure naturally degrades the efficiency and convergence rates. Nevertheless the size of the global solve can be reduced and several methods have been proposed for $L^\infty$ coefficients that offer significant dimension reduction. There is a huge literature and a variety of methods have been introduced. These methods include upscaling based on harmonic coordinates and elliptic inequalities \cite{OZ07}, \cite{OZ07b}, elliptic solvers based on H-matrices \cite{B07,H15}, explicit solution of local computations through Bayesian numerical homogenization \cite{O15},
dimension reduction methods based on global changes of coordinates and MS-FEM for upscaling porous media flows  \cite{EGH06}, \cite{EH06}, the heterogeneous multiscale methods \cite{WE07}, \cite{WMZ05}, \cite{ES08}, and an adaptive coarse scale--fine scale projection method \cite{NPP08}. Additional contemporary methods include numerical homogenization based on the flux norm for $L^{\infty}$ coefficients \cite{BO10}, rough polyharmonic splines \cite{BO18},  subgrid upscaling methods \cite{AB06} and global Galerkin projection schemes for problems with $L^{\infty}$ coefficients and homogeneous Dirichlet boundary data  \cite{M00}. For  a coarse mesh of diameter $H$ local bases that deliver order $H$ convergence with $O((log{(1/H)})^{d+1}$ approximation functions are developed in \cite{MP14}. For comparison the method presented here and in \cite{BL11}, \cite{BH14}, show that the coarse mesh can be fixed arbitrarily and for a given relative error $\tau$ one needs $O((log{(1/\tau)})^{d+1}$ $n$-width approximation functions.


\section{Problem Formulation} \label{sec:problem_formulation}
\subsection{Variational Formulation of the Problem}
To fix ideas we consider the scalar problem over a bounded domain $\Omega \subset \mathbb{R}^2$ with piecewise $C^1$-boundary given by
\begin{equation}\label{eq:globalpde}
   -\div(A(x)\nabla u(x)) = f(x), \text{ } x\in \Omega
\end{equation}
with Neumann boundary conditions prescribed on the boundary $\partial \Omega_N \subset \partial \Omega$
\begin{equation}\label{eq:globalneumann}
   n\cdot A(x)\nabla u(x) = g(x), \text{ }x\in \partial \Omega_N,
\end{equation}
where $n$ is the unit outer normal vector, and Dirichlet boundary conditions on $\partial \Omega_D \subset \partial \Omega$
\begin{equation}\label{eq:globaldirichlet}
   u(x) = h(x) \text{, } x \in \partial \Omega_D,
\end{equation}
such that $\partial \Omega_D \cap \partial \Omega_N = \emptyset$ and $\overline{\partial \Omega_D }\cup \overline{\partial \Omega_N} = \partial \Omega$.  Here $A(x)$ is the $2\times 2$ conductivity matrix with rough coefficients $a_{ij}(x) \in L^\infty(\Omega)$, and satisfies  the standard ellipticity and boundedness conditions
\begin{equation}\label{eq:ellipticity}
   0 \leq \alpha v \cdot v \leq A(x) v\cdot v \leq \beta v \cdot v < \infty, \text{ }\forall x \in \Omega \text{ and } v\in\mathbb{R}^2.
\end{equation}

The unique weak solution of \eqref{eq:globalpde} belongs to the convex space
\begin{equation}\label{eq:convexspace}
H^1_{D}(\Omega) = \{ u \in H^1(\Omega) : u = h \text{ on } \partial \Omega_D \},
\end{equation}
and satisfies
\begin{equation}\label{eq:weakform}
   B\left(u, v \right) = F(v)
\end{equation}
for all $v$ in the energy space
\begin{equation}\label{eq:energyspace}
H^1_{0D}= \{ v \in H^1(\Omega) : v = 0 \text{ on } \partial \Omega_D \}, 
\end{equation}
where
\begin{equation*}
   B(u,v) = \int_\Omega A(x) \nabla u \cdot \nabla v \, dx \text{ and } F(v) = \int_\Omega f v \, dx + \int_{\partial \Omega_N} gv \, ds.
\end{equation*}
The energy norm is given by $\norm{u}_{\mathcal{E}(\Omega)} = (B(u,u))^{1/2}$.



\subsection{MS-GFEM}
\label{Global}
In this paper the numerical solution to \eqref{eq:weakform} is computed using MS-GFEM.  We outline the a priori convergence estimates and motivation for the method.  To illustrate the ideas we consider a computational domain $\Omega \subset \mathbb{R}^2$ containing smooth inclusions separated by a prescribed minimum distance. The MS-GFEM is a version of GFEM \cite{BC94} and is a domain decomposition given by a partition of unity of overlapping subdomains.  Over each subdomain a local approximation space of shape functions is constructed that is adapted to the coefficient of the PDE restricted to the subdomain.  In MS-GFEM the local shape functions are characterized by an optimal local approximation space obtained from  oversampling.  The local computations associated with construction of local approximation spaces are independent and can  be performed in parallel.  The resulting global stiffness matrix can be several orders of magnitude smaller than the stiffness matrix obtained by applying FEM directly (\cite{BL11}).  

We begin by describing the partition of unity.
Let $\{\omega_i\}_{i=1}^{m}$ be a collection of open sets covering the domain $\Omega$ such that $\cup_{i=1}^{N} \omega_i = \Omega$ and let $\phi_i \in C^1(\omega_i)$, $i=1,\ldots,m$, be a partition of unity subordinate to the open covering.  Here the maximum number of sets $\omega_i$ that any point $x$ in $\Omega$ can belong to is at most $\kappa$. The partition of unity functions satisfy the following properties:
\begin{gather*}
   0 \leq \phi_i \leq 1, \,\,\,\,\, i=1,\ldots,N, \label{eq:pu1}\\
   \phi_i(x) = 0, \, x\in \Omega\setminus\omega_i, \,\,\,\,\, i=1,\ldots,N, \label{eq:pu2}\\
   \sum_{i} \phi_{i}(x) = 1 \,\,\,\,\, \forall x\in \Omega, \label{eq:pu3}\\
   \max_{x\in \Omega} \abs{\phi_{i}(x)} \leq C_1, \,\,\,\,\, i=1,\ldots,N, \label{eq:pu4}\\
   \max_{x \in \Omega} \abs{\nabla \phi_{i}(x)} \leq \frac{C_2}{\diam(\omega_i)}, \,\,\,\,\, i=1,\ldots,N, \label{eq:pu5}
\end{gather*}
where $C_1$ and $C_2$ are bounded positive constants and $\diam(\omega_i)$ is the diameter of the set $\omega_i$.  The partition of unity functions are chosen to be flat-topped so that local approximation spaces are linearly independent and for ensuring good conditioning of the global stiffness matrix (\cite{GS07}).   

Now we describe the construction of the optimal local approximation spaces that characterize the MS-GFEM method. These spaces are constructed through oversampling and are spanned by spectral bases given in terms of the eigenspaces of a restriction operator on A-harmonic functions  described below, see \cite{BL11} and \cite{BH14}.
To construct the local bases first let $\{\omega_i^*\}_{i=1}^{N}$ be a second collection of open sets  such that each $\omega_i $ is contained within the larger open set  $\omega_i^*$. For interior domains we require $\dist(\partial\omega^\ast_i,\partial\omega_i)>0$. For $\omega_i$ touching the boundary of $\Omega$, i.e., $\partial\omega_i\cap\partial\Omega\not=\emptyset$, we require $\dist(\partial\omega_i^*\cap \Omega,\partial\omega_i\cap \Omega)>0$.  In what follows we will refer to the subdomains $\omega_i$ and $\omega_i^*$ as patches.  
The local approximations are constructed from a local affine space spanned by a particular solution (see equation \eqref{eq:localparticular} below) and a local approximation space $V^{m_i}_{\omega_i}$.  Here the dimension of the local approximation space is denoted by $m_i$.
On each subdomain $\omega_i$ the local approximation space is created by first constructing a finite dimensional space of functions on $\omega_i^*$ denoted by $V^{m_i}_{\omega_i^*}$ and then restricting them to $\omega_i$ to form the approximation space $V^{m_i}_{\omega_i}$.  In MS-GFEM the local shape functions describing a basis for $V^{m_i}_{\omega_i}$ are  the eigenfunctions associated with the singular values of the restriction operator (see section \ref{sec:localspaces}) acting on A-harmonic functions, i.e., functions $\xi$ that satisfy
\begin{equation}\label{eq:aharmpde}
   -\div(A(x)\nabla \xi(x)) = 0 \text{ on } \omega_i^*.
\end{equation}
The space of A-harmonic functions defined on $\omega_i^*$ is written $H_A(\omega_i^*)$.  For patches that do not share a boundary with $\Omega$ functions are chosen such that they are equivalent up to a constant. The associated quotient space is written $H_A(\omega_i^*)/\mathbb{R}$.  For boundary patches sharing  Neumann conditions on $\partial \omega_i^* \cap \partial \Omega_N$, the local functions  are taken to  satisfy homogeneous Neumann  boundary conditions on $\partial \omega_i^* \cap \partial \Omega_N$. For boundary patches sharing non-homogeneous Dirichlet conditions on $\partial \omega_i^* \cap \partial \Omega_
D$, the local functions are taken to satisfy homogeneous Dirichlet boundary conditions on $\partial \omega_i^* \cap \partial \Omega_D$.  Here all boundary patches either share a Dirichlet boundary or Neumann boundary but not both. The construction of A-harmonic approximation spaces for boundary patches are described in detail in \cite{BL11} and \cite{BH14}. For interior and boundary patches with Neumann data the local approximation space is augmented with the constant functions. In all cases the dimension of the local approximation space over $\omega_i$ is denoted by $m_i$.
The global approximation space is constructed from the local approximation spaces and is defined by
\begin{equation}\label{eq:globalapprox}
   V^N = \left\{ \sum_{i=i}^N \phi_i \xi_i : \xi_i \in V_{\omega_i}^{m_i} \right\},
\end{equation}
and one verifies as in \cite{BB04} that $V^N$ is a subspace of the energy space $H_{0D}^1(\Omega)$.
The boundary data \eqref{eq:globalneumann}, \eqref{eq:globaldirichlet}, and the right hand side of \eqref{eq:globalpde} are satisfied by local particular solutions. These solutions $\chi_i \in H^1(\omega_i^*)$ are defined by
\begin{equation}\label{eq:localparticular}
   -\div(A(x)\nabla \chi_i(x)) = f(x), \text{ } x\in \omega_i^*,
\end{equation}
with Dirichlet  data $\chi_i = 0$ on $\partial \omega_i^*$ for interior patches.  For boundary patches sharing non-homogeneous Neumann data the particular solutions $\chi_i$ also satisfy the boundary data $\eqref{eq:globalneumann}$ and $\chi_i = 0$ on $\partial \omega_i^*\cap \Omega$. For boundary patches sharing non-homogeneous Dirichlet data the particular solutions $\chi_i=\chi_i^R+\chi_i^D$ where $\chi_i^R$ solves \eqref{eq:localparticular} with $\chi_i^R=0$ on $\partial\omega_i^*$ and $\chi_i^D$ satisfy the boundary data $ \eqref{eq:globaldirichlet}$ on $\partial\omega_i^* \cap \partial \Omega$ and homogeneous Neumann data on $\partial\omega^*_i\cap\Omega$ and \eqref{eq:localparticular} with $f=0$.  The global particular solution $u^F$ is then defined by {\em pasting} together the local particular solutions, i.e.,
\begin{equation}\label{eq:globalparticular}
   u^F = \sum_{i=1}^N \phi_i \chi_i.
\end{equation}

The finite dimensional approximate solution of \eqref{eq:weakform} is posed over the convex space $K^N = V^N + u^F$.  The problem becomes a variational inequality over the convex space $K^N$. We seek a solution $u^G \in V^N$ to the following problem for all $v \in V^N$
\begin{equation} \label{eq:galerkin}
   B(u^G, v) = F(v) - B(u^F, v).
\end{equation}
The MS-GFEM approximate solution to \eqref{eq:weakform} is given by $u_0 = u^G+u^F$ and
\begin{equation}\label{argmin}
u_0=\argmin\left\{\Vert u-w\Vert_{{\mathcal{E}}(\Omega)}\,:\,w\in V^N\right\}
\end{equation}
follows from Galerkin orthogonality.
The existence of a unique solution $u_{0} \in K^{N}$ follows from the standard theory of variational inequalities see, e.g., \cite{D89}.  Given any function $u$ in $H^1_{D}(\Omega)$ and trial field $u_T=\sum_{i=1}^N\phi_i\xi_i+u^F$ a nontrivial calculation following the proofs of Theorems 3.2 and 3.3 of \cite{BB04}  shows that if the local error on every $\omega_i$ satisfies
\begin{equation}\label{locerror1}
\Vert u-(\xi_i+\chi_i)\Vert_{{\mathcal E}{(\omega_i})}\leq \epsilon_i\Vert u\Vert_{{\mathcal E}{(\omega_i^*})}\leq\epsilon_i\Vert u\Vert_{\mathcal{E}(\Omega)},
\end{equation}
for some $\xi_i$ in $V_{\omega_i}^{m_i}$ and local particular solution $\chi_i$,
then the global (i.e., total) approximation error of the trial field is bounded by
\begin{equation}\label{globalerror1}
\Vert u-u_T\Vert_{{\mathcal E}(\Omega)}\leq (\kappa)^{1/2}C\left(\sum_{i}^N\epsilon_i^2\right)^{1/2}\Vert u\Vert_{{\mathcal E}{(\Omega})},
\end{equation}
where $C>0$ is independent of $u$ in $H^1_{D}(\Omega)$. Now from Galerkin minimality \eqref{argmin} together with \eqref{globalerror1} we see that the global error of the Galerkin solution $u_0$  is controlled by the local errors, this is the hallmark of the GFEM method. The calculations behind \eqref{locerror1} and \eqref{globalerror1} for MS-GFEM are provided in the appendix for completeness. In the next section we show how to construct optimal finite dimensional local approximation spaces to produce local approximations $\xi_i$ that satisfy \eqref{locerror1} when the local particular solution $\chi_i$ is given. 

\subsection{Optimal Local Approximation Spaces and MS-GFEM}\label{sec:localspaces}
We now describe the optimal local approximation spaces developed in \cite{BL11} and \cite{BH14}.  We begin by restricting attention to a single subdomain $\omega$ and omit the subscript. The oversampling problem is defined over the larger subdomain $\omega^* \supset \omega$. The space of A-harmonic functions $H_A(\omega^*)$ is the space defined as
\begin{equation}\label{aharmonic}
H_A(\omega^*) = \{ v \in H^1(\omega^*) : B(v, w) = 0 \text{ for all } w\in H^1_0(\omega^*)\}.
\end{equation}
Note we can add a constant function to any function in $H_A(\omega^*)$ and still satisfy \eqref{aharmonic}. Hence we consider the quotient space made up of all equivalence classes of functions in $H_A(\omega^*)$ that are the same up to a constant and denote this as
$H_A(\omega^*)/\mathbb{R}$.  Elements $u$ of $H_A(\omega^*)/\mathbb{R}$ restricted to $\omega$ can be approximated using an optimal local spectral basis \cite{BL11} and \cite{BH14}.  The restriction operator $P$ is defined by $Pu(x)=u(x)$ for $x$ in $\omega$. The restriction is a compact map from
$H_A(\omega^*)/\mathbb{R}$ into $H_A(\omega)/\mathbb{R}$  and the optimal local basis on $\omega$ is given by the span of the eigenfunctions $\{\xi_j\}_{i=1}^\infty$ associated with the singular values $\{\lambda_j\}_{j=1}^\infty$ of the restriction operator
\begin{equation}
P^*P\xi_j=\lambda_j\xi_j,
\label{restricteigen}
\end{equation}
see \cite{BL11}. 
This space provides local shape functions and we see that MS-GFEM is a spectral approximation method. The optimal local approximation space is given by the span of the spectral basis and denoted by $V^m_{\omega}={\rm span}\{\xi_1,\xi_2,\ldots,\xi_m\}$.
The optimality of this oversampling space is seen from the theory of Kolmogorov $n$-widths,  \cite{P85}.  The best accuracy of approximation over all $n$-dimensional subspaces $S(n)\subset H_A(\omega)/\mathbb{R}$ of a function $u\in H_A(\omega^*)/\mathbb{R}$  is measured by the Kolmogorov $n$-width as
\begin{equation}\label{eq:nwidth}
   d_n(\omega, \omega^*) = \inf_{S(n)} \sup_{u\in H_A(\omega^*)/\mathbb{R}} \inf_{v\in S(n)} \frac{\norm{Pu - v}_{\mathcal{E}(\omega)}}{\norm{u}_{\mathcal{E}(\omega^*)}},
\end{equation}
and most importantly $d_n(\omega, \omega^\ast)=\sqrt{\lambda_{n+1}}$. Thus the square root of the $n^{th}$ eigenvalue gives the approximation error incurred when approximating the A-harmonic part of the solution using $V_{\omega}^m$.
The optimal subspace achieving the infimum of \eqref{eq:nwidth} for $n=m$ is written as $V_{\omega}^{m}$ and is defined as the span of the first $m$ eigenfunctions satisfying \eqref{restricteigen}, see \cite{P85}.  These eigenfunctions can be computed explicitly and are characterized by the following theorem.

\begin{theorem}[\cite{BL11}, Theorem 3.2]\label{thm:nwidth2}
The optimal approximation space is given by $V^m_{\omega}=\spanset \{ \xi_1, \ldots,$ $ \xi_m \}$, where $\xi_j = P \phi_j$ and $\lambda_j$ and $\phi_j$ are the largest $m$ eigenvalues and corresponding eigenfunctions that satisfy
\begin{equation}\label{eq:fullspectralproblem}
   (\phi_j, \delta)_{\mathcal{E}(\omega)} = \lambda_j (\phi_j, \delta)_{\mathcal{E}(\omega^*)} \text{ } \forall \delta \in H_A(\omega^*).
\end{equation}
The functions $\{\phi_j\}_{j=1}^\infty$ form a complete orthonormal set in $H_A(\omega^*)/\mathbb{R}$.  
\end{theorem}
\begin{definition}\label{spectbasis}
The optimal local basis $\{\xi_i\}_{i=1}^m$ is referred to as the ``$n$-width'' functions or spectral basis and $V^m_\omega={\rm span}\{\xi_1,\xi_2,\ldots,\xi_m\}$ as the spectral subspace.
\end{definition}
An asymptotically exponential decay of the $n$-widths $\sqrt{\lambda_{m+1}}$ for $n=m$ is proved in \cite{BL11}.
\begin{theorem}\label{thm:nwidthdecay}
   The accuracy has nearly exponential decay for $m$ sufficiently large, i.e., for any small $\epsilon$ one has
\begin{equation}\label{eq:nwidth11}
    d_m(\omega, \omega^*) = \sqrt{\lambda_{m+1}} \leq e^{-m^{\frac{1}{1+d}-\epsilon}}.
\end{equation}
  It follows that for $u_{A} \in H_A(\omega^*)$ there exists $\xi \in V_{\omega}^m$ such that
\begin{equation}\label{eq:nwidth2}
   \norm{u_{A} - \xi}_{\mathcal{E}(\omega)} \leq \sqrt{\lambda_{m+1}}\norm{u_{A}}_{\mathcal{E}(\omega^*)}.
\end{equation}
\end{theorem}
Now for a solution $u$ of \eqref{eq:globalpde}  in $H^1_{D}(\Omega)$ we construct a trial $u_T$ with a $\xi_i$ on each $\omega_i^*$ such that \eqref{locerror1} holds.  Here we will  write $\lambda_{m+1}^i$ to correspond to domain $\omega_i^*$. Note first that for any $\chi_i$ and element $\xi_i\in V^{m_i}_{\omega_i}$ we have 
\begin{equation}\label{eq:har}
u-(\xi_i+\chi_i)=(u-\chi_i)-\xi_i,
\end{equation}
where $u-\chi_i$ is A-harmonic. Thus from Theorem \ref{thm:nwidthdecay} we can choose a $\xi_i$ such that
\begin{equation}\label{eq:est}
   \norm{u-\chi_i - \xi_i}_{\mathcal{E}(\omega_i)} \leq \sqrt{\lambda^i_{m+1}}\norm{(u-\chi_i)}_{\mathcal{E}(\omega_i^*)}.
\end{equation}
Now note from \eqref{eq:localparticular} that for interior patches
\begin{equation}\label{eq:patchlocint}
\begin{aligned}
\Vert\chi_i\Vert^2_{\mathcal{E}(\omega_i^*)}&=\int_{\omega_i^*}A(x)\nabla\chi_i\cdot\nabla\chi_i\,dx\\
&=\int_{\omega_i^*}f\chi_i\,dx=\int_{\omega_i^*}A(x)\nabla u\cdot\nabla\chi_i\,dx \leq\Vert u\Vert_{\mathcal{E}(\omega_i^*)}\Vert \chi_i\Vert_{\mathcal{E}(\omega_i^*)},
\end{aligned}
\end{equation}
so
\begin{equation}\label{eq:patchlocintest}
\begin{aligned}
&\Vert\chi_i\Vert_{\mathcal{E}(\omega_i^*)}\leq\Vert u\Vert_{\mathcal{E}(\omega_i^*)},
\end{aligned}
\end{equation}
and a simple calculation shows
\begin{equation}\label{eq:patchlocintesquadt}
\begin{aligned}
&\Vert u-\chi_i\Vert_{\mathcal{E}(\omega_i^*)}\leq 2\Vert u\Vert_{\mathcal{E}(\omega_i^*)}.
\end{aligned}
\end{equation}
From this we get
\begin{equation}\label{est2}
   \norm{u-\chi_i - \xi_i}_{\mathcal{E}(\omega_i)} \leq \epsilon_i\norm{u}_{\mathcal{E}(\omega_i^*)}.
\end{equation}
with $\epsilon_i=2\sqrt{\lambda^i_{m+1}}$.
Proceeding in the same way for boundary patches sharing Neumann data we note $u-\chi_i$ is A-harmonic and choose $\xi_i$ to  obtain 
\begin{equation}\label{est3}
   \norm{u-\chi_i - \xi_i}_{\mathcal{E}(\omega_i)} \leq \epsilon_i\norm{u}_{\mathcal{E}(\omega_i^*)},
\end{equation}
with $\epsilon_i=2\sqrt{\lambda^i_{m+1}}$.
For Dirichlet data similarly $u-\chi_i$ is A-harmonic and we may choose $\xi_i$ to obtain
\begin{equation}\label{est4}
   \norm{u-\chi_i - \xi_i}_{\mathcal{E}(\omega_i)} \leq \epsilon_i\norm{u}_{\mathcal{E}(\omega^*_i)},
\end{equation}
with $\epsilon_i=2\sqrt{5\lambda^i_{m+1}}$.
These calculations are provided in the Appendix for completeness.
The global error estimate \eqref{globalerror1} now follows from \eqref{est2}, \eqref{est3}, and \eqref{est4}.
Therefore it is evident that the a priori estimates showing asymptotic exponential accuracy are obtained when using the optimal local bases  in the construction of the global approximation space $V^N$. The global error on $\Omega$ is controlled by the local errors and we have an a priori asymptotic exponential decay of the MS-GFEM approximation.
The numerical experiments given here show an immediate exponential decay of the $n$-widths for this class of coefficients, see section \ref{sec:localspectral}. The numerical experiments also show that the global approximation error is exponentially decaying independent of contrast, see sections \ref{sec:numerics}, \ref{sec:convgstudy}.

We summarize noting that the $n$-width eigenfunctions are the optimal oversampled approximation functions and are computed for the particular geometry and microstructure $\Omega$.  They carry information about the local fields and allow for a drastic reduction in the dimension of the global stiffness matrix.  The total computational complexity of this method is the same as a direct application of the FEM method and is of the same order as the number of inclusions in the domain $\Omega$, see \cite{BH14}.  But now most all of the work is done in parallel and offline. This results in a global solve given by the inversion of a matrix that is orders of magnitude smaller than the total number of inclusions.


\section{Computational Method}
We outline the implementation of the MS-GFEM algorithm and describe the construction of the global stiffness matrix, particular solutions, and right hand side. Recall that we wish to solve \eqref{eq:weakform} using the optimal local bases pasted together over the computational domain $\Omega \subset \mathbb{R}^2$. As data we are given a matrix of material properties $A(x)$ together with Neumann and Dirichlet boundary conditions specified over the boundary of $\Omega$.  We begin by first creating a standard finite element mesh on the entire domain $\Omega$ which will be used in all computations.  Future implementations will use independent  meshes over each patch $\omega_i$ enhancing the parallel nature of the local computations.

The steps in the MS-GFEM algorithm are:
\begin{enumerate}
   \item Define a covering of $\Omega$ by square patches $\{\omega_i\}_{i=1}^N$ and expand each patch $\{\omega_i^*\}_{i=1}^N$ such that $\omega_i \subset \omega_i^*$ are concentric and the ratio of side edges of $\omega_i$ to $\omega_i^*$ is less than 1.
   \item Construct partition of unity functions $\{\phi_i\}_{i=1}^{N}$ subordinate to $\{\omega_i\}_{i=1}^N$.  In this work we use the standard bilinear and quadratic shape functions over triangles and squares to construct the partition of unity functions.
   \item Solve local particular solutions $\{\chi_i\}_{i=1}^N$ satisfying \eqref{eq:localparticular}.  The boundary data on interior patches is given by $\chi_i = 0$ on $\partial \omega_i^*$.  On boundary patches $\chi_i$ satisfies Neumann boundary data \eqref{eq:globalneumann} or Dirichlet data \eqref{eq:globaldirichlet} on $\partial \omega_i^*\cap \partial \Omega$ with $\chi_i=0$ on $\partial \omega_i^* \cap \Omega$.
   \item Construct right hand side of \eqref{eq:galerkin}, denoted as $\mathbf{r}$.
   \item Construct finite dimensional local spaces $S^{n_i}_{\omega_i^*} \subset H_A(\omega_i^*)$. Here $S^{n_i}_{\omega_i^*}$ is thought of a discrete approximation to $H_A(\omega^*_i)$. The dimension of these spaces are $n_i$ and the elements of this space are  A-harmonic extensions of boundary data given by the piece wise linear hat functions defined on $\partial \omega_i^*$. Together these piece wise linear functions form a partition of unity on $\partial\omega_i^*$. 
   \item Solve the local spectral problems over $S^{n_i}_{\omega_i^*}$ to build the local spectral bases that generate the spaces $V_{\omega_i}^{m_i}$ with $n_i>m_i$.
   \item Construct the global stiffness matrix $\mathrm{G}$ in \eqref{eq:galerkin} with entries given by $B(u, v)$ for $u, v \in V^N$.
   \item Solve the Galerkin formulation given by $\mathrm{G}\mathbf{x}=\mathbf{r}$.
   \item The approximate solution to \eqref{eq:convexspace} is given by $u_0=u^G+u^F$.
\end{enumerate}

We elaborate on the elements of this algorithm below and illustrate the computations for a simple partition of unity comprised of two subdomains in section \ref{sec:numerics}.  For details in computing finite dimensional subspaces of $H^0_A(\omega_i^*)$ and the generation local spectral bases see section \ref{sec:aharmonic}.  

\subsection{Construction of Local Solution Spaces $S^{n_i}_{\omega_i^*}$.}
The material domain $\Omega$ is covered by square or hexahedral patches $\omega_1, \ldots, \omega_N$ with an expanded covering $\omega_1^*, \ldots, \omega_N^*$ such that $\omega_i \subset \omega_i^*$.
We use the following notations.  Let $\mathcal{H}=\{ E_e \}_{e=1}^{n_{el}}$ be a finite element mesh of $\Omega$, with elements $E_e$, which conforms to the boundaries of each patch $\omega_i$ and $\omega_i^*$, $i=1,\dots, N$.  The piecewise  bilinear and quadratic shape functions are denoted $H_n(\alpha)$ where ${\alpha}=(\alpha_1, \alpha_2)$ are coordinates on the reference element and ${x}=(x_1,x_2)$ are the Cartesian coordinates on the domain $\Omega$.  The map $\alpha_e: E_e \rightarrow \square$ is the change of coordinates map from the $e$-th element $E_e$ to the reference element $\square$.  We begin by generating a discrete A-harmonic solution space $S^{n_i}_{\omega_i^*} \subset H_A(\omega_i^*)$ of A-harmonic FEM solutions of
\begin{equation}\label{eq:aharm}
\left\{
\begin{aligned}
   -\div (A(x)\nabla w_{i}^k(x)) = 0 \text{ for }x\in\omega_i^* \\
   w_{i}^k(x) = h_i^k(x) \text{ for }x \in \partial \omega_i^*.
\end{aligned}
\right.
\end{equation}
Here $h_i^k(x)$ is the $k$-th hat function defined by $h_i^k = 1$ at the $k$-th boundary node on $\partial \omega_i^* \cap \Omega$ and $h_i^k = 0$ at all other boundary nodes. Here $1\leq k\leq n_i$.  
For boundary patches, $\partial\omega_i^* \cap \partial \Omega \neq \emptyset$, only hat functions corresponding to interior nodes, $\partial\omega_i^*\cap\Omega$, are used as boundary data, and $h_i^k(x)$ satisfies homogeneous Neumann data on $\partial\omega_i^*\cap\partial\Omega_N$ and homogeneous Dirichlet data on $\partial\omega_i^*\cap\partial\Omega_D$.
We write the $k$-th A-harmonic function over $\omega_i^*$ as
\begin{equation} \label{eq:aharmfn}
   {w}_{i}^k({x}) = \sum_{e=1}^{n_{el}}\sum_{l=1}^{n_{en}}\mathbf{w}_{i,e}^{k,l} H_{l}({\alpha}_e \left({x})\right),
\end{equation}
where $\mathbf{w}_{i,e}^{k}$ is the vector of nodal values at the $e$-th element, $1 \leq e \leq n_{el}$, contained in $\omega_i^*$, with $\mathbf{w}_{i,e}^{k,l}$ the component at the $l$-th node, $1\leq l \leq n_{en}$.

\subsection{Construction of Local Spectral Bases, $V^{m_i}_{\omega_i}$}
Now we construct the spectral basis that delivers the local approximation with optimal approximation properties. We introduce the span $V_{\omega_{i}^{*}}^{m_{i}}=\spanset\left\{\xi_{i}^{1},\dots\xi_{i}^{m_{i}}\right\}$, where the $\xi_{i}^{j}$ are the eigenfunctions corresponding to the $m_{i}$ largest eigenvalues of
\begin{equation} \label{eq:spectralMatrices}
\mathrm{Q}_{i}\mathbf{x}=\lambda \mathrm{P}_{i}\mathbf{x}.
\end{equation}
Here the entries of the matrices $\mathrm{Q}_{i}$ and $\mathrm{P}_{i}$ are given as
\begin{align}\label{eq:Q_ijk}
\mathrm{Q}_{i}^{jk} &= (w_{i}^{j}, w_{i}^{k})_{\mathcal{E}(\omega_i)} = \int_{\omega_{i}} A\nabla w_{i}^{j}\nabla w_{i}^{k} dx = \int_{\partial\omega_{i}} ({n}\cdot A \nabla w_{i}^{j})  w_{i}^{k} ds, \\[0.5ex] \label{eq:P_ijk}
\mathrm{P}_{i}^{jk} &= (w_{i}^{j}, w_{i}^{k})_{\mathcal{E}(\omega_i^*)} = \int_{\omega_{i}^{*}} A\nabla w_{i}^{j}\nabla w_{i}^{k} dx = \int_{\partial\omega_{i}^{*}} ({n}\cdot A \nabla w_{i}^{j})  w_{i}^{k} ds,
\end{align}
where we have used $S^{n_i}_{\omega_i^*}$ as an approximate basis for $H_A(\omega_i^*)$ in the local spectral problems defined by \eqref{eq:fullspectralproblem}.  These functions give the discrete approximation to the $n$-width functions when $\xi^j_i$'s are restricted to $\omega_i$ and the space $V^{m_i}_{\omega_i}=\spanset\left\{\xi_{i}^{1},\dots\xi_{i}^{m_{i}}\right\}$ gives the discrete approximation to the spectral subspace on $\omega_i$. The $k$-th entry of the vector $\mathbf{x}$ in \eqref{eq:spectralMatrices} is defined to be the coefficient of the $k$-th basis function $w_i^k(x) \in S^{n_i}_{\omega_i^*}$.  Thus, the $q$-th, $1 \leq q \leq m_i$, eigenfunction over $\omega_i^*$ is written as
\begin{equation}\label{eq:nwidthfn}
   {\xi}_i^q({x}) = \sum_{k=1}^{m_i} \xi_i^{q,k} {w}_{i}^{k}({x}),
\end{equation}
with coefficients $\mathbf{x}_k = \xi_i^{q,k}$.  
Here the eigenfunctions are listed as $\{\xi_i^1, \xi_i^2, \ldots, \xi_i^{m_i}\}$, where the corresponding eigenvalues are listed in descending order $1 > \lambda_i^1 \geq \lambda_i^2 \geq \cdots \lambda_i^{m_i} > 0$ ($\mathrm{P}_i$ and $\mathrm{Q}_i$ are symmetric positive definite).

Combining \eqref{eq:aharmfn} and \eqref{eq:nwidthfn} gives the $n$-width functions written in terms of the underlying finite element discretization
\begin{equation*}\label{eq:nwidthdiscret}
   {\xi}_i^q ({x})= \sum_{e=1}^{n_{el}}\sum_{k=1}^{m_i} \sum_{l=1}^{n_{en}} \xi_i^{q,k}\mathbf{w}_{i,e}^{k,l} H_l({\alpha}_e ({x})).
\end{equation*}

The entries of the spectral matrices \eqref{eq:Q_ijk} may be computed quickly by rewriting the middle integrals in \eqref{eq:Q_ijk} as
\begin{align*}
    Q_i^{jk} &= \int_\omega A \nabla w^j \cdot \nabla(v_\omega + \hat{w}^k) \, ds \\
    P_i^{jk} &= \int_{\omega^*} A \nabla w^j \cdot \nabla(v_{\omega^*} + \tilde{w}^k) \, ds
\end{align*}
where $w^k=v_\omega + \hat{w}^k$ on $\omega$ with $v_\omega \in H^1_0(\omega)$ vanishing on $\partial \omega$ and $w^k=v_{\omega^*}+\tilde{w}^k$ with $v_{\omega^*} \in H^1_0(\omega^*)$ vanishing on $\partial \omega^*$. In the finite element discretization $\hat{w}^k$ is zero on all nodes in the interior of $\omega$ and $\hat{w}^k=w^k$ on nodes lying on $\partial \omega$. Similarly $\tilde{w}^k$ is zero on all nodes interior to $\omega^*$ and $\tilde{w}^k=w^k$ on nodes lying on $\partial\omega^*.$ Then the spectral matrices may be computed by
\begin{align*}
    Q_i^{jk} &= \int_\omega A \nabla w^j \cdot \nabla(\hat{w}^k) \, ds \\
    P_i^{jk} &= \int_{\omega^*} A \nabla w^j \cdot \nabla(\tilde{w}^k) \, ds
\end{align*}
where the integrand is nonzero only on elements away from the boundary. This follows since $H_A^0 \perp H^1_0$ in the $\norm{\cdot}_{\mathcal{E}}$-norm by definition. This integral reduces the number of elements to be integrated over from quadratic to linear in two dimensions and cubic to quadratic in three dimensions. The reduction in memory usage is the same, allowing for a much larger FEM mesh to be used underneath these computations.

\subsection{Global Stiffness Matrix and Right Hand Side}
Using the bases for the local trial fields $V^{m_i}_{\omega_i}$ the global trial field $V^N$ is defined as in \eqref{eq:globalapprox} by 
\begin{equation*}\label{eq:globalnwidths}
   V^N = \spanset\{ \phi_i \xi_i^q : 1 \leq i \leq N,\text{ } 1\leq q \leq m_i, \text{ and } \xi_i^q \in V^{m_i}_{\omega_i} \}.
\end{equation*}
The local spaces $V^{m_i}_{\omega_i}$ are augmented by the constant functions for interior patches or patches with Neumann boundary conditions. 
For boundary patches on which Dirichlet boundary conditions are enforced, $\partial \omega_i \cap \partial \Omega_D \neq \emptyset$, constant functions need not be added to $V^{m_i}_{\omega_i}$.
The Galerkin discretization defined by \eqref{eq:galerkin} gives the matrix equation $\mathrm{G}\mathbf{x}=\mathbf{r}$ where the stiffness matrix $\mathrm{G}$ is defined entrywise by
\begin{equation}\label{eq:globstiffmat}
   \mathrm{G}_{iqjr} = (\phi_i { \xi}_i ^{q}, \phi_j { \xi}_j^{r})_{\mathcal{E}(\Omega)}.
\end{equation}
The $ij$-th block of $G$ corresponds to cross products of functions from $\phi_i V^{m_i}_{\omega_i}$ and $\phi_j V^{m_j}_{\omega_j}$.  The $iq$-th element of the unknown vector $\mathbf{x}$ is then the coefficient of the global trial function $\phi_i \xi_i^q$, so that we have the part of the solution corresponding to the local A-harmonic functions written
\begin{equation}\label{eq:aharmpart}
   u^G = \sum_{i=1}^N \sum_{q=1}^{m_i} \mathbf{x}_{i}^q \phi_i \xi_i^q.
\end{equation}
Likewise, the right hand side of \eqref{eq:galerkin} is written entrywise as
\begin{equation}\label{eq:rhs}
   \mathbf{r}_{iq} =  F(\phi_i \xi_i^q) - B(u^F, \phi_i \xi_i^q) = \int_{\Omega} f \phi_i \xi_i^q \, dx + \int_{\partial \Omega_N} g \phi_i \xi_i^q \, dS - \int_\Omega u^F \phi_i \xi_i^q \, dx.
\end{equation}

With $\mathrm{G}\mathbf{x}=\mathbf{r}$ defined as above, solving for $\mathbf{x}$ gives the A-harmonic part of the solution \eqref{eq:aharmpart}.
The approximate solution to \eqref{eq:globalpde} is given by
\begin{equation*}
   u_0 = u^G + u^F = \sum_{i=1}^N\sum_{q=1}^{m_i} \mathbf{x}_i^q \phi_i \xi_i^q + \sum_{i=1}^N \phi_i \chi_i
\end{equation*}
for local particular solutions $\chi_i$ as in \eqref{eq:localparticular}.


\section{Numerical Implementation for a Benchmark Problem}\label{sec:numerics}

For the implementation of the MS-GFEM algorithm multiple partial differential equations must be solved. It should be noted that any method to solve PDEs may be used, here we use the finite element method. The problem domain $\Omega$ has an underlying finite element mesh and overlapping patches $\omega_i$ to which the mesh conforms. 
All computations are performed over the finite elements. A concrete numerical example of the MS-GFEM algorithm is demonstrated below.

\subsection{Problem Formulation} \label{subsec:problem_formulation}

\begin{figure}[htbp] 
\centering
\includegraphics[height=5cm]{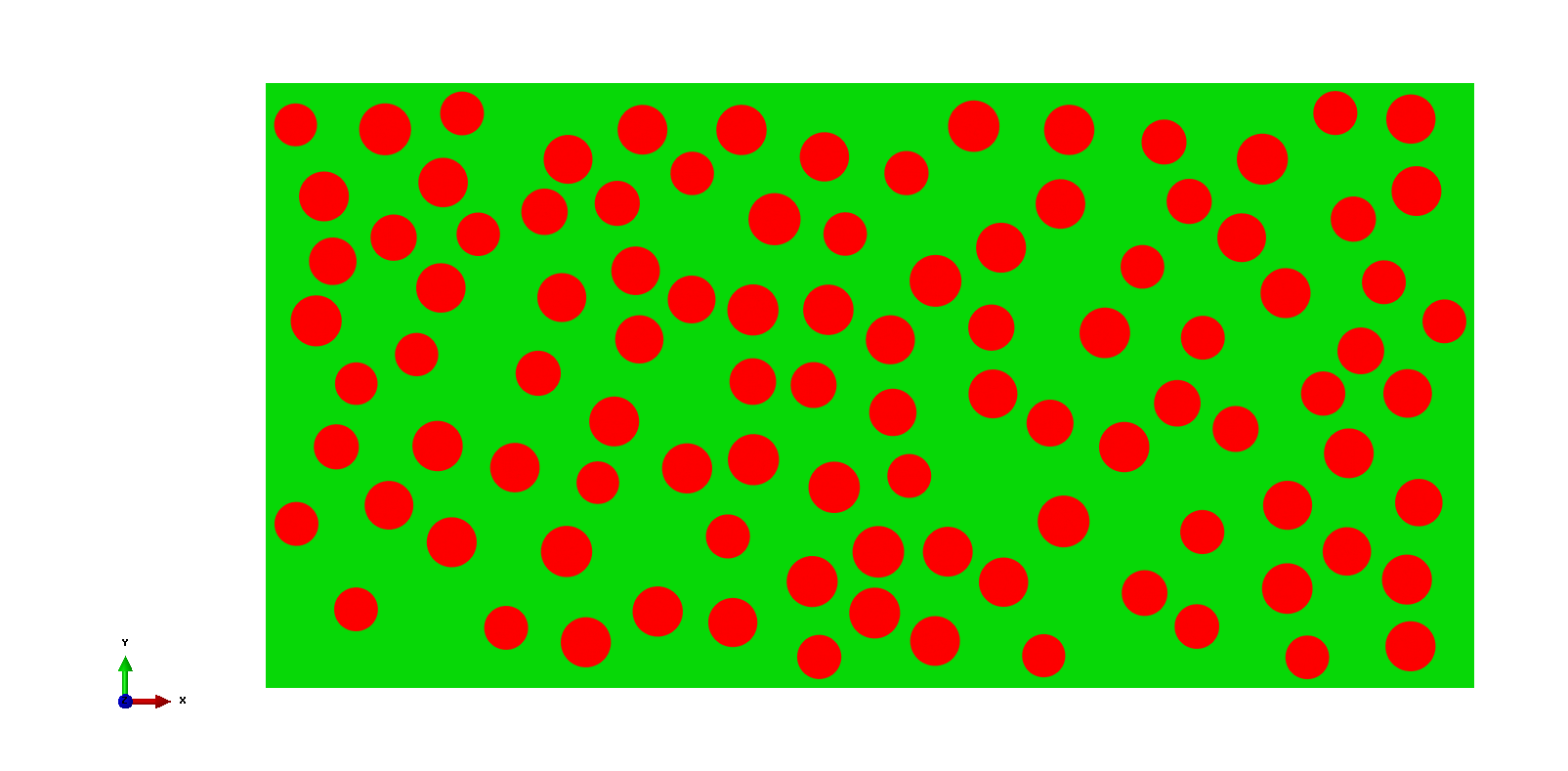} 
\caption{Rectangular domain of size $20\times10$ with 100 circular inclusions for problem \eqref{eq:pde_2}.}
\label{fig:100_inclusions_01}
\end{figure}

All computations in this work are done on the same rectangular domain with 100 separated circular inclusions with small variations in radii (figure~\ref{fig:100_inclusions_01}). For demonstration of the numerical implementation we consider the following problem.
\begin{equation}\label{eq:pde_2}
\left\{
\begin{aligned}
   -\div(A(x_1,x_2)\nabla u(x_1,x_2)) &= 0, \text{ } (x_1,x_2) \in \Omega \\
   u(x_1,x_2)& = 0 \text{, } (x_1,x_2) \text{ on the left of } \partial\Omega  \\
   u(x_1,x_2)& = 1 \text{, } (x_1,x_2) \text{ on the right of } \partial \Omega \\
   {n}\cdot A\nabla u(x_1,x_2) &= 0, \text{ } (x_1,x_2) \text{ on the top and bottom of }\partial\Omega.
\end{aligned}
\right.
\end{equation}
Dirichlet boundary conditions are applied on the left and right of the domain while the top and bottom have homogeneous Neumann boundary conditions.  
The matrix material has property $A_{1}$ and the inclusion material has property $A_{2}$. They are given as
\begin{equation*}\label{eq:A}
A_{1} = 
\begin{pmatrix*}
\alpha & 0 \\
0 & \alpha   
\end{pmatrix*}, \,\,\, 
A_{2} = 
\begin{pmatrix*}
\beta & 0 \\
0 & \beta  
\end{pmatrix*}.
\end{equation*}
In the section \ref{sec:convgstudy} we will study a large range of material contrasts given by the ratio $\alpha/\beta$.  For the purpose of illustration we shall take $\alpha=1$ and $\beta=100$ in this section.

\subsection{Mesh Selection} \label{subsec:mesh refinement}
\begin{figure}
    \centering
    \includegraphics[width=4in]{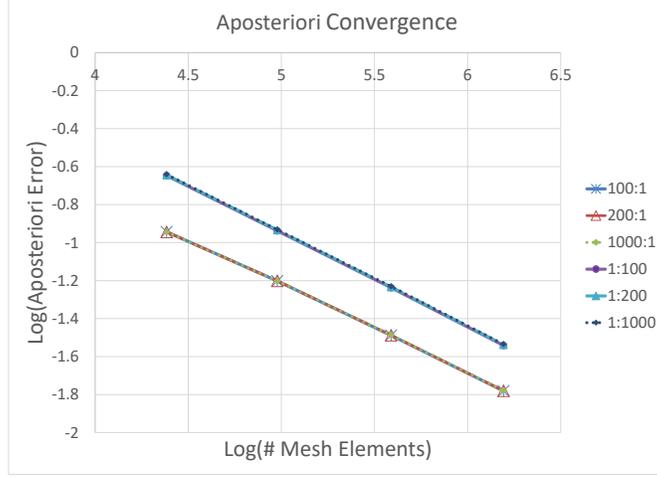}
    \caption{Log plot of a posteriori error estimate against the log of the number of mesh elements The green curve is the a posteriori contrast for the matrix conductivity 1000 and the inclusion 1, while the blue curve is for matrix conductivity 1 and inclusion conductivity 1000. For other contrasts the convergence graphs are essentially the same.}
    \label{fig:aposteriori}
\end{figure}

All computations in this work are performed on the same rectangular domain. Here we will compare the accuracy of the MS-GFEM implementation to the standard FEM method.  In order to compare the methods we select the underlying mesh used for the MS-GFEM method by solving for the standard FEM solution for the same problem over progressively finer meshes for the highest contrast used in our study. The FEM error for a given mesh is measured using the a posteriori estimate
\begin{equation}\label{eq:aposteriori}
    \epsilon_h^2 = \frac{\norm{u_{h/2}}^2 - \norm{u_h}^2}{\norm{u}^2}\frac{1}{1-2^{-2\gamma}},
\end{equation}
where we use bilinear elements and assume $\epsilon_h^2=Ch^{2\gamma}$ for $1/2\leq \gamma\leq 1$ (see  \cite{BH14}). Here the a posteriori estimate \eqref{eq:aposteriori} is standard and we have chosen $\gamma=1/2$. The finest mesh which gives an acceptable FEM solution over all material contrasts will be the same one used for comparison of the MS-GFEM method with the FEM method. Here we will carry out the study for contrasts ranging from $1/1000$ to $1000$. From Table \ref{tab:aposteriori1} we choose the finest mesh with $1,560,058$ elements as the common mesh in which to compare MS-GFEM solutions. The a posteriori estimate for the FEM using finer and finer mesh seeds is shown in figure~\ref{fig:aposteriori}. For the purposes of this study we denote the FEM solution with the finest mesh as the {\em overkill} solution for a given contrast and denote this solution by $u$. The MS-GFEM solution will now be computed over this same mesh. The following subsections describe the MS-GFEM implementation.

\begin{table}[!ht]
\caption{Log of a posteriori errors computed according to \eqref{eq:aposteriori}.}
\centering
\begin{tabular}{| c | c | c | c | c | c | c | c |}
\cline{3-8}
\multicolumn{2}{ c }{\rule{0pt}{3ex}}  &  \multicolumn{6}{|p{60ex}|}{\centering A posteriori Error by Contrast (fiber:matrix)}  \\ \hline
\rule{0pt}{3ex}\# Elements & Log(\# Elements)    
& { 100:1}
& { 200:1}
& { 1000:1} 
& { 1:100}  
& { 1:200} 
& { 1:1000} \\ \hline
\rule[-.3\baselineskip]{-3pt}{3ex}
6686    & 3.825 & --        & --       & --       & --       & --       & --\\
24187   & 4.384 & -0.94312  & -0.94368 & -0.94418 & -0.64733 & -0.64384 & -0.64103 \\
95017   & 4.978 & -1.20170  & -1.20145 & -1.20129 & -0.93778 & -0.93424 & -0.93140 \\
388189  & 5.589 & -1.48816  & -1.48767 & -1.48730 & -1.23698 & -1.23345 & -1.23061 \\ 
1560058 & 6.193 & -1.78071  & -1.78001 & -1.77948 & -1.54125 & -1.53772 & -1.53489 \\\hline
\end{tabular}
\label{tab:aposteriori1}
\end{table}

\subsection{Construction of Partition of Unity}
\begin{figure}[ht] 
\centering
\begin{tabular}{cc}
\includegraphics[height=2.5cm]{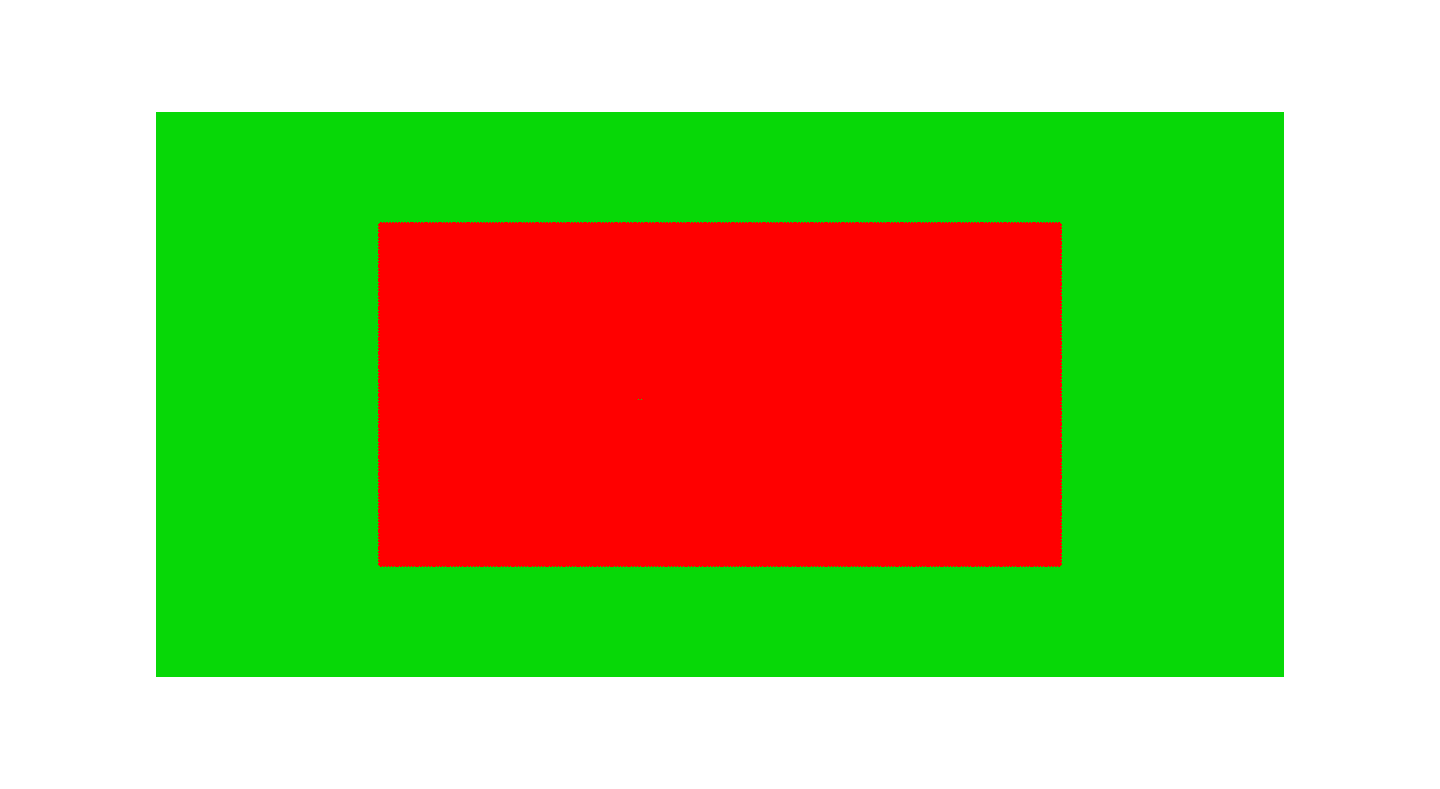} & \includegraphics[height=2.5cm]{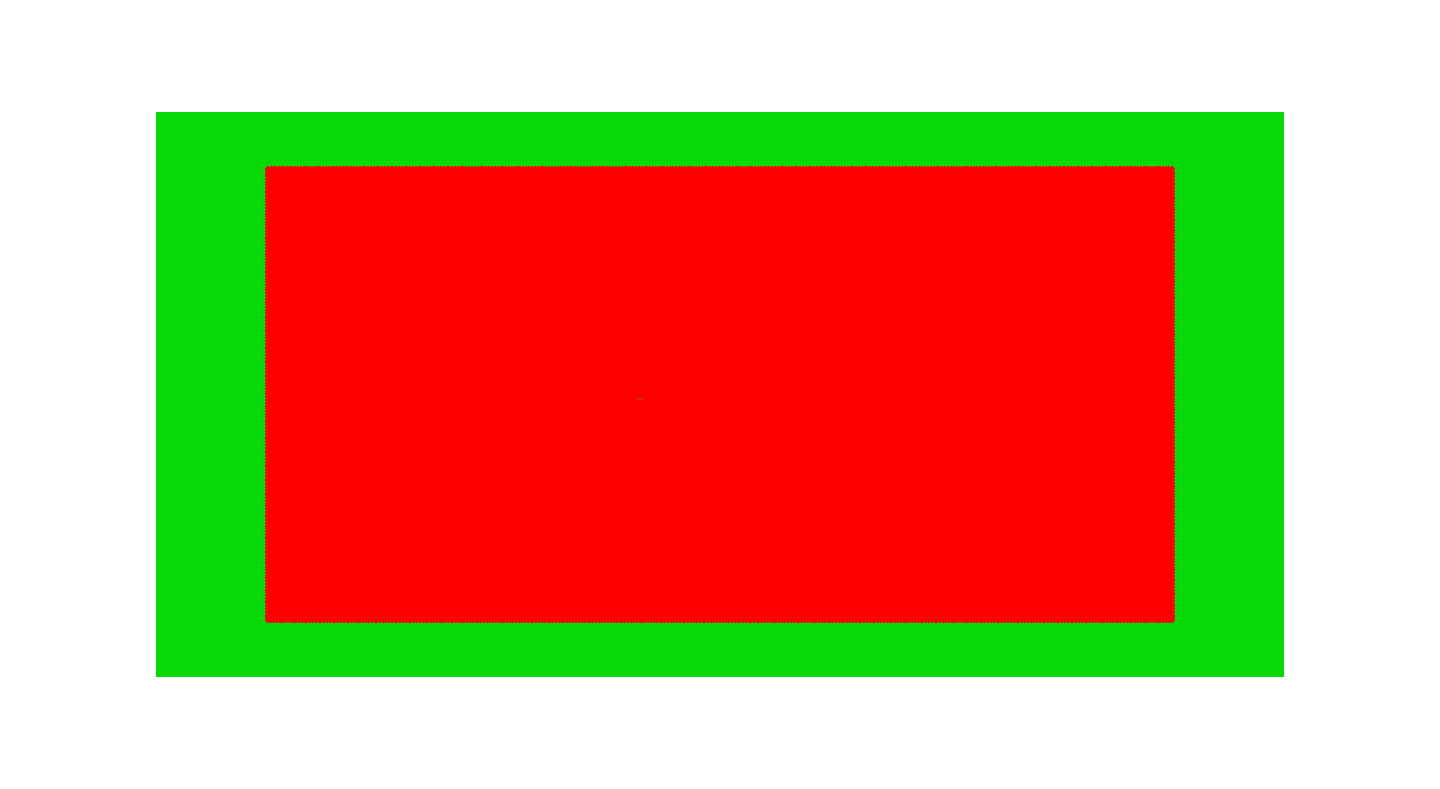} \\ \includegraphics[height=2.5cm]{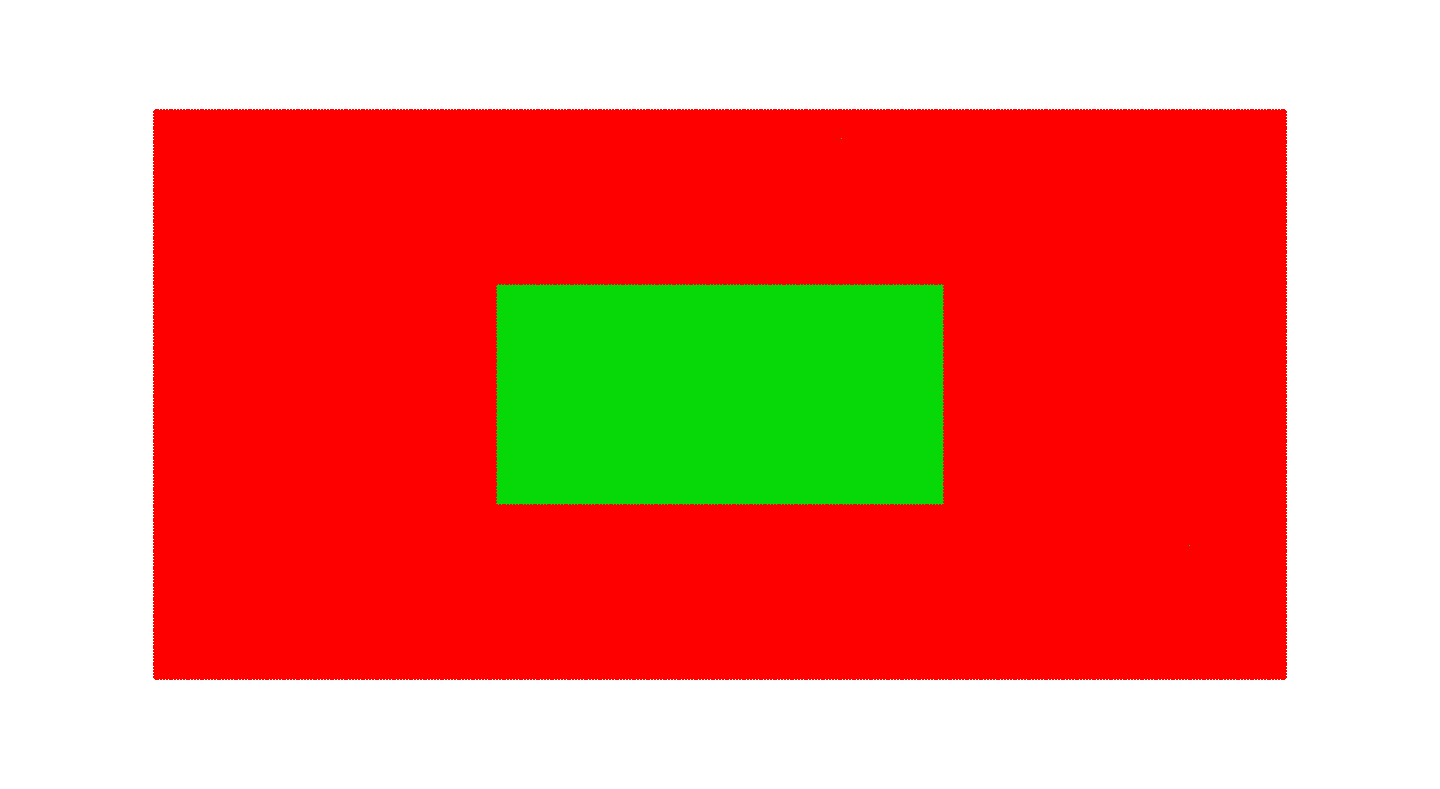} & \includegraphics[height=2.5cm]{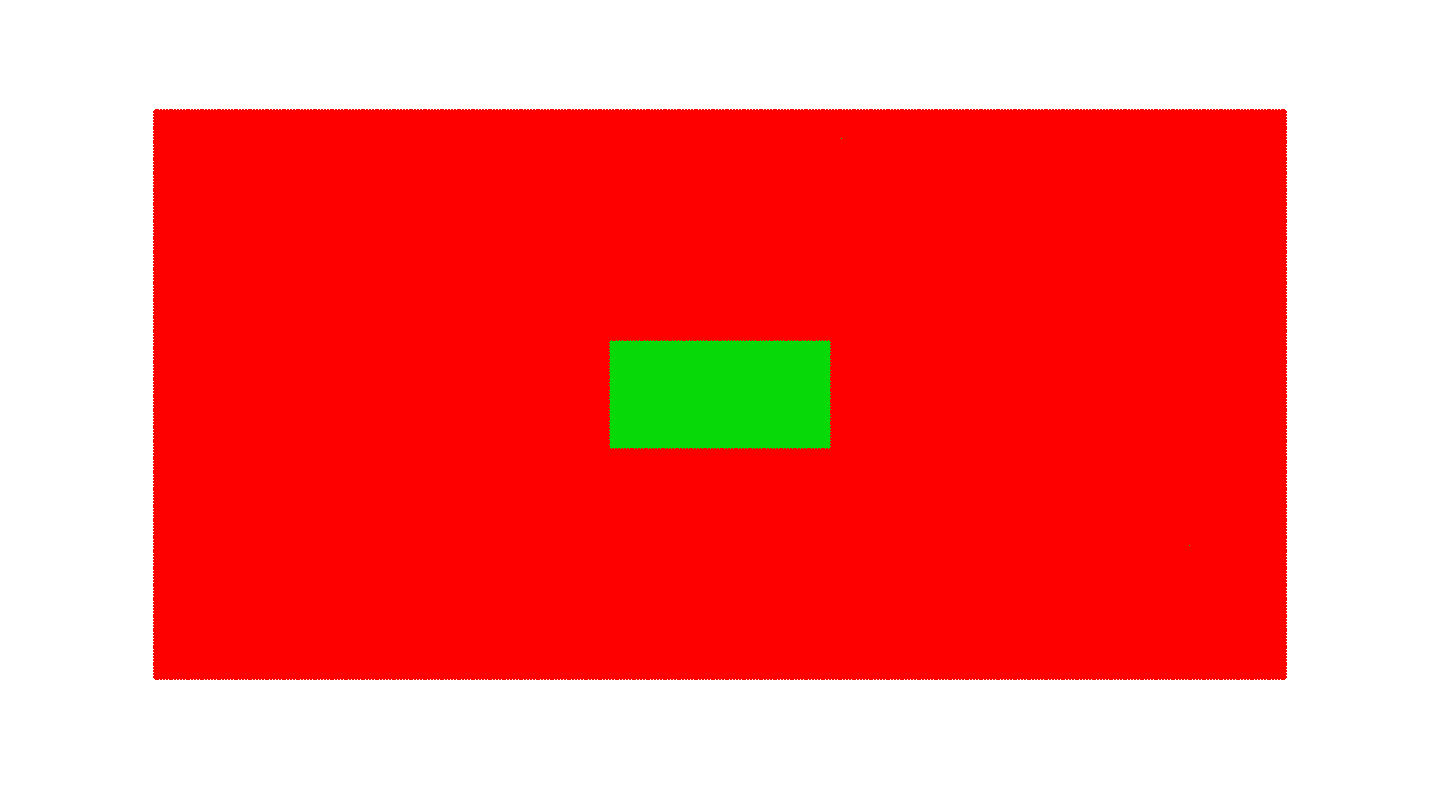}
\end{tabular}
\caption{Top row: The patches $\omega_{1}$, $\omega_{1}^{*}$, respectively. Bottom row: $\omega_{2}$, $\omega_{2}^{*}$, respectively.  Red indicates the subdomain, green is the area of $\Omega$ outside the subdomain.}
\label{fig:covering}
\end{figure}

The domain is a $20\times10$ rectangle with 100 circular inclusions with small variation in radii. To fix ideas we work with a two patch covering $\{\omega_1, \omega_2\}$ of $\Omega$ by a rectangular patch $\omega_1 = 12\times6$ and a rectangular annulus $\omega_2 = 20\times 10 - 8\times4$. The expanded patches are defined as $\omega_1^* = 16\times8$ and $\omega_2^* = 20\times10 - 4\times2$. The covering is shown in Figure \ref{fig:covering}.

\begin{figure}[ht] 
\centering
\begin{tabular}{cc}
\includegraphics[height=2.5cm]{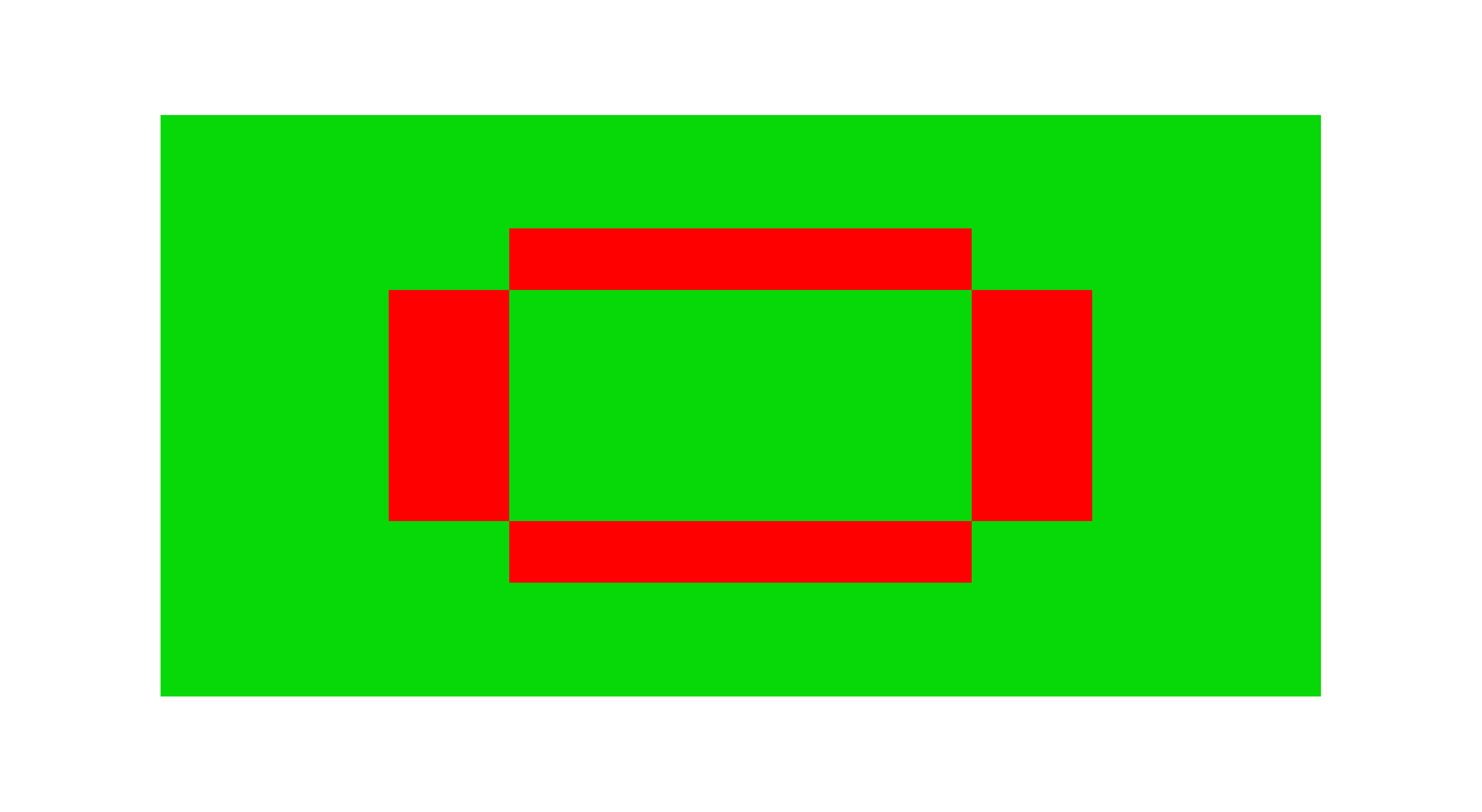} & \includegraphics[height=2.5cm]{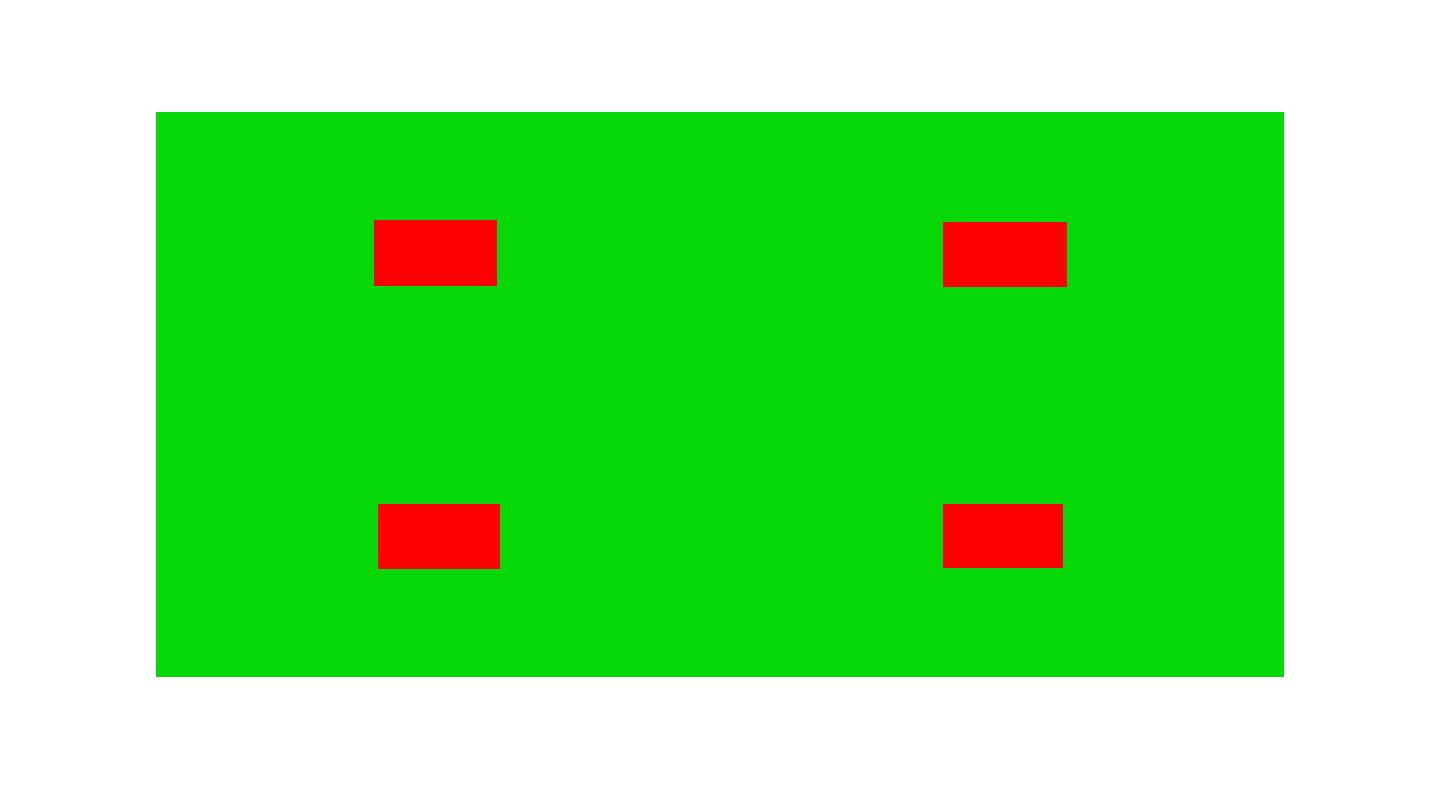}
\end{tabular}
\caption{Sets in the overlap of $\omega_{1}\cap\omega_{2}$ shown in red.  On the left region $\phi_1$ is linear.  On the right $\phi_1$ is bilinear.}
\label{fig:pou_sets}
\end{figure}

\begin{figure}[ht] 
\centering
\begin{tabular}{cc}
\includegraphics[height=2.75cm]{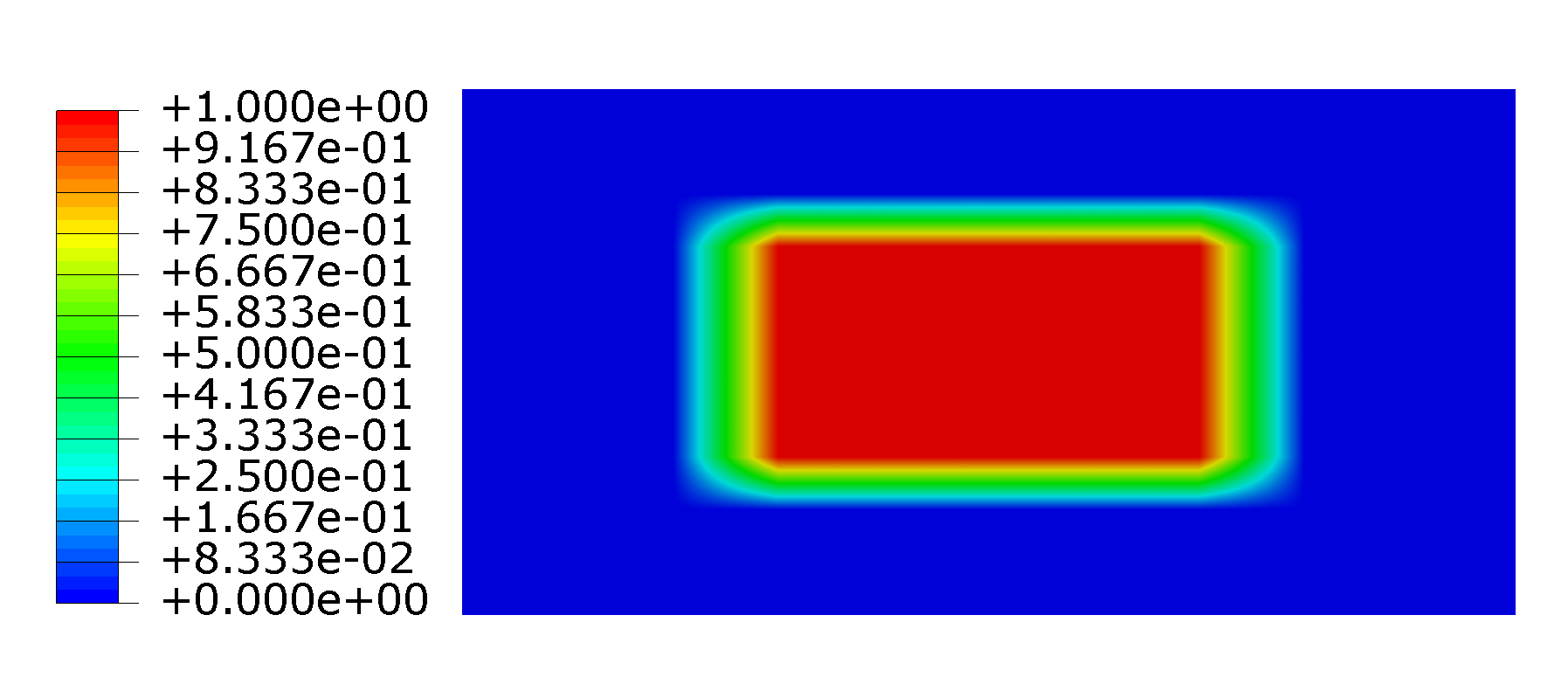} &
\includegraphics[height=2.75cm]{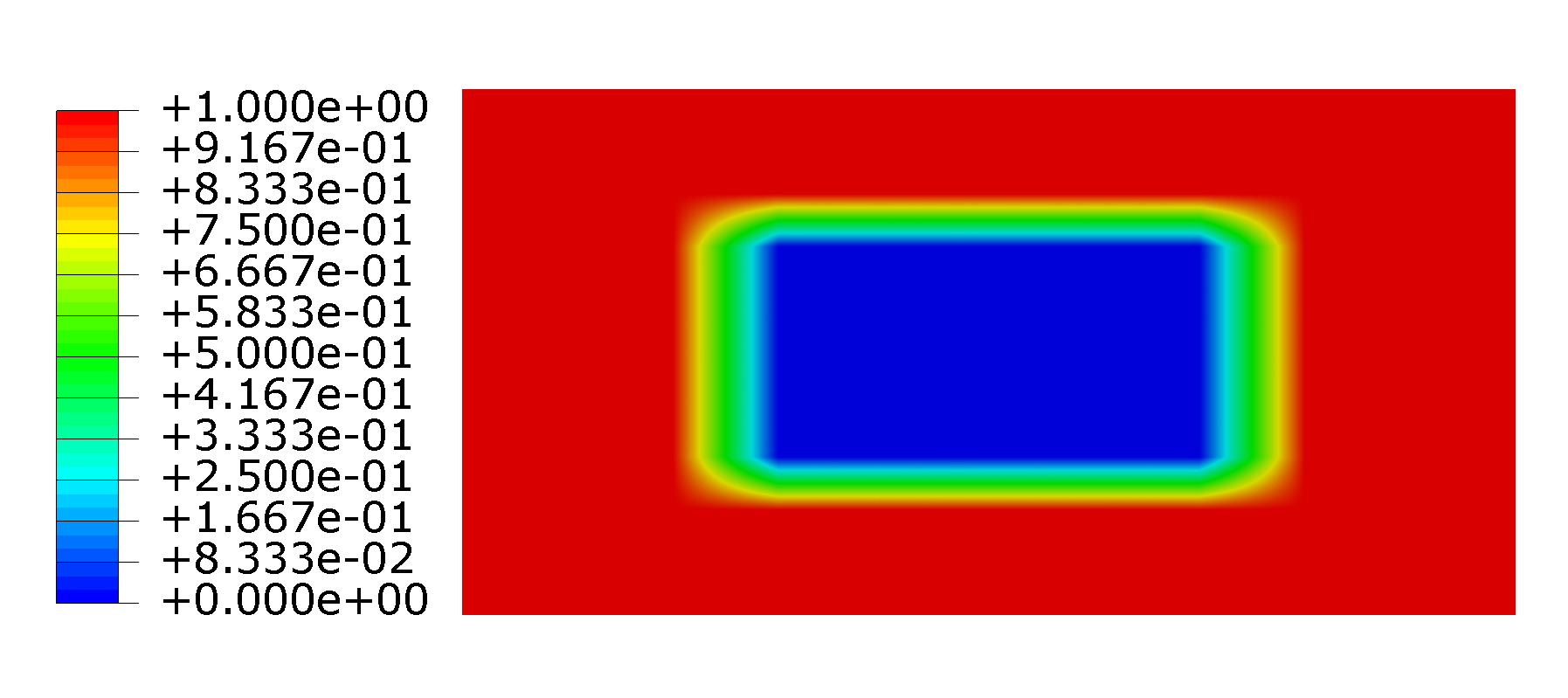}\\
\includegraphics[height=2.75cm]{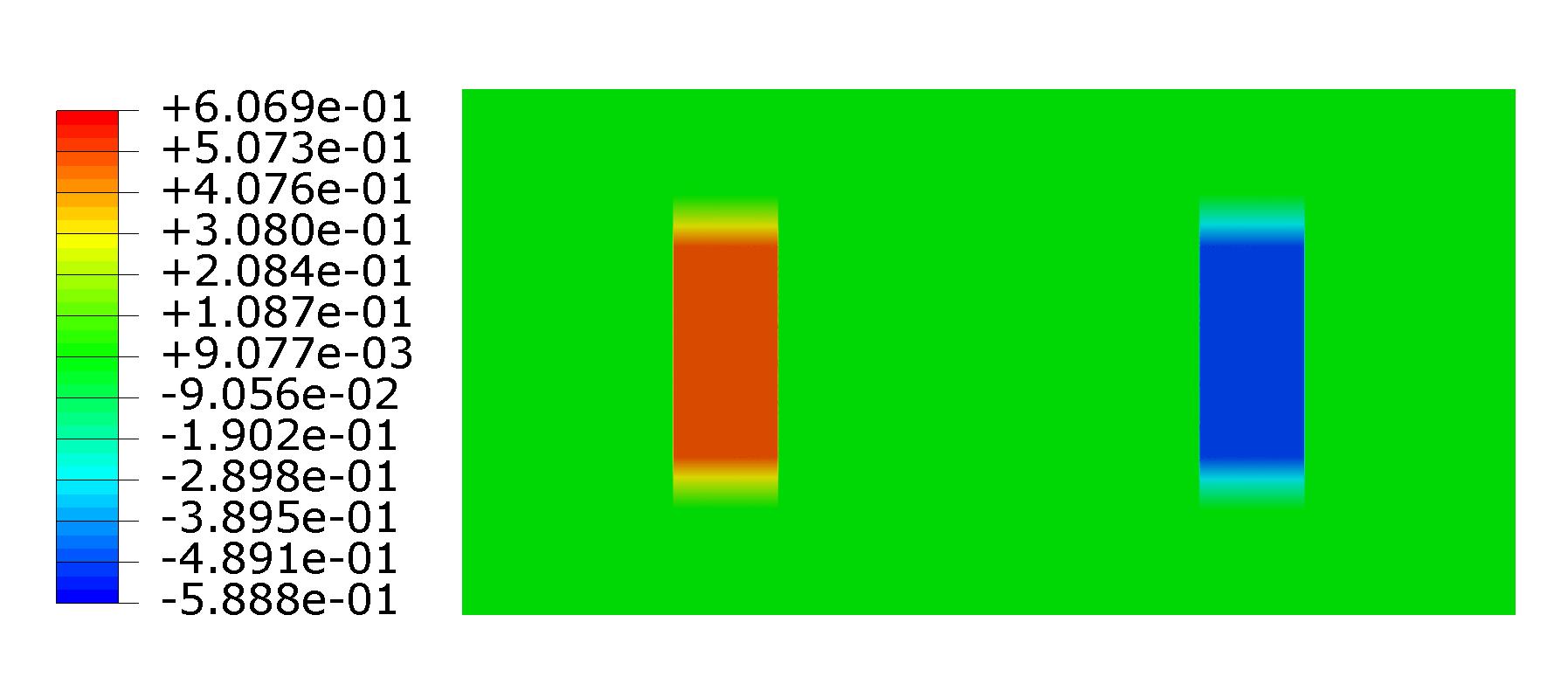} & \includegraphics[height=2.75cm]{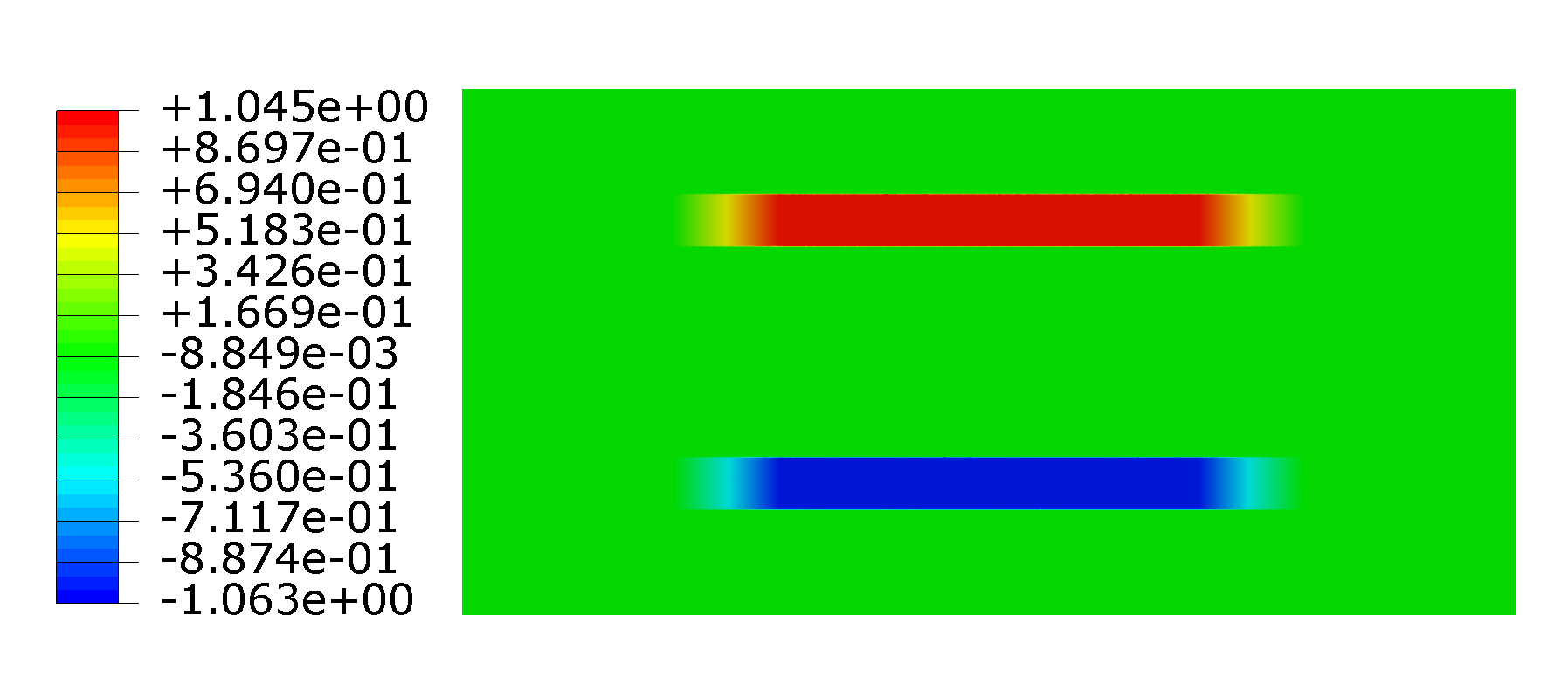}
\end{tabular}
\caption{The partition of unity functions $\phi_1$ (top left), $\phi_2$ (top right), and the derivatives with respect to $x_1$ (bottom left) and $x_2$ (bottom right) over $\Omega$.}
\label{fig:pou_1}
\end{figure}

To construct the partition of unity function $\phi_{1}$ over $\omega_{1}$, the overlap between $\omega_{1}$ and $\omega_{2}$ is determined and split into two sets as shown in Figure \ref{fig:pou_sets}. On $\omega_1 - \omega_2$ (the $8\times4$ green rectangle on the bottom left of Figure~\ref{fig:covering}) the partition of unity function is 1, outside $\omega_1$ it is 0. On the overlap shown on the left of Figure~\ref{fig:pou_sets}, $\phi_1$ is linear.  On the corners of the overlap, shown on the right, it is bilinear. The values of $\phi_1$ are computed at the nodes and after that interpolated to the integration points using the finite element shape functions. The partition of unity function $\phi_{2}$ is then just $\phi_{2}=1-\phi_{1}$. The partition of unity functions and their derivatives are shown in Figure \ref{fig:pou_1}.

\subsection{Local Particular Solutions}

\begin{figure}[ht] 
\centering
\begin{tabular}{c}
\includegraphics[height=3.5cm]{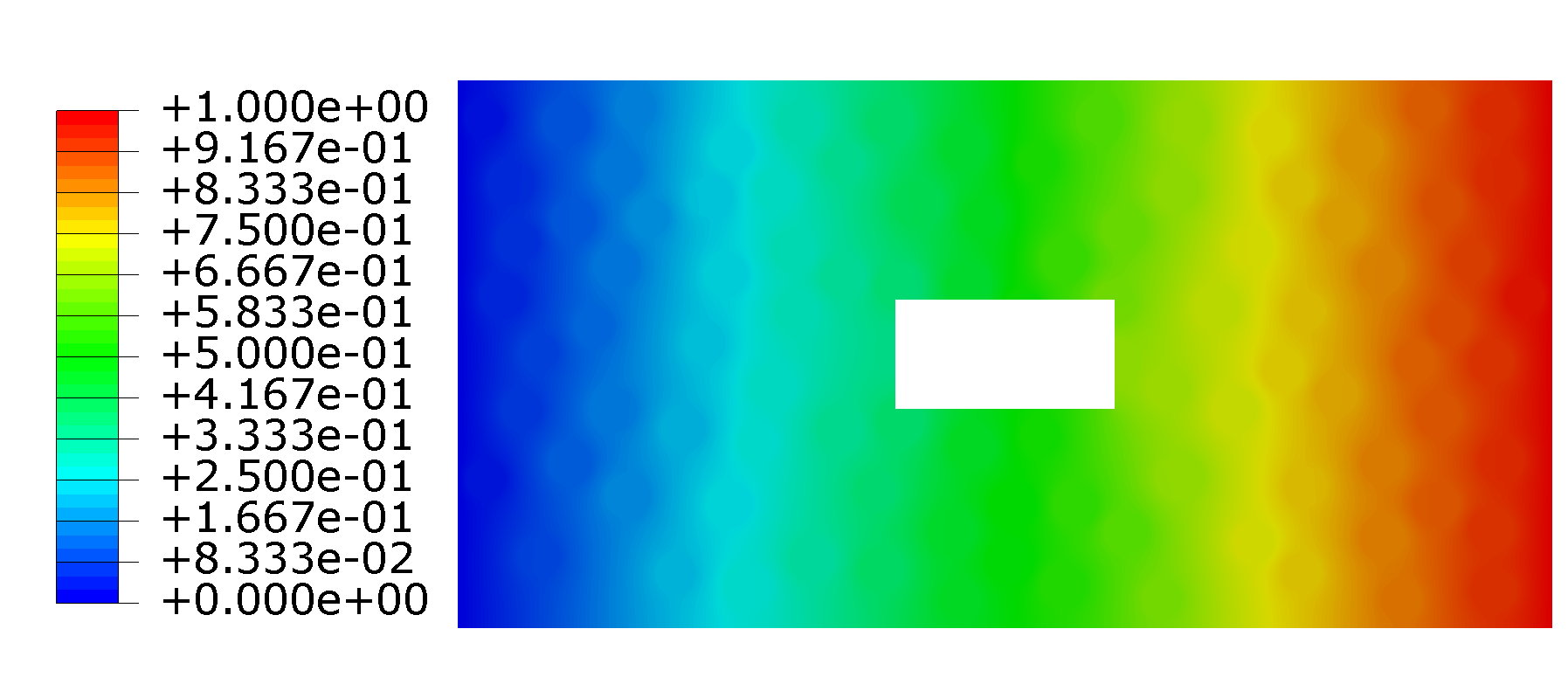}
\end{tabular}
\caption{The particular solution $\chi_{2}$ displayed on $\omega_2^*$. Here $\chi_{1}\equiv0$ since $\omega_1^*$ is an interior domain ($\partial\Omega \cap \partial\omega_1^*=\emptyset$) and the PDE has 0 right hand side.}
\label{fig:partSol}
\end{figure}

We find the local particular solutions $\chi_{1}$ and $\chi_{2}$ over $\omega_1^*$ and $\omega_2^*$, respectively. On boundaries of the patches $\partial \omega_i^*$ which do not coincide with the outer boundary $\partial \Omega$ we enforce $\chi_i=0$ and on boundaries which coincide with outer boundaries the given boundary conditions in \eqref{eq:pde_2} are applied. Since the PDE is homogeneous and the boundary data is homogeneous Dirichlet data, $\chi_1\equiv0$. 
For $\omega_2^*$ the problem becomes
\begin{equation}\label{eq:partSol2}
\left\{
\begin{aligned}
-\div(A\nabla \chi_2(x_1,x_2)) &= 0 \text{, } (x_1,x_2) \in \omega_{2}^{*} \\
-n\cdot A\nabla\chi_2(x_1,x_2) &= 0 \text{, } (x_1,x_2) \in \partial \omega_{2}^{*}\cap\Omega \\
\chi_2(x_1,x_2) &= 0 \text{, } (x_1,x_2) \text{ on the left side of } \partial \omega_2^* \cap \partial\Omega_{D} \\
\chi_2(x_1,x_2) &= 1 \text{, } (x_1,x_2) \text{ on the right side of } \partial \omega_2^* \cap \partial\Omega_{D} \\
-{n}\cdot A\nabla \chi_2(x_1,x_2) &= 0 \text{, } (x_1,x_2) \in \partial\omega_2^*\cap\partial \Omega_N.
\end{aligned}
\right.
\end{equation}
The particular solution $\chi_2$ is shown in Figure \ref{fig:partSol}. The finite element analysis for $\chi_2$ is used in the assembly of the right hand side, see \eqref{eq:rhs}.

\subsection{Construction of Local Solution Spaces $S^{n_i}_{\omega_i^*}$}

\begin{figure}[ht] 
\centering
\begin{tabular}{cc}
\includegraphics[height=3.5cm]{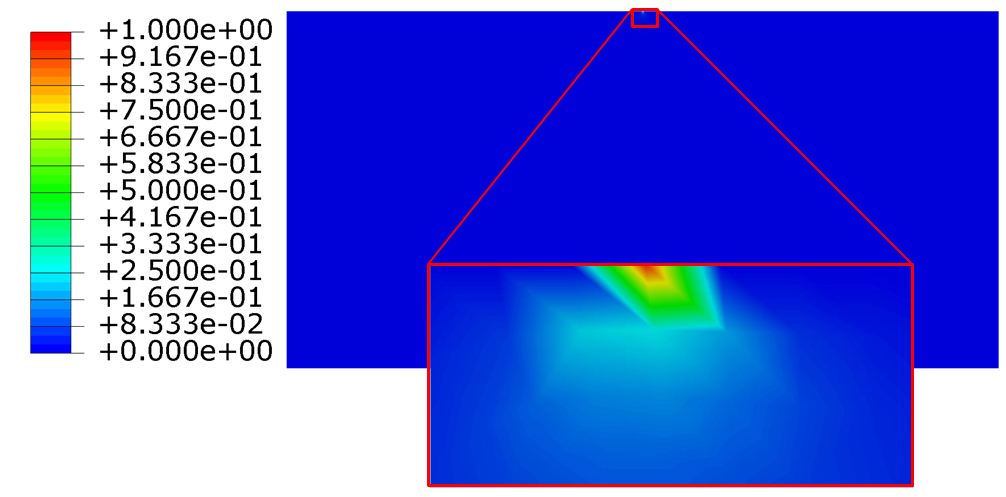} &
\includegraphics[height=3.6cm]{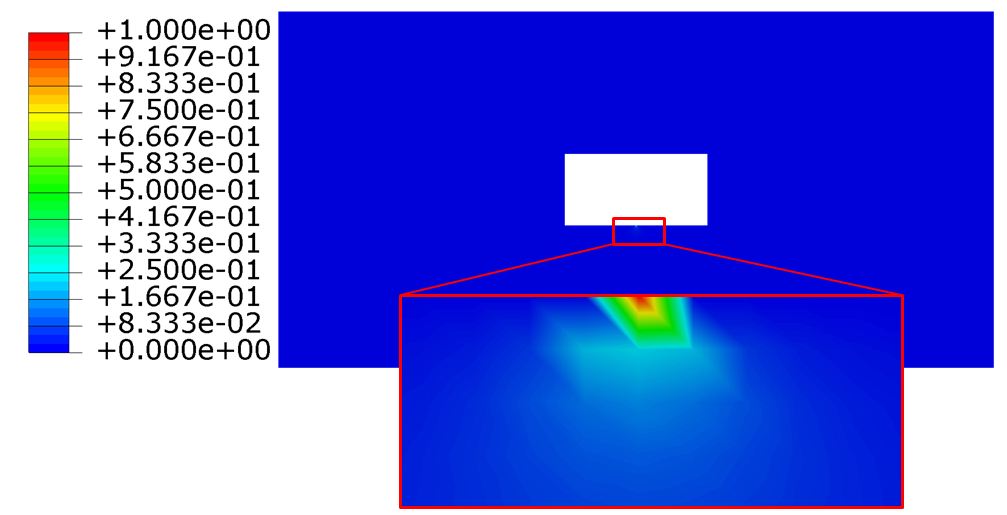}
\end{tabular}
\caption{Two of the computed functions $w_1^j(x,y)$ (left) and $w_2^j(x,y)$ (right) shown over $\omega_1^*$ and $\omega_2^*$, respectively, using a hat function with a single node as support for boundary data. We have re-scaled the region near the boundary datum to make the field visible.}
\label{fig:hat1}
\end{figure}

\begin{figure}[ht] 
\centering
\begin{tabular}{cc}
\includegraphics[height=3.5cm]{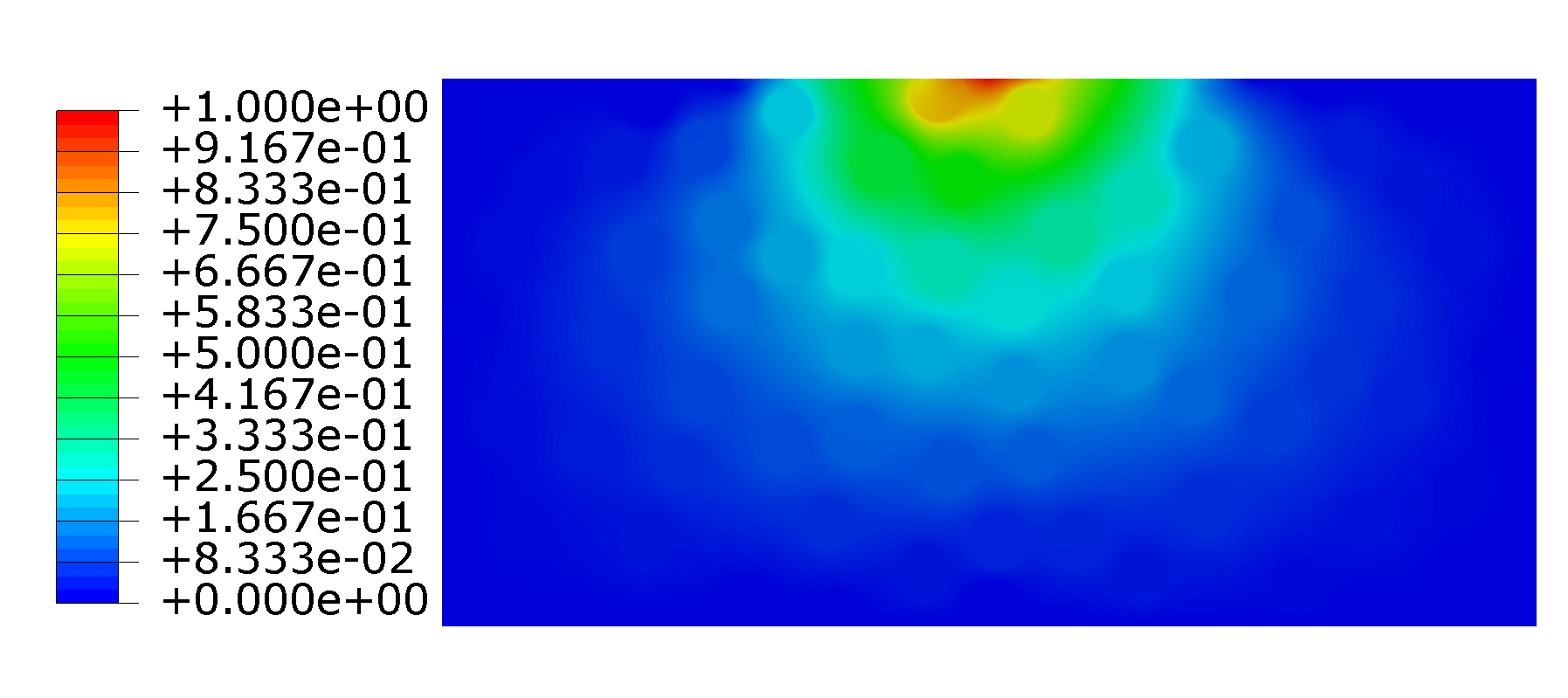} &
\includegraphics[height=3.5cm]{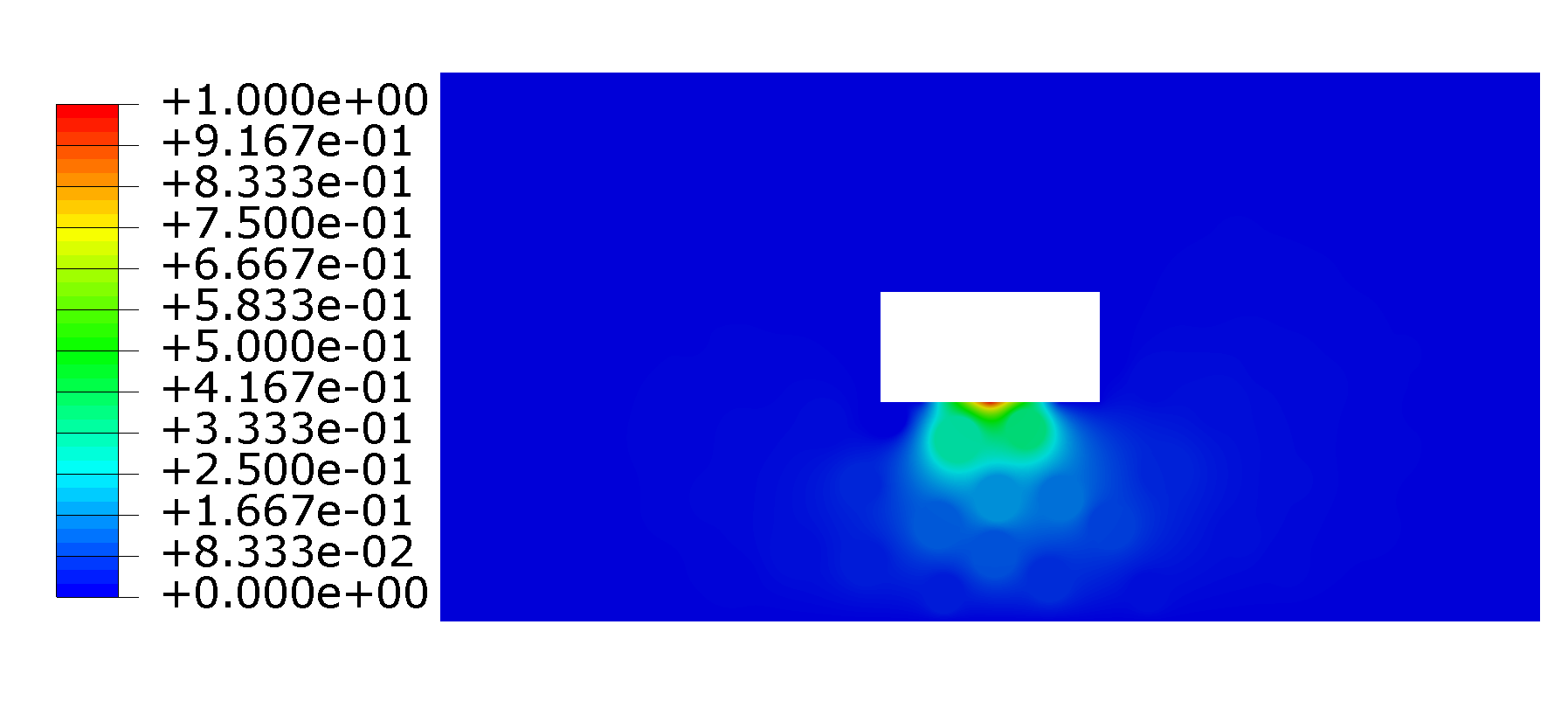}
\end{tabular}
\caption{An example of A-harmonic functions $w_i^j$ on $\omega_1^*$ (left) and $\omega_2^*$ (right) using a hat function with 301 nodes (left) and 101 nodes (right) as support for boundary data.} 
\label{fig:hat5}
\end{figure}

Next we generate the local spaces $S^{n_1}_{\omega_1^*}$ and  $S^{n_2}_{\omega_2^*}$. Here we use as boundary data the piece wise linear hat functions defined on $\partial \omega_i^*$. Taken together these piece wise linear functions form a partition of unity on $\partial\omega_i^*$. 
To make an element of $S^{n_1}_{\omega^*_1}$ we consider each of the $1\leq j \leq n_1$ nodes on $\partial\omega_1^*$. For each node we solve the discrete Dirichlet problem
\begin{equation}\label{eq:solSpaces1}
\left\{
\begin{aligned}
-\div(A\nabla w_{1}^{j}(x,y)) &= 0, \text{ } (x,y) \in \omega_{1}^{*} \\
w_{1}^{j}(x,y) &= 1 \text{, } \text{on the } j \text{-th node of } \partial\omega_1^* \\
w_{1}^{j}(x,y) &= 0 \text{, } \text{every other node on }\partial\omega_1^*.
\end{aligned}
\right.
\end{equation}
the collection of A-harmonic functions constructed this way gives the basis for $S^{n_1}_{\omega^*_1}$.
On $\omega_2^*$, $S^{n_2}_{\omega_2^*}$ consists of functions corresponding to hat functions on $\partial\omega_2^*\cap\Omega$ as boundary data.  Homogeneous Dirichlet and Neumann conditions are enforced on $\partial\omega_2^*\cap\partial\Omega_D$ and $\partial\omega_2^*\cap\partial\Omega_N$, respectively.  For each of the $1\leq j \leq n_2$ nodes on $\partial\omega_2^*\cap\Omega$ the following problem is solved.
\begin{equation}\label{eq:solSpaces2}
\left\{
\begin{aligned}
-\div(A\nabla w_{2}^{j}(x,y)) &= 0, \text{ } (x,y) \in \omega_{2}^{*} \\
w_{2}^{j}(x,y) &= 1 \text{, } \text{on the } j \text{-th node of } \partial\omega_2^*\cap\Omega \\
w_{2}^{j}(x,y) &= 0 \text{, } \text{every other node on } \partial\omega_2^*\cap\Omega \\
w_{2}^{j}(x,y) &= 0 \text{, } (x,y) \in \partial\omega_2^* \cap \partial \Omega_{D} \\
-{n}\cdot A\nabla w_{2}^{j}(x,y) &= 0, \text{ } (x,y) \in \partial\omega_2^* \cap \partial \Omega_N.
\end{aligned}
\right.
\end{equation}
Figure \ref{fig:hat1} shows two computed functions $w_i^j(x,y)$ over $\omega_1^*$ and $\omega_2^*$.
The figure shows that the functions $w_i^j(x,y)$ are localized and decay rapidly away from the boundary when the hat function used as boundary data has support of only one node.  We note that if we consider broader hat functions with support over several boundary nodes instead of just one node then the associated solutions of \eqref{eq:solSpaces1} and \eqref{eq:solSpaces2} have energy densities that extend and are nonzero further inside $\omega_1^*$ and $\omega_2^*$. To see this we provide a computation using a hat function with support over 301 nodes on the boundary of $\omega_1^*$ and 101 nodes on the boundary of $\omega_2^*$ instead of one node, results are shown in Figure~\ref{fig:hat5}.  In section \ref{sec:aharmonic} we investigate the use of hat function boundary data with different support to generate the local solution spaces. In this approach we pick a collection of hat functions all of the same width and taken together form a partition of unity on the boundary. It is clear that using fewer boundary datum reduces the number of problems to solve in generating $S^{n_i}_{\omega^*_i}$. However the interesting observation found in section \ref{sec:aharmonic} is that the reduction of the number of boundary datum used in generating the space done in this way does not significantly effect the convergence rate of the spectral basis. This is because A-harmonic extensions of boundary data with large support on the boundary have energy that decays slowly into the domain relative to A-harmonic extensions of boundary data with small support sets, see Figures \ref{fig:hat1} and \ref{fig:hat5} and \cite{BLS08} and \cite{LS16}. We say that solutions with slow energy decay {\it penetrate} into the domain.

In conclusion the equations \eqref{eq:solSpaces1} and \eqref{eq:solSpaces2} will be solved using bilinear or quadratic FEM, the functions $w_{i}^{j}$ are saved at the nodes of the underlying finite element mesh and represented by finite element shape functions and the stiffness matrices $P_i$ and $Q_i$ are assembled.

\subsection{Construction of Local Spectral Bases}
\label{sec:localspectral}

Now we use $S_{\omega_i^*}^{n_i}$ to construct the local spectral basis with span given by the space $V_{\omega_{i}^{*}}^{m_{i}}={\rm span}\left\{\xi_{i}^{1},\dots\xi_{i}^{m_{i}}\right\}$, where the $\xi_{i}^{j}$ are the eigenfunctions corresponding to the $m_{i}$ largest eigenvalues of
\begin{equation} \label{eq:Q_ixlP_ix}
\mathrm{Q}_{i}\mathbf{x}=\lambda \mathrm{P}_{i}\mathbf{x}
\end{equation}
with entries of $\mathrm{P}$ and $\mathrm{Q}$ defined as in \eqref{eq:Q_ijk} and \eqref{eq:P_ijk}.  For this study the integration is performed over the domain rather than the boundary using the underlying finite elements, and, due to symmetry, only for the entries $j\geq k$.  Equation \eqref{eq:Q_ijk} is written as
\begin{equation}\label{eq:I_i}
\int_{\omega_{i}} A\nabla w_{i}^{j}\nabla w_{i}^{k} \, dx=\sum_{e=1}^{N^{e}_{\omega_{i}}}{(\mathbf{w}_{i,e}^{j})^T} \mathrm{K}_{e}\mathbf{w}_{i,e}^{k},
\end{equation}
where $T$ denotes the transpose, $N^{e}_{\omega_{i}}$ is the number of elements in patch $\omega_{i}$, $\mathbf{w}_{i,e}^{j}$ and $\mathbf{w}_{i,e}^{k}$ are the vectors of nodal values of the functions $w_{i}^{j}$ and $w_{i}^{k}$ at the element $e$, and $\mathrm{K}_{e}$ is the element stiffness matrix.
The integral (\ref{eq:P_ijk}) over $\omega_i^*$ is computed by adding (\ref{eq:I_i}) to the integral over the set $\omega_{i}^{*}-\omega_{i}$. The summation is done in parallel using OpenMP. To ensure that the matrix $\mathrm{P}_{i}$ is positive definite, each row $j$ of $\mathrm{P}_i$ and $\mathrm{Q}_i$ is normalized by $\mathrm{P}_{i}^{jj}$.

\begin{figure}[ht] 
\centering
\begin{tabular}{cc}
\includegraphics[trim={1.5cm 2.25cm 1.5cm 1.5cm},clip,height=5cm]{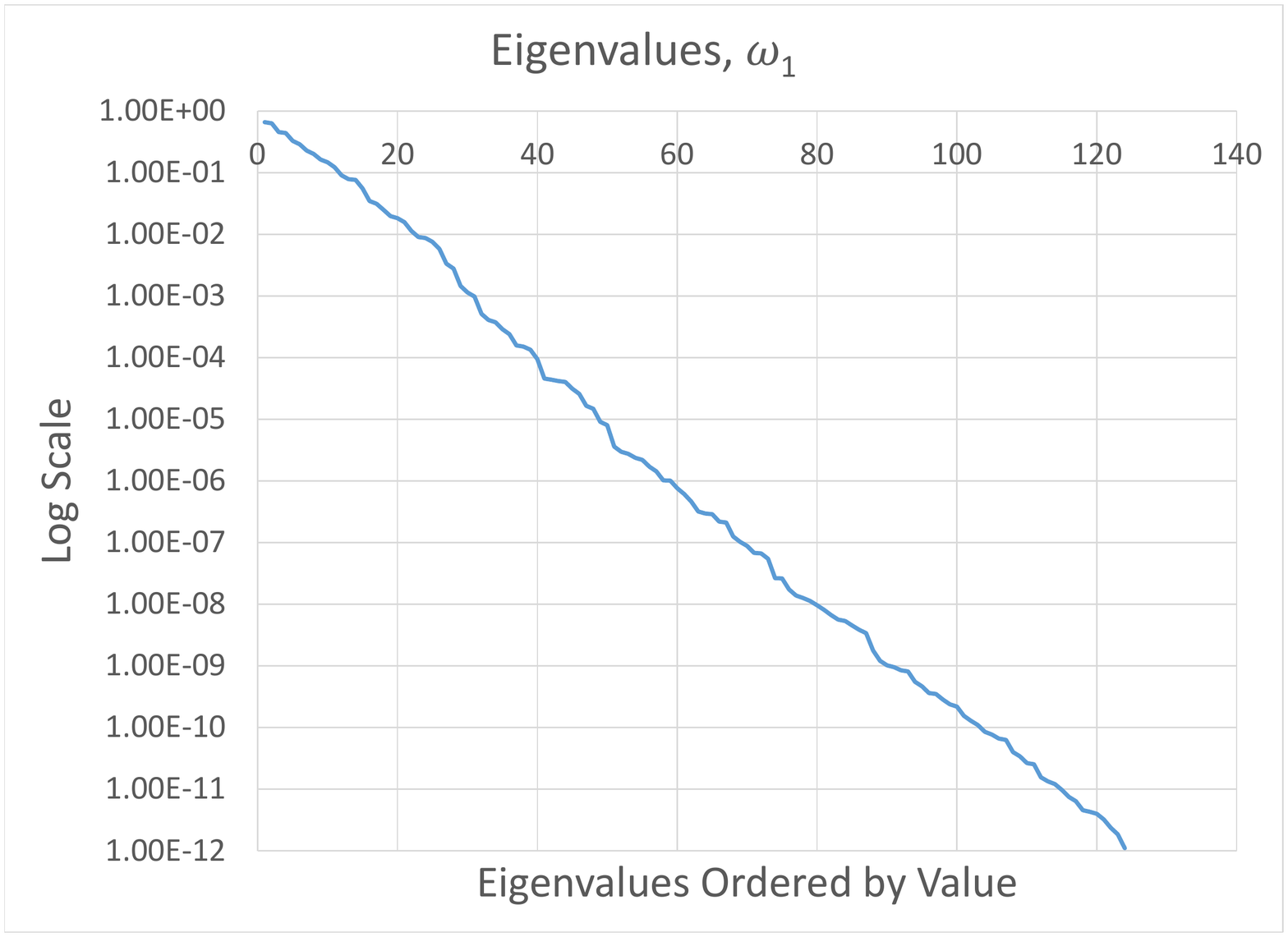} &
\includegraphics[trim={1.5cm 2.25cm 1.5cm 1.5cm},clip,height=5cm]{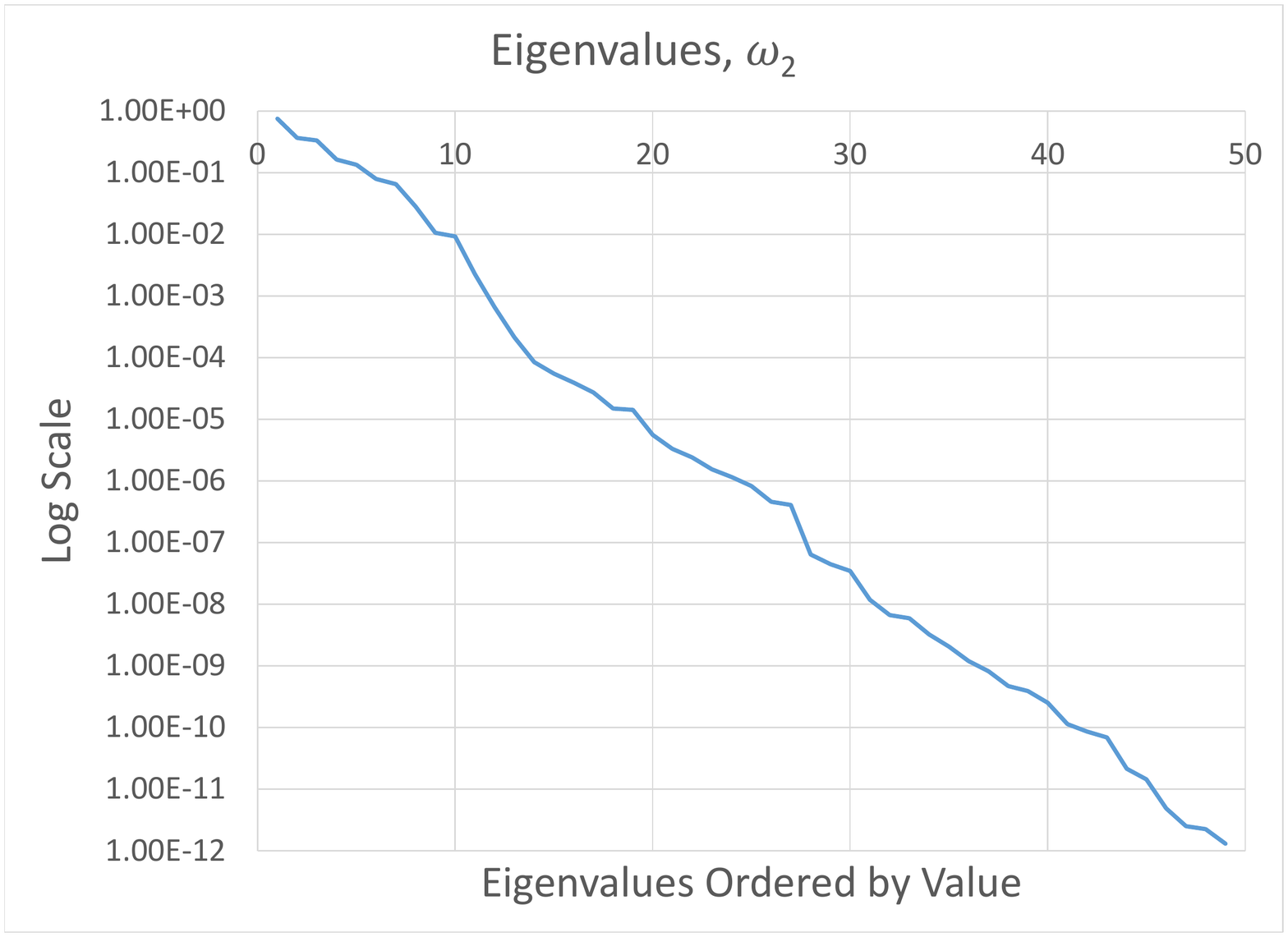}
\end{tabular}
\caption{Log plot of eigenvalues for $\omega_1$ (left) and $\omega_2$ (right) for conductivity in the matrix 100 and in the fibers 1.}
\label{fig:eigen12}
\end{figure}

The generalized eigenvalue problems are then solved using routines from the Intel\textsuperscript{\textregistered} Math Kernel Library (MKL). The problems are reduced to standard symmetric eigenvalue problems. In the following description of the methodology, we neglect the index $i$ keeping in mind that the problem has to be solved for each patch. Considering that the matrix $\mathrm{P}$ is symmetric and positive definite, a Cholesky factorization is performed
\begin{equation}\label{eq:cholesky}
\mathrm{P}=\mathrm{U}^{T}\mathrm{U}
\end{equation}
where $\mathrm{U}$ is an upper triangular matrix. Equation (\ref{eq:cholesky}) is then inserted in (\ref{eq:Q_ixlP_ix})
\begin{equation*}\label{eq:Q_ixlP_ix_2}
\mathrm{Q}\mathbf{x}=\lambda \mathrm{U}^{T}\mathrm{U}\mathbf{x}.
\end{equation*}
The inverse of $\mathrm{U}^{T}$ is computed and multiplied on both sides. On the left side there is also an identity matrix $\mathrm{I}=\mathrm{U}^{-1}\mathrm{U}$ inserted between $\mathrm{Q}$ and $\mathbf{x}$
\begin{equation}\label{eq:Q_ixlP_ix_3}
{\mathrm{U}^{T}}^{-1}\mathrm{Q}\mathrm{U}^{-1}\mathrm{U}\mathbf{x}=\lambda \mathrm{U}\mathbf{x}.
\end{equation}
With $\mathbf{y}=\mathrm{U}\mathbf{x}$ and $\mathrm{C}={\mathrm{U}^{T}}^{-1}\mathrm{Q}\mathrm{U}^{-1}$ (\ref{eq:Q_ixlP_ix_3}) becomes a standard eigenvalue problem
\begin{equation}\label{eq:C_iyly}
\mathrm{C}\mathbf{y}=\lambda \mathbf{y}.
\end{equation}
The eigenvalues for problem (\ref{eq:C_iyly}) are the same as for (\ref{eq:Q_ixlP_ix}), the eigenvectors $\mathbf{x}$ are obtained by solving $\mathrm{U}\mathbf{x}=\mathbf{y}$.

\begin{figure}[ht] 
\centering
\begin{tabular}{cc}
\includegraphics[height=3.15cm]{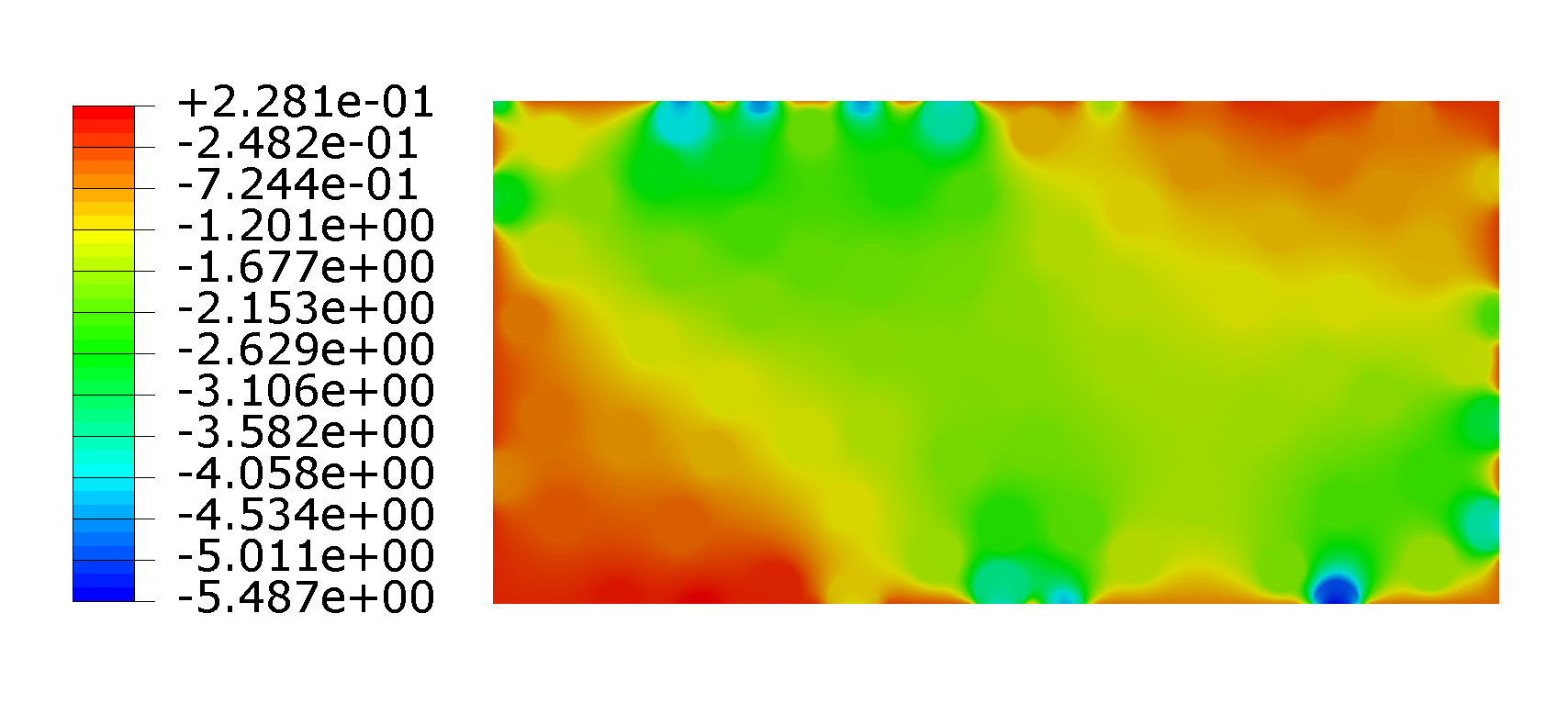}&  \includegraphics[height=3.15cm]{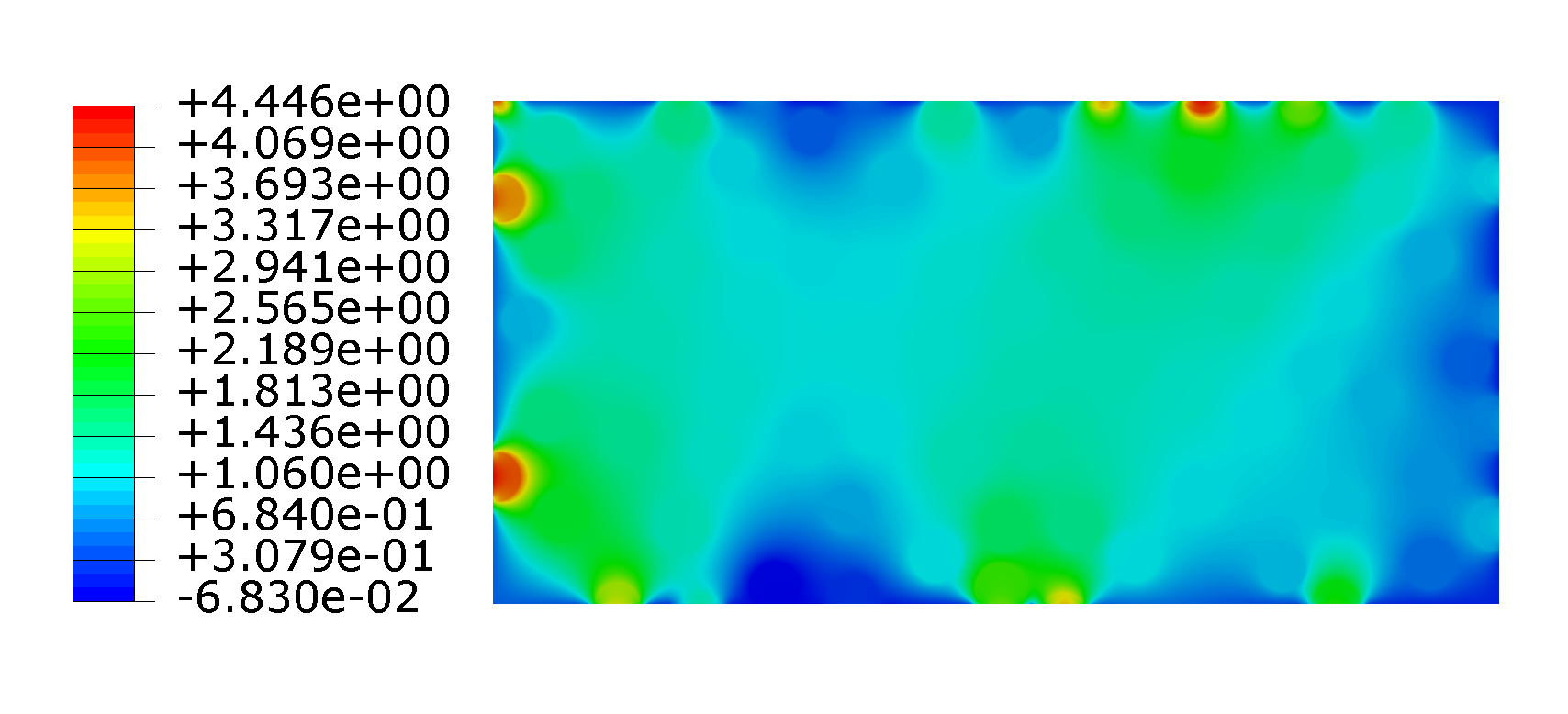}\\
\includegraphics[height=3.0cm]{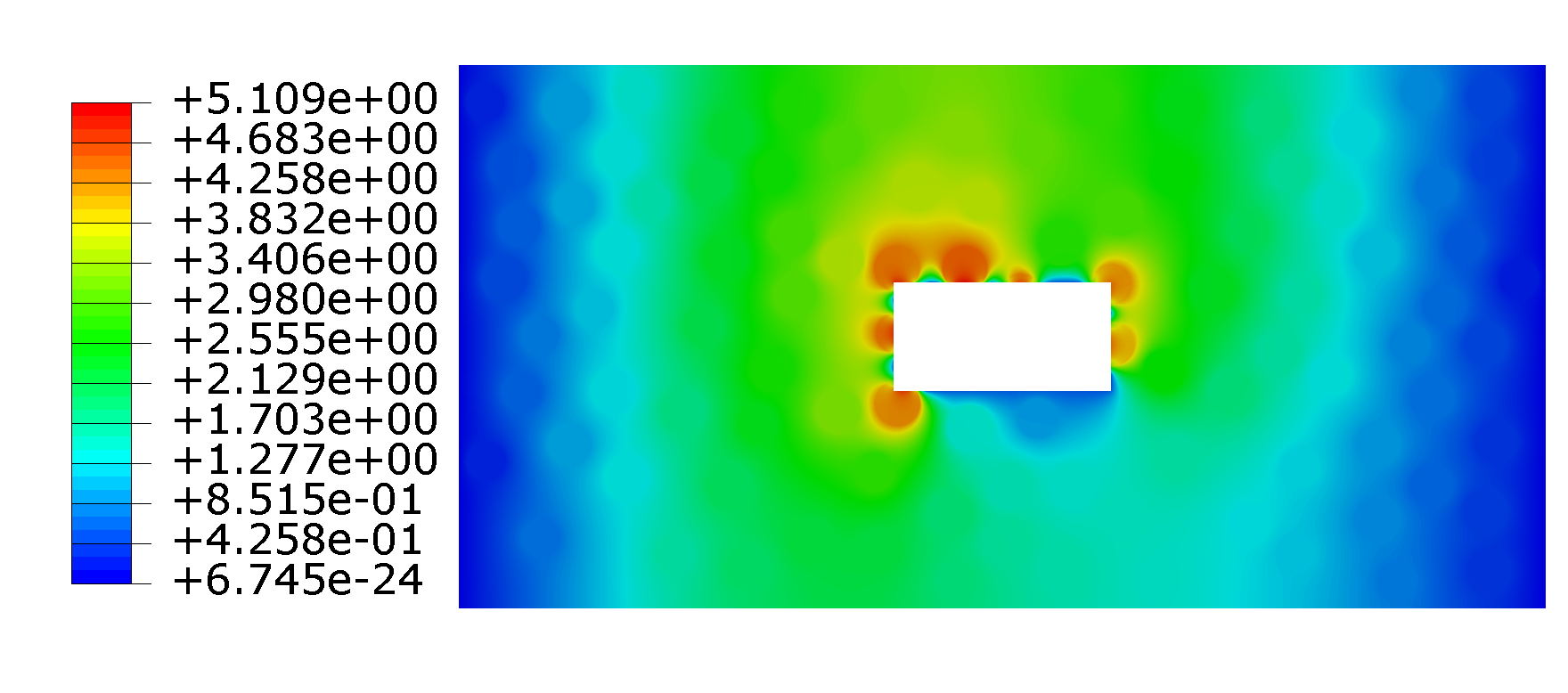}&  \includegraphics[height=3.0cm]{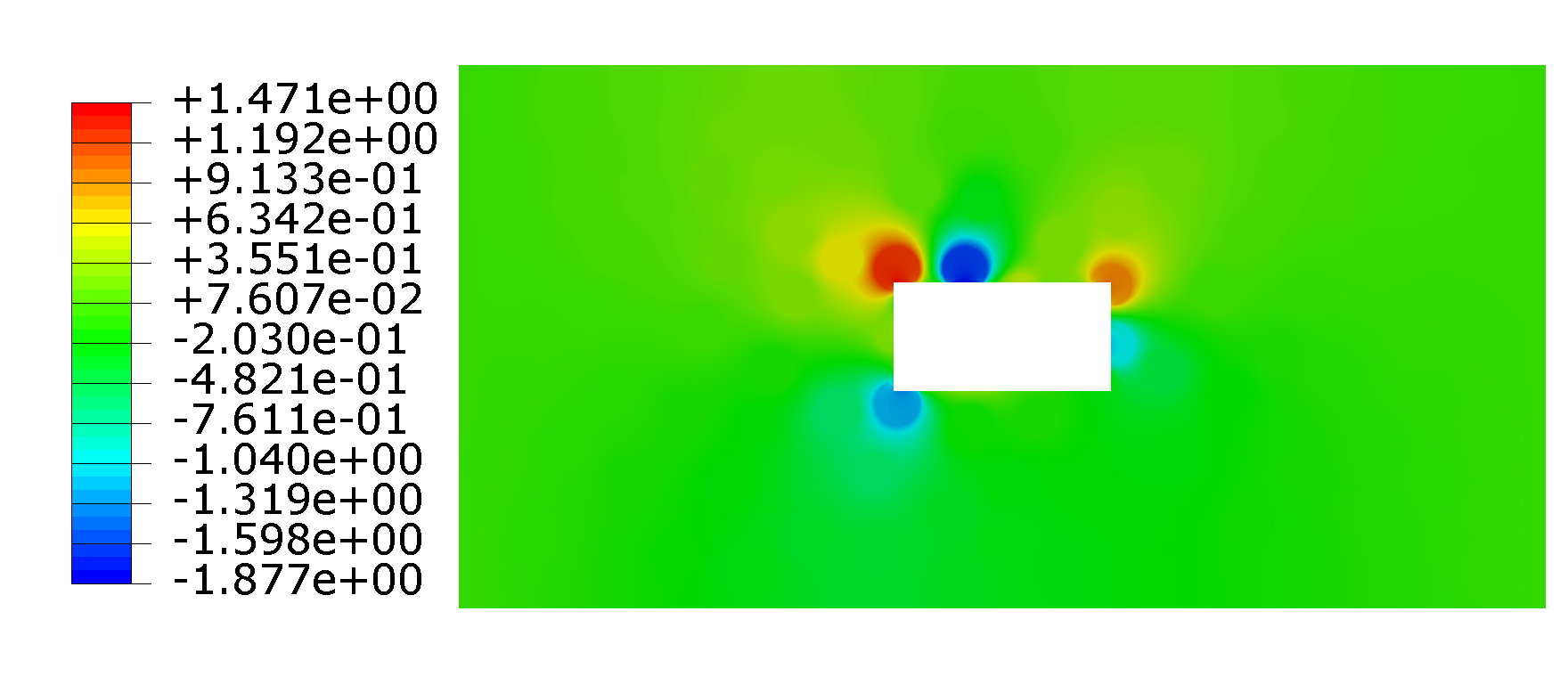}
\end{tabular}
\caption{Local spectral basis functions $\xi_1^5$ (top left) and $\xi_1^7$ (top right) over $\omega_1^*$ and $\xi_2^1$ (bottom left) and $\xi_2^6$ (bottom right) over $\omega_2^*$.}
\label{fig:shape12}
\end{figure}

The eigenvalues for the two patches are shown in Figure \ref{fig:eigen12}. Eigenfunctions corresponding to eigenvalues which are smaller than $10^{-12}$ are discarded. One can see in the graph on the right that the eigenvalues come in pairs, due to symmetry, except for the first one of $\omega_{2}$.  Figure \ref{fig:shape12} shows $\xi_1^4$, $\xi_1^8$, $\xi_2^1$, and $\xi_2^6$.  In addition to the local spectral basis functions, $V^{m_1}_{\omega_1^*}$ is augmented by the constant functions.

\subsection{Solving the Global System, $\mathrm{G}\mathbf{x}=\mathbf{r}$}

\begin{figure}[ht] 
\centering
\begin{tabular}{cc}
\includegraphics[height=3.4cm]{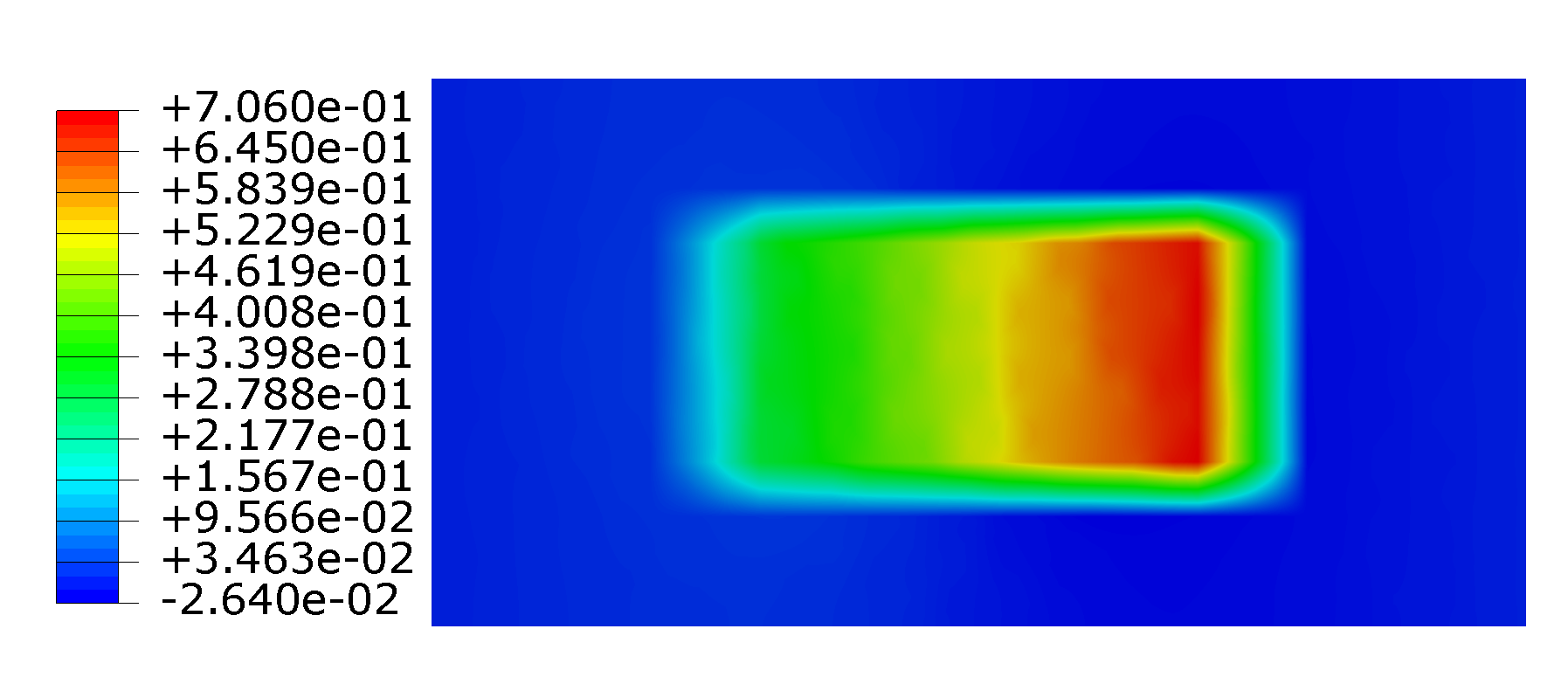} & \includegraphics[height=3.4cm]{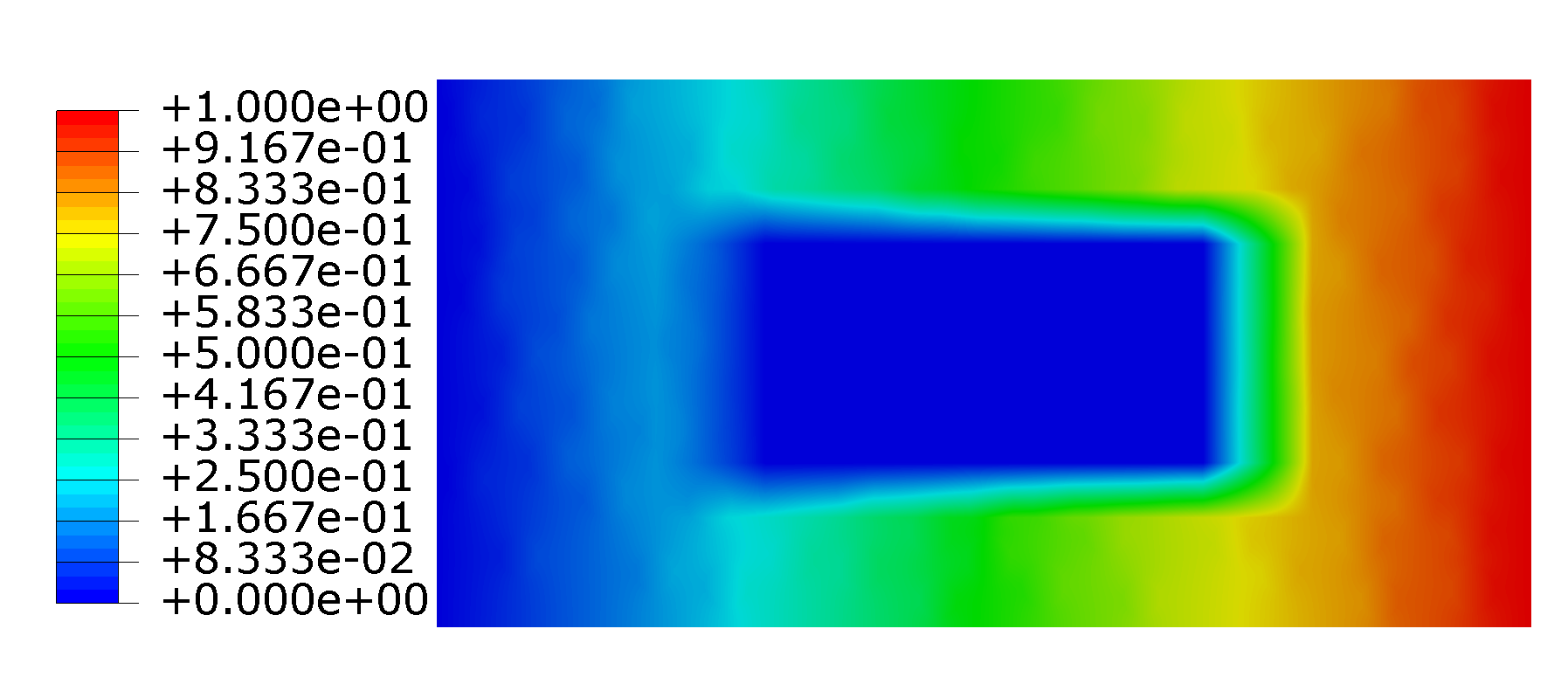}
\end{tabular}
\caption{The solutions $u^{G}$ (left) and $u^{F}$ (right).}
\label{fig:u_g_part}
\end{figure}

\begin{figure}[ht] 
\centering
\begin{tabular}{cc}
\includegraphics[height=3.4cm]{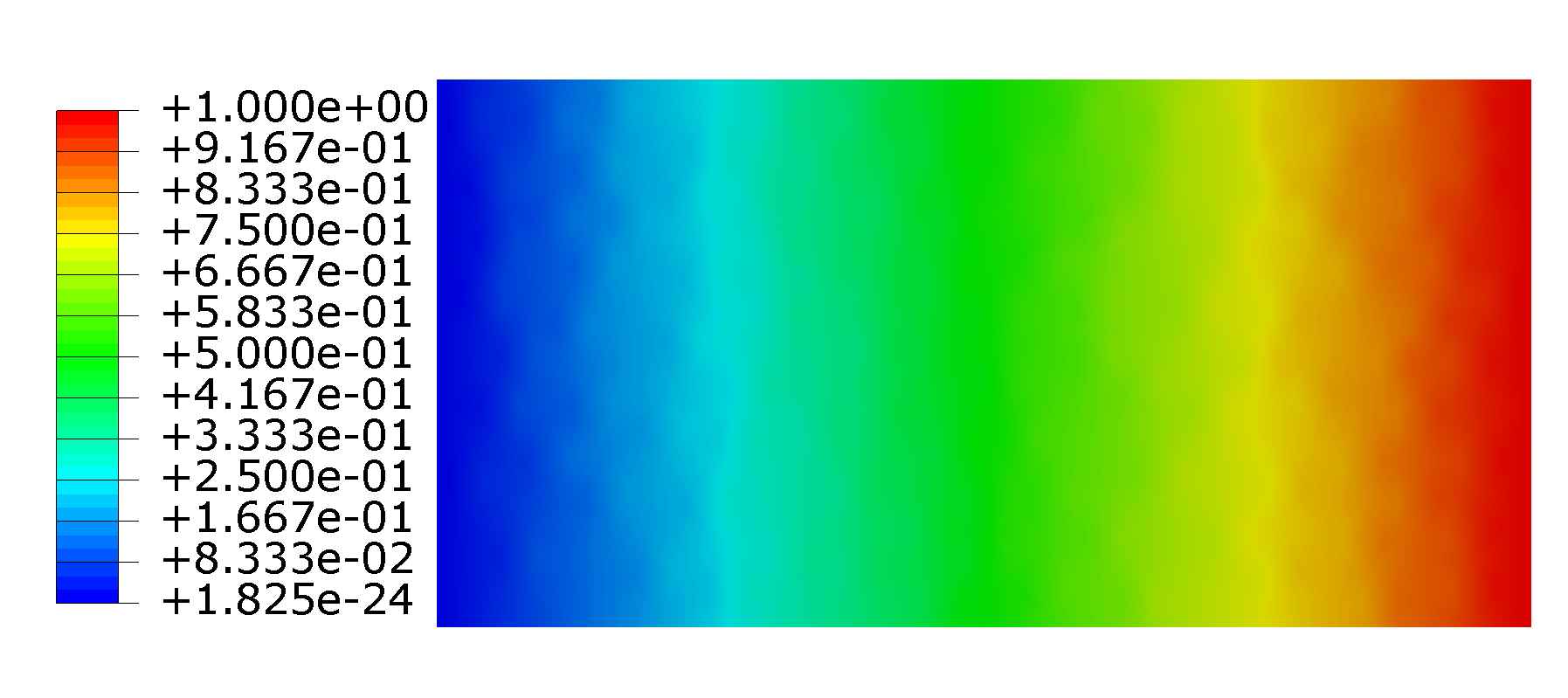} & \includegraphics[height=3.4cm]{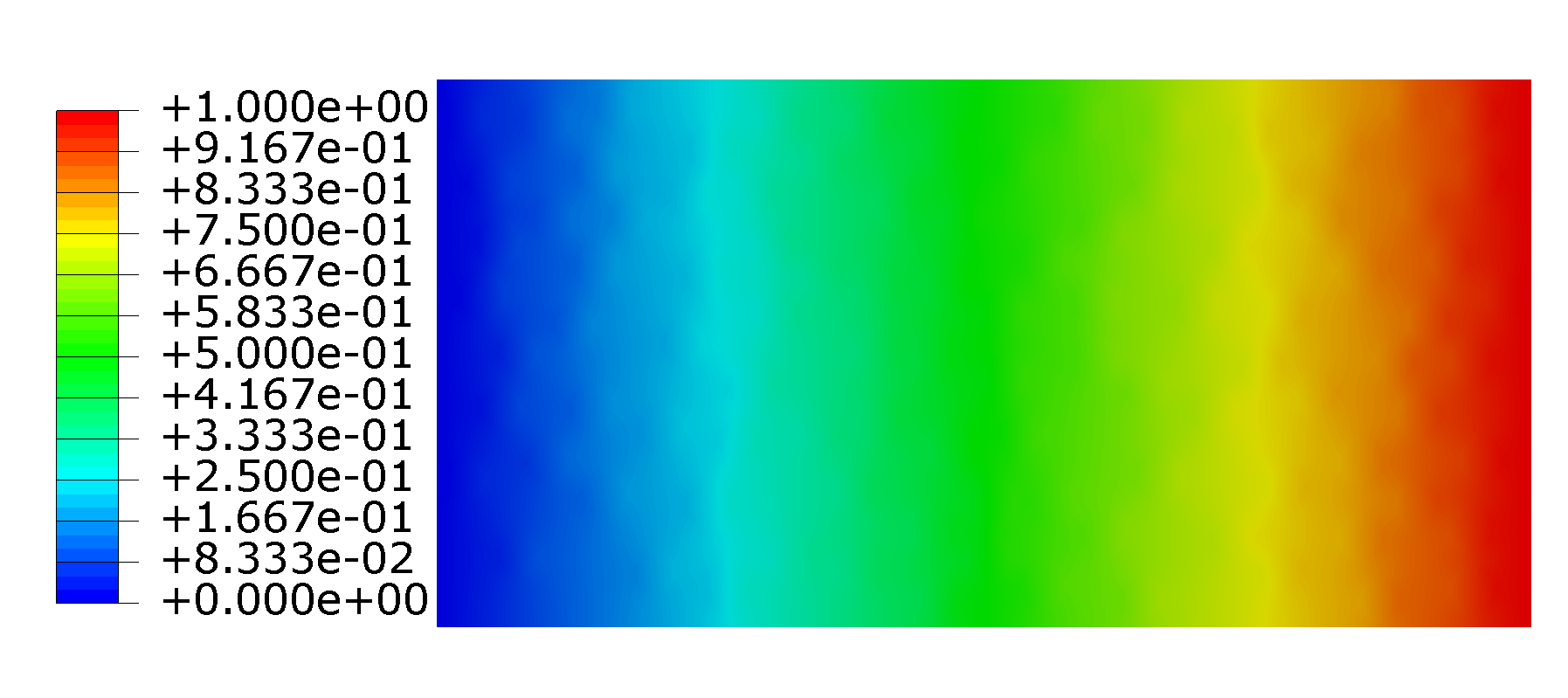}
\end{tabular}
\caption{The final solution $u_0$ computed by the MS-GFEM (left) and an FEM (right).}
\label{fig:final}
\end{figure}

The global stiffness matrix and right hand side entries, \eqref{eq:globstiffmat} and \eqref{eq:rhs}, are computed using integration on the underlying finite element mesh similarly to (\ref{eq:I_i}) using the saved element stiffness matrices from the particular solution runs. The constructions of the global stiffness matrix and right hand side are performed in parallel using OpenMP.
The global system $\mathrm{G}\mathbf{x}=\mathbf{r}$ is solved using the Pardiso solver from the Intel\textsuperscript{\textregistered} MKL. The solution vector $\mathbf{x}$ is used to recover $u^G$ as in \eqref{eq:aharmpart}.  The solution $u^G$ is shown in Figure \ref{fig:u_g_part} (left). The particular solution $u^F$ over the whole domain $\Omega$ consisting of the particular solutions over the patches glued together with the partition of unity functions is shown in Figure \ref{fig:u_g_part} (right).
The final solution is then recovered from $u_0=u^G+u^F$ shown in Figure \ref{fig:final} (left). The solution provided directly by a finite element method is displayed in Figure~\ref{fig:final} (right) for comparison with the MS-GFEM solution.
The MS-GFEM computed solution has energy $\norm{u_0}_{\mathcal{E}(\Omega)}^2=0.799227$.  The FEM, ``overkill,''solution has energy given by $\norm{u}_{\mathcal{E}(\Omega)}^2=0.799227.$  The relative error of the MS-GFEM solution compared to the FEM solution is given by
\begin{equation*}
   \frac{\norm{u-u_0}_{\mathcal{E}(\Omega)}}{\norm{u}_{\mathcal{E}(\Omega)}} = 2.52\times10^{-5}.
\end{equation*}


\section{Contrast Independent Convergence Study}
\label{sec:convgstudy}
In this section we present numerical results indicating that the MS-GFEM converges to the ``overkill'' FEM solution at an exponential rate that is independent of contrast between material properties of the matrix and the inclusions.
To demonstrate the contrast independent convergence rate we consider again problem \eqref{eq:pde_2}. We compute the MS-GFEM solution $u_0$ for different material contrasts between the matrix and inclusion materials. We adopt the convention $A:B$ where the first number is the matrix conductivity and the number following the colon is the inclusion conductivity. We carry out the simulation for 1:1, 1:1.5, 1:2, 1:10, 1:100, 1:200, and 1:1000 and 1.5:1, 2:1, 10:1, 100:1, 200:1, and 1000:1. 
All errors are computed relative to the FEM ``overkill'' solution $u$ for each choice of contrast. The relative error with respect to the energy norm is computed as
\begin{equation}\label{eq:convergence_error}
   E = \frac{\norm{u_0-u}_{\mathcal{E}(\Omega)}^2}{\norm{u}_{\mathcal{E}(\Omega)}^2}.
\end{equation} 
Tables \ref{tab:ee_x_1} and \ref{tab:ee_x_2} show the relative errors of the MS-GFEM solution $u_0$ computed across contrasts. The number of spectral basis functions from $\omega_1$ and $\omega_2$ used in computing $u_0$ are listed under their respective headings in the tables.
\begin{table}[!htbp]
\centering
\begin{tabular}{| r | r | c | c | c | c | c | c | c | c |}
\hline
\rule[-.3\baselineskip]{-3pt}{3ex}
$\omega_{1}$ & $\omega_{2}$&	$1000:1$ & $200:1$	& $100:1$	& $10:1$ &	$2:1$	& $1.5:1$	& $1:1$ \\ \hline
\rule[-.3\baselineskip]{-3pt}{3ex}
5	& 2  & 1.25E-01 & 1.25E-01 & 1.24E-01 & 1.23E-01 & 1.08E-01 & 1.01E-01 & 8.98E-02 \\
10	& 4  & 2.64E-02 & 2.61E-02 & 2.58E-02 & 2.29E-02 & 1.81E-02 & 1.58E-02 & 2.52E-02 \\
15	& 6  & 6.11E-03 & 6.04E-03 & 5.98E-03 & 6.60E-03 & 7.80E-03 & 6.86E-03 & 1.75E-03 \\
20	& 8  & 1.86E-03 & 1.85E-03 & 1.86E-03 & 3.35E-03 & 3.73E-03 & 2.78E-03 & 1.65E-03 \\
25	& 10 & 1.05E-03 & 1.06E-03 & 1.08E-03 & 1.72E-03 & 2.18E-03 & 1.74E-03 & 1.28E-03 \\
30	& 12 & 2.23E-04 & 2.83E-04 & 3.60E-04 & 1.31E-03 & 1.56E-03 & 1.08E-03 & 3.58E-04 \\
35	& 14 & 8.79E-05 & 1.19E-04 & 1.64E-04 & 8.61E-04 & 1.10E-03 & 7.30E-04 & 5.33E-05 \\
40	& 16 & 3.20E-05 & 4.63E-05 & 7.06E-05 & 4.82E-04 & 6.61E-04 & 4.43E-04 & 4.12E-05 \\
45	& 18 & 2.72E-05 & 3.69E-05 & 5.43E-05 & 3.35E-04 & 4.15E-04 & 2.63E-04 & 2.77E-05 \\
50	& 20 & 2.42E-05 & 3.04E-05 & 4.30E-05 & 2.37E-04 & 1.94E-04 & 1.23E-04 & 1.22E-05 \\
55	& 22 & 2.42E-05 & 2.75E-05 & 2.86E-05 & 1.19E-04 & 9.26E-05 & 5.48E-05 & 6.89E-06 \\
60	& 24 & 2.51E-05 & 2.49E-05 & 2.55E-05 & 4.39E-05 & 4.18E-05 & 2.87E-05 & 6.70E-06 \\
65	& 26 & 2.68E-05 & 2.57E-05 & 2.52E-05 & 3.52E-05 & 2.53E-05 & 1.69E-05 & 6.25E-06 \\
70	& 28 & 2.81E-05 & 2.67E-05 & 2.59E-05 & 2.59E-05 & 2.04E-05 & 1.46E-05 & 5.91E-06 \\ \hline
\end{tabular}
\vspace*{0.3cm}
\caption{Relative error versus dimension of spectral basis for different contrasts computed over a mesh of 1,560,058 elements. Here the A-harmonic local space $S^{n_i}_{\omega_i}$ is generated by hat functions defined on $\partial\omega^*_i$ of width 1 node.}
\label{tab:ee_x_1}
\end{table}
\begin{table}[!htbp]
\centering
\begin{tabular}{| r | r | c | c | c | c | c | c | c | c |}
\hline
\rule[-.3\baselineskip]{-3pt}{3ex}
$\omega_{1}$ & $\omega_{2}$ & $1:1000$ & $1:200$ & $1:100$ & $1:10$ & $1:2$ & $1:1.5$ & $1:1$ \\ \hline
\rule[-.3\baselineskip]{-3pt}{3ex}
5	& 2  & 2.89E+00 & 1.23E+00 & 8.49E-01 & 2.16E-01 & 8.01E-02 & 8.11E-02 & 8.98E-02 \\ 
10	& 4  & 1.39E+00 & 6.08E-01 & 4.34E-01 & 1.17E-01 & 2.94E-02 & 2.35E-02 & 2.52E-02 \\ 
15	& 6  & 8.48E-01 & 4.40E-01 & 3.31E-01 & 9.30E-02 & 1.52E-02 & 9.46E-03 & 1.75E-03 \\ 
20	& 8  & 2.09E-01 & 1.09E-01 & 6.83E-02 & 1.17E-02 & 5.28E-03 & 3.35E-03 & 1.65E-03 \\ 
25	& 10 & 2.87E-02 & 1.48E-02 & 1.12E-02 & 8.01E-03 & 3.13E-03 & 1.98E-03 & 1.28E-03 \\ 
30	& 12 & 4.05E-03 & 3.52E-03 & 3.42E-03 & 4.34E-03 & 1.85E-03 & 1.14E-03 & 3.58E-04 \\ 
35	& 14 & 2.57E-03 & 2.52E-03 & 2.52E-03 & 2.57E-03 & 1.12E-03 & 6.99E-04 & 5.33E-05 \\ 
40	& 16 & 9.46E-04 & 9.65E-04 & 1.01E-03 & 1.29E-03 & 7.21E-04 & 4.31E-04 & 4.12E-05 \\ 
45	& 18 & 5.21E-04 & 4.72E-04 & 4.32E-04 & 4.81E-04 & 3.61E-04 & 2.43E-04 & 2.77E-05 \\ 
50	& 20 & 4.55E-04 & 3.93E-04 & 3.66E-04 & 2.66E-04 & 1.31E-04 & 8.51E-05 & 1.22E-05 \\ 
55	& 22 & 5.59E-05 & 1.35E-04 & 3.18E-04 & 2.04E-04 & 8.73E-05 & 5.36E-05 & 6.89E-06 \\ 
60	& 24 & 3.40E-05 & 4.89E-05 & 6.34E-05 & 1.67E-04 & 5.33E-05 & 3.06E-05 & 6.70E-06 \\ 
65	& 26 & 2.63E-05 & 3.33E-05 & 3.86E-05 & 4.86E-05 & 2.81E-05 & 1.78E-05 & 6.25E-06 \\ 
70	& 28 & 2.50E-05 & 2.24E-05 & 2.33E-05 & 3.37E-05 & 1.83E-05 & 1.22E-05 & 5.91E-06 \\ \hline 
\end{tabular}
\vspace*{0.3cm}
\caption{Relative error versus dimension of spectral basis for different contrasts computed over a mesh of 1,560,058 elements. Here the A-harmonic local space $S^{n_i}_{\omega_i}$ is generated by hat functions defined on $\partial\omega_i^*$ of width 1 node.}
\label{tab:ee_x_2}
\end{table}
\begin{table}[!htbp]
\centering
\begin{tabular}{| r | r | c | c | c | c | c | c |}
\hline
\rule[-.3\baselineskip]{-3pt}{3ex}
$\omega_{1}$ & $\omega_{2}$&	$1000:1$ & $100:1$	& $1:1$ & $1:100$ &	$1:1000$  \\ \hline
\rule[-.3\baselineskip]{-3pt}{3ex}
5	& 2  & 1.25E-01 & 1.24E-01 & 9.03E-02 & 8.02E-01 & 2.75E+00 \\
10	& 4  & 2.64E-02 & 2.59E-02 & 2.49E-02 & 2.45E-01 & 7.74E-01 \\
15	& 6  & 6.08E-03 & 5.97E-03 & 1.74E-03 & 1.40E-01 & 3.80E-01 \\
20	& 8  & 1.85E-03 & 1.83E-03 & 1.63E-03 & 5.84E-02 & 2.11E-01 \\
25	& 10 & 1.03E-03 & 1.05E-03 & 1.48E-03 & 1.32E-02 & 2.07E-02 \\
30	& 12 & 2.08E-04 & 2.87E-04 & 3.71E-04 & 4.74E-03 & 5.50E-03 \\
35	& 14 & 8.32E-05 & 1.23E-04 & 3.92E-05 & 2.44E-03 & 3.19E-03 \\
40	& 16 & 3.19E-05 & 5.88E-05 & 3.79E-05 & 7.73E-04 & 7.48E-04 \\
45	& 18 & 1.97E-05 & 3.93E-05 & 3.02E-05 & 3.06E-04 & 3.89E-04 \\
50	& 20 & 9.85E-06 & 2.56E-05 & 1.57E-05 & 2.87E-04 & 3.66E-04 \\
55	& 22 & 2.98E-06 & 9.70E-06 & 3.23E-06 & 2.45E-04 & 3.45E-05 \\
60	& 24 & 1.45E-06 & 6.43E-06 & 2.54E-06 & 4.27E-05 & 2.29E-05 \\
65	& 26 & 7.42E-07 & 4.22E-06 & 1.81E-06 & 1.96E-05 & 1.03E-05 \\
70	& 28 & 6.05E-07 & 2.37E-06 & 1.04E-06 & 1.04E-05 & 7.39E-06 \\ \hline
\end{tabular}
\vspace*{0.3cm}
\caption{Relative errors versus dimension of spectral basis of the MS-GFEM solution using quadratic finite elements on a mesh with 388,189 elements and 1,156,882 nodes. Here the A-harmonic local space $S^{n_i}_{\omega_i}$ is generated by hat functions defined on $\partial\omega_i^*$ of width 3 nodes.}
\label{tab:ee_x_1_quad}
\end{table}

\begin{figure}[htb] 
\centering
\includegraphics[height=5cm]{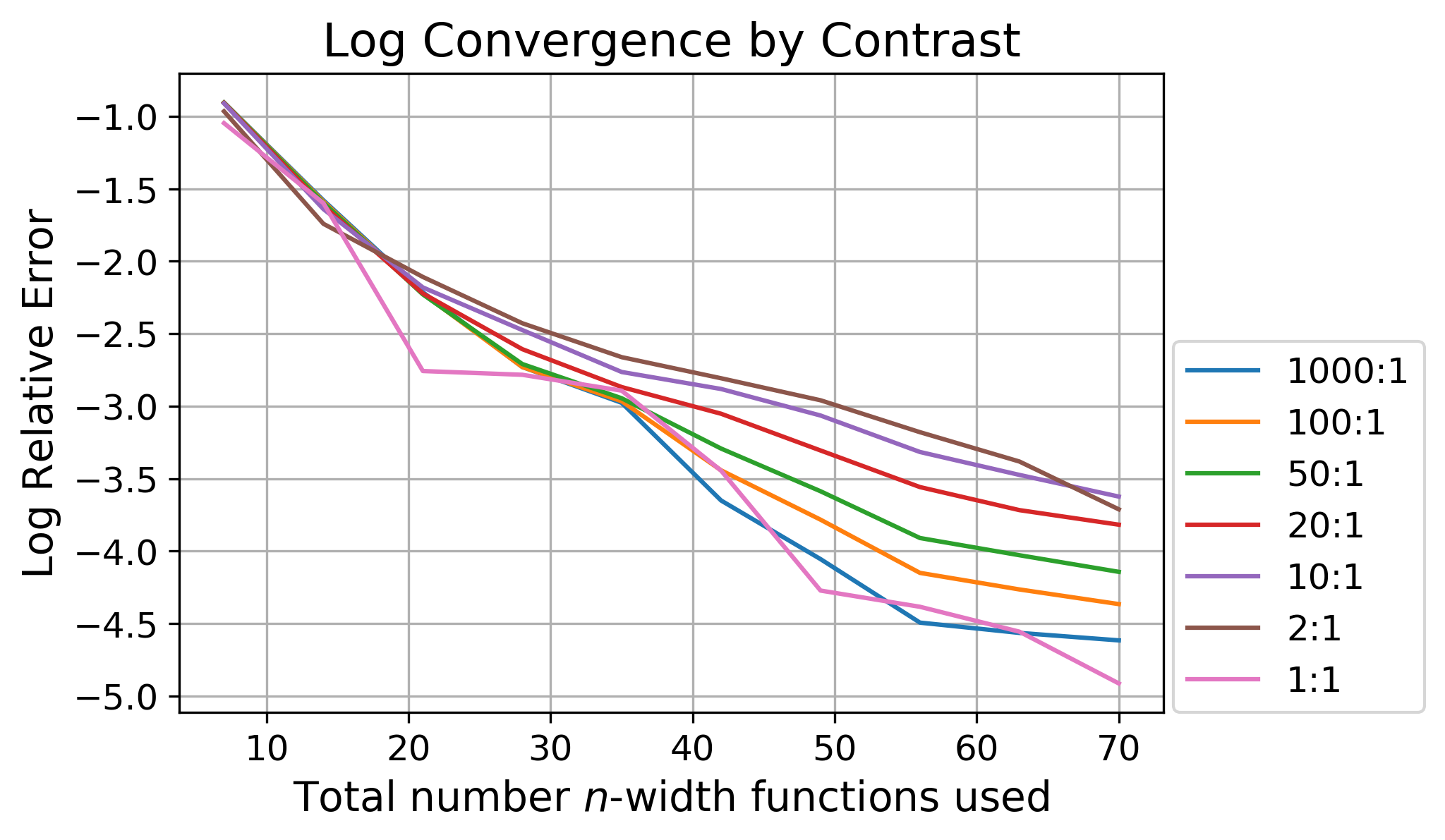}
\caption{Log plot of relative error of MS-GFEM solution (computed with bilinear FEM) versus number of local basis functions used.  Each line is for a different matrix conductivity and inclusions with conductivity 1.}
\label{fig:ee_x_1}
\end{figure}

\begin{figure}[!htbp] 
\centering
\includegraphics[height=5cm]{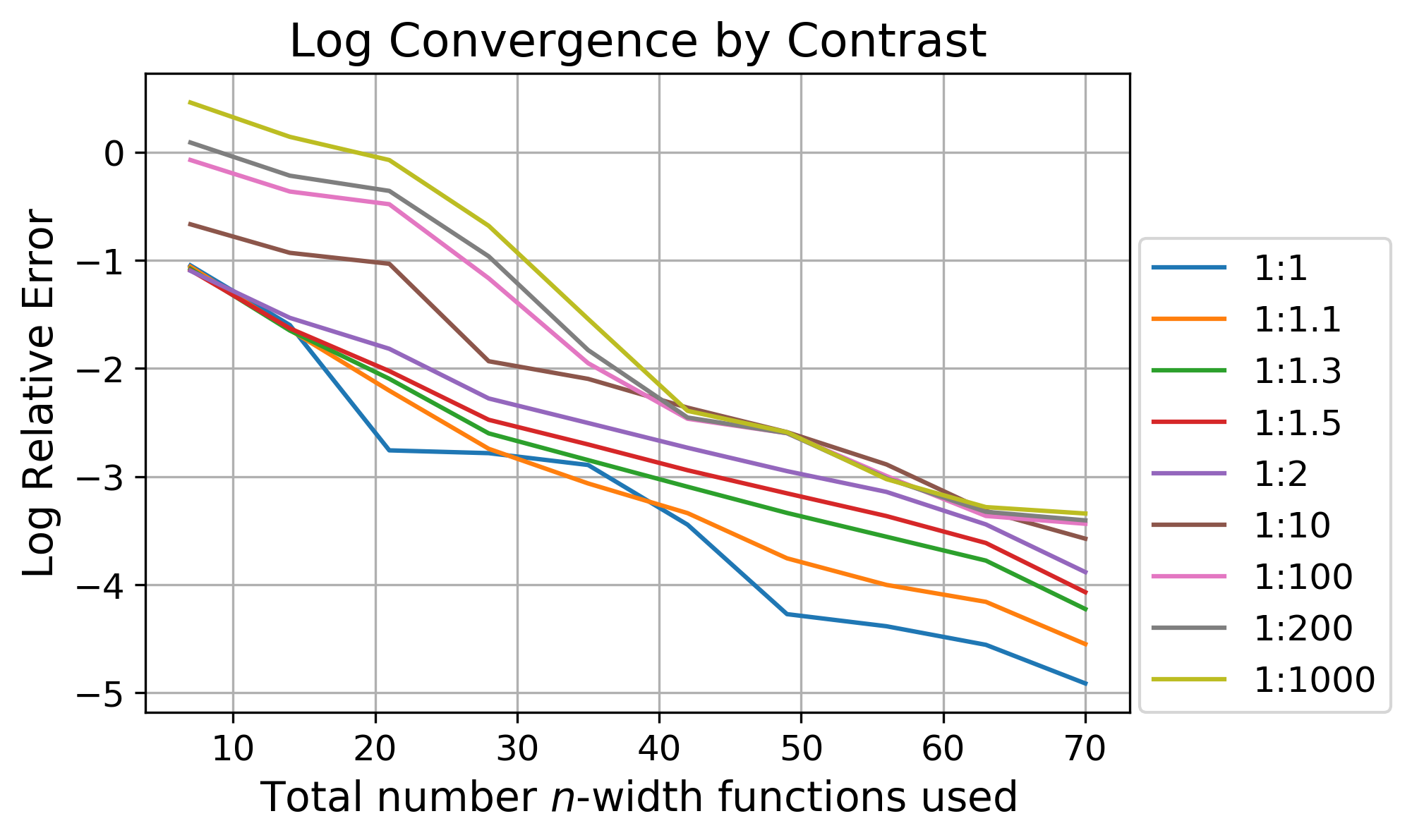}
\caption{Log plot of relative error of MS-GFEM solution (computed with bilinear FEM) versus number of local basis functions used.  Each line is for matrix of conductivity 1 and inclusions of different conductivity.}
\label{fig:ee_1_x}
\end{figure}

\begin{figure}[!htb] 
\centering
\includegraphics[height=5cm]{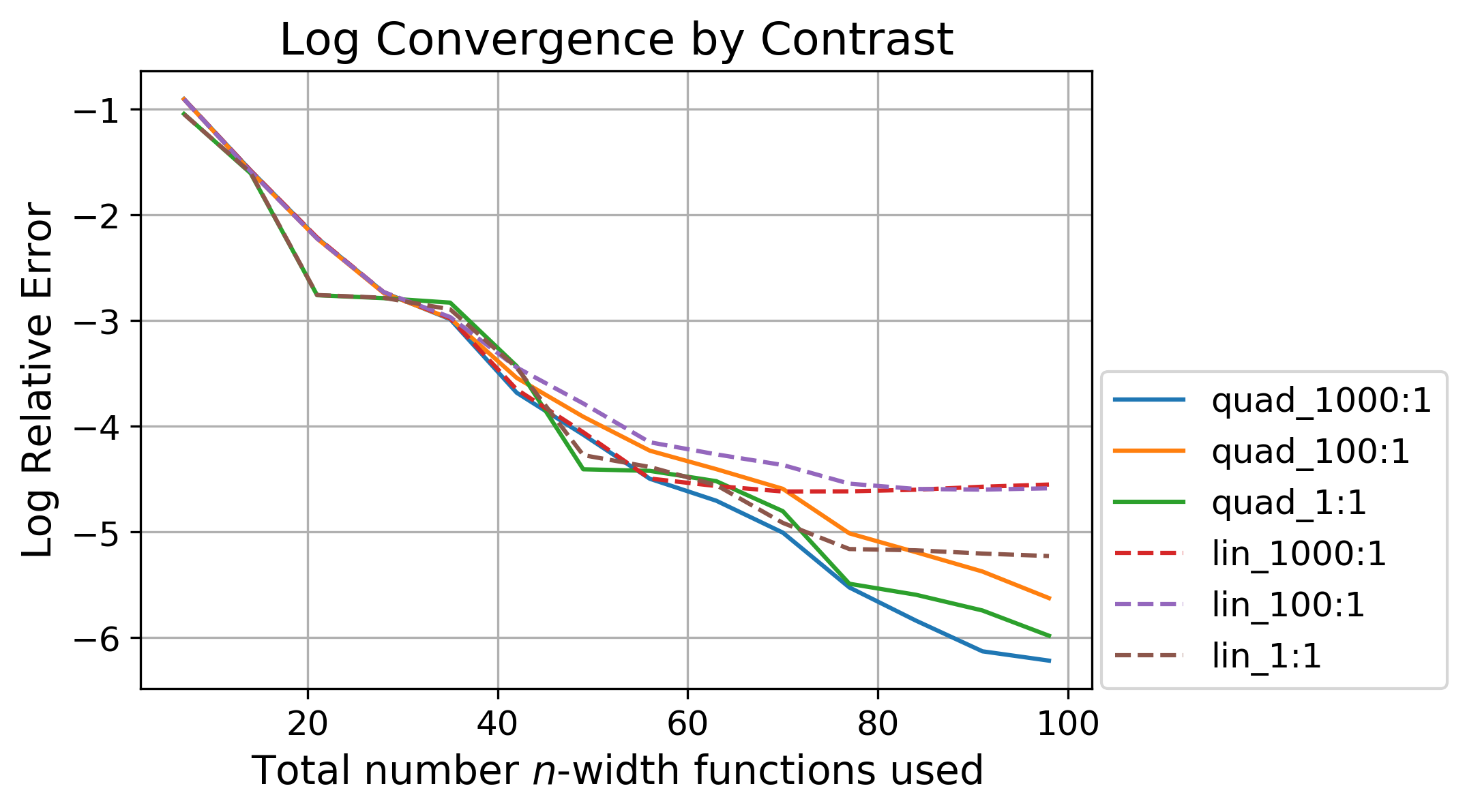}
\includegraphics[height=5cm]{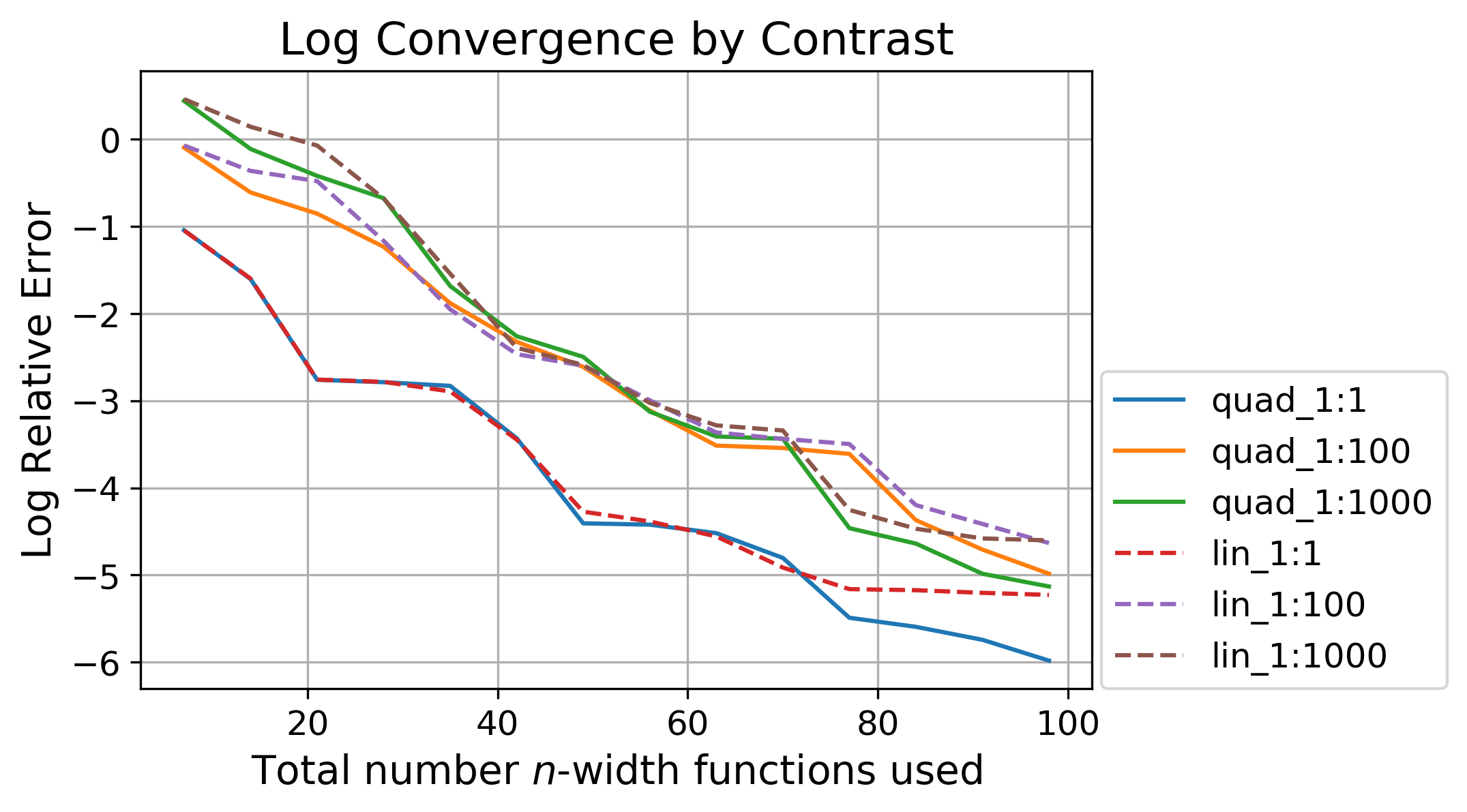}
\caption{Comparison of quadratic and linear finite elements for computing MS-GFEM global relative errors. Log plot of relative error versus number of local basis functions used.  Each line is the relative error for quadratic or linear elements of contrast listed with matrix conductivity followed by inclusion conductivity.}
\label{fig:quad_vs_lin}
\end{figure}

In Figure~\ref{fig:ee_x_1} the log of the relative errors is plotted against the dimension of the spectral basis for inclusions with conductivity 1 and matrix with conductivity 1, 2, 10, 20, 50, 100, and 1000. 
Figure \ref{fig:ee_1_x} shows the same for inclusions with conductivity 1, 1.1, 1.3, 1.5, 2, 10, 100, 200, and 1000, and matrix with conductivity 1. These figures show the exponential convergence of the error seen with respect to the dimension of the local approximation space. Here the contrast between matrix and inclusions is seen to not influence the convergence rate.

Next we increase the order of finite elements to quadratic to see what effect it has on the rate of convergence of the relative error versus dimension of the local spectral bases. The relative errors of the MS-GFEM solution using quadratic finite elements on a mesh with 388,189 elements and 1,156,882 nodes are now used. The spectral bases are now generated by an A-harmonic approximation space $S_{\omega_i^*}^{n_i}$ with boundary data given by hat functions defined on $\partial\omega^*_i$ that are supported over three boundary nodes. The convergence rates for the associated MS-GFEM for several differenct contrasts between matrix and particles are displayed in Figure \ref{fig:quad_vs_lin}. The convergence rate for MS-GFEM using quadratic elements is also compared to the convergence rate for MS-GFEM using linear elements in Figure \ref{fig:quad_vs_lin}. The simulations show that the change of element type has no effect of the convergence rates associated with the dimension of the local spectral basis up to a dimension of about fifty. Beyond fifty the convergence rate of MS-GFEM using linear elements the convergence rate goes to zero, whereas the convergence rate for MS-GFEM with quadratic elements continues to converge exponentially with dimension of up to eighty. The relative error decreases below $10^{-5}$ for all cases using spectral bases of dimension greater than seventy. Here the convergence for MS-GFEM using linear elements flattens and is zero after dimension fifty due to the linear FEM approximation error incurred in the generation of the discrete A-harmonic spaces $S_{\omega^*_i}^{m_i}$.

\subsection{Eigenvalue Decay and Size of Oversampling Domain}\label{sec:eval decay}
\begin{figure}
    \centering
    \includegraphics[trim={1cm 1.75cm 1cm 1cm},clip,width=0.7\linewidth]{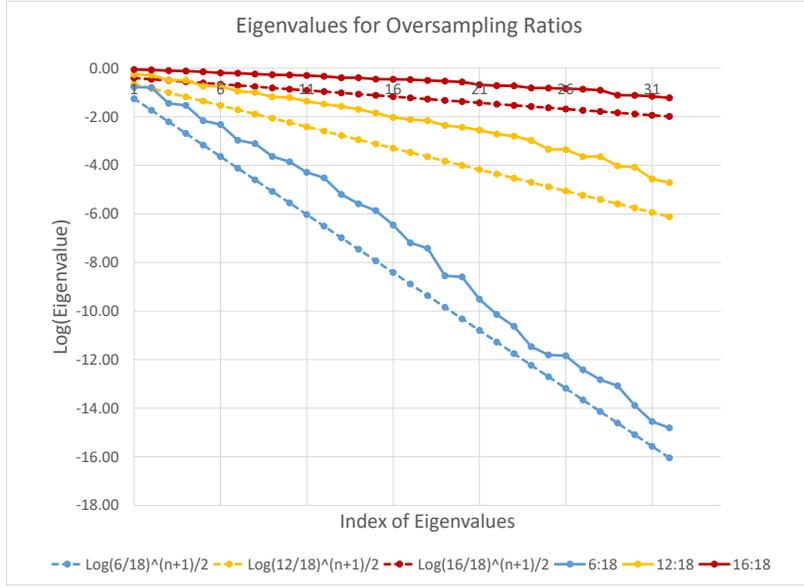}
    \caption{Log plot of eigenvalues for different subdomains $\omega_1$ holding $\omega_1^*$ fixed. In this case the side lengths of $\omega_1$ are proportionately changed. Conductivity 1 in matrix and 1000 in inclusions. Eigenfunctions generated with 19-node width hat functions as boundary data for $H_A(\omega^*)$.}
    \label{fig:domain ratio fixed outer}
\end{figure}

The key feature of MS-GFEM is the construction of optimal local solution spaces through oversampling. The size of the oversampled patch $\omega^*$ as compared to $\omega$ has a direct impact on the decay rate of the eigenvalues which bound the local errors (cf. Theorem~\ref{thm:nwidthdecay}). We investigate the size of the oversampling domain $\omega^*$ relative to the domain $\omega$ used in approximation. We first consider oversampeling  domain $\omega^*$ given by a disk of radius $R$ and the smaller concentric disk $\omega$ of radius $\rho$, $R>\rho$. We look at the A-harmonic space of solutions of \eqref{eq:aharmpde} with a homogeneous coefficient $A(x)=I$ for $\omega_i^*=\omega^*$. A straight forward calculation shows that the $n$-width eigenvalues are given by 
\begin{equation}\label{rate}
\lambda_{n} = \left(\frac{\rho}{R}\right)^{2n}.
\end{equation}
Now keeping in mind that the earlier numerical results  given here for particle composites illustrate the contrast independent convergence of the MS-GFEM, and so the $n$-width eigenvalues should also be independent of contrast.  Therefore we argue that a relationship like $\lambda_{n} = \left(\frac{\rho}{R}\right)^{2n}$ should also hold for the heterogeneous case. Our computations show that this is approximately true. To illustrate the point we consider rectangular domains with coefficient matrix $A(x)$ associated with particle composites as in Figure \ref{fig:100_inclusions_01} with contrast $1$ in the matrix and $1000$ in the inclusions. We denote by $\rho$ and $R$ the longest side lengths of $\omega$ and $\omega^*$, respectively. Figure \ref{fig:domain ratio fixed outer} shows the decay of the $n$-width eigenvalues with respect the ratio $\rho/R$  to for three values $16/18$, $12/18$, and $6/18$. We see from these simulations that there is good agreement between the formula \eqref{rate} given by the solid lines and the numerically computed convergence rate given by the data points. 


\section{Reduction of Computational Work by Penetration; an Effective Computation of Finite Dimensional Subspaces of $H_A(\omega_i^*)$} \label{sec:aharmonic}

The primary numerical work in generating the best local basis is not in the solution of the eigenvalue problem but in the numerical generation of the local A-harmonic subspace over which the eigenvalue problem is solved.
Numerically generating the A-harmonic subspaces $S^{n_i}_{\omega_i^*}\subset H_{A}(\omega_i^*)$ can be computationally expensive for fine meshes.  Defining the boundary data as in \eqref{eq:aharm} with a single boundary element for the support of the hat function $h_i^j$ results in as many problems to solve as there are boundary nodes on $\partial\omega_i^*$.  Most importantly, we see from Figure~\ref{fig:hat1} that the resulting function $w_i^j$ given by the A-harmonic extension of the boundary hat functions with support over one boundary element has very small support and gradient inside $\omega^\ast$.  However for A-harmonic extensions of boundary hat functions having 301 nodes as support we see from Figure~\ref{fig:hat5} that the support of this function penetrates further into in $\omega^\ast$.  This is the phenomenon of penetration where A-harmonic extensions of boundary data with less oscillation have a gradient that penetrates further into $\omega^*$, than A-harmonic extensions of more oscillatory boundary data, see \cite{BLS08} and \cite{LS16}. 
In light of the exponential decay of the $n$-width eigenfunctions in the energy norm over $\omega$ it becomes clear that the space spanned by the lower $n$-width modes are generated by linear combinations of A-harmonic extensions of coarse boundary data, ie., boundary hat functions of large support. This motivates the usage of significantly wider hat functions on the boundary such that their  A-harmonic extension will have significant variation that penetrates into $\omega$.  This means that less functions need to be used for the computation of the local space $S^{n_i}_{\omega^*_i}$ if we use boundary data given by hat functions having support over several boundary nodes.

We see from Theorem \ref{thm:nwidthdecay} that the local error controlled by the $n$-width eigenvalues, decay at worst exponentially so the number of local spectral basis functions $\xi_i^j$ can be reduced with little effect on the error. For the example of section \ref{sec:numerics}, the eigenvalues decay rapidly, see Figure~\ref{fig:eigen12}, and the 100-th and 50-th eigenvalues dip below $10^{-15}$ for $\omega_1$ and $\omega_2$, respectively. Therefore only a handful of the associated $n$-width eigenfunctions need to be used in the local approximation spaces to achieve excellent accuracy; the rest can be discarded.

Thus if we want to approximate $M$ dimensional local subspaces using $n$ spectral basis functions ($n<M$)  we are motivated to use $M$ wider boundary boundary hat functions and their A harmonic extensions.  This results in fewer problems to solve and smaller spectral matrices in \eqref{eq:spectralMatrices} since the number of A-harmonic functions $w_i^j$ can be of the same order of wide hat functions used for boundary data.   This gives a strategy for reducing the computational work in computing the local basis and for computing the spectral basis.

When the number of boundary hat functions is reduced from $N$ to $M$ by defining boundary hat functions with larger support then the computational time for constructing $S^{m_i}_{\omega_i^*}$ is reduced by a factor of $\frac{M}{N}$.  The time for computing each entry of the matrices (\ref{eq:Q_ijk}) and (\ref{eq:P_ijk}) is approximately the same.  Due to symmetry of the matrices, only the upper triangular part and the diagonal need to be computed. If N hat functions were used then $\frac{N(N+1)}{2}$ entries have to be computed. If that number is reduced to $M$ then
\begin{equation*}
\frac{M(M+1)}{2}=\frac{N\rho\left(N\rho+1\right)}{2}
\end{equation*}
entries have to be computed, where $\rho = M/N$. For large values of $N$  
\begin{equation*}
\frac{M(M+1)}{2} \sim \frac{N^2\rho^2}{2}
\end{equation*} 
and the computational time for generating the matrices is reduced by a factor of $\left(\frac{M}{N}\right)^{2}$.  The number of hat functions used for each successive widening progresses as $\{N,N/2,N/3,N/4,\ldots\}$.  For example, if the support of the hat functions is increased from 1 to 3 nodes, then the total number of boundary value problems \eqref{eq:aharm} to solve is $M=N/2$ and the time to fill the spectral matrices is about $1/4$ the time as with single node hat functions.  We illustrate these observations with the following example.  

We solve the problem \eqref{eq:pde_2}
as before
on the domain shown in Figure~\ref{fig:100_inclusions_01} using hat functions with varying widths to generate the local solution spaces $S^{n_i}_{\omega_i^*}$.  We label solutions $u_{k},\ k=1,3,5,7,9,\dots$, where $k$ represents the number of nodes in the support of the hat function used as boundary data for problem \eqref{eq:aharm}.  

The matrix material has conductivity $A_{1}$ and the inclusion material has conductivity $A_{2}$.  We compute solutions to \eqref{eq:pde_2} using two different sets of material properties.  For Case 1
\begin{equation*}\label{eq:A_i}
A_{1} = 
\begin{pmatrix*}
1 & 0 \\
0 & 1   
\end{pmatrix*}\ \ \ 
A_{2} = 
\begin{pmatrix*}
100 & 0 \\
0 & 100   
\end{pmatrix*}
\end{equation*}
and for Case 2 we reverse the material properties to
\begin{equation*}\label{eq:A_i2}
A_{1} = 
\begin{pmatrix*}
100 & 0 \\
0 & 100   
\end{pmatrix*}\ \ \ 
A_{2} = 
\begin{pmatrix*}
1 & 0 \\
0 & 1   
\end{pmatrix*}.
\end{equation*}

\begin{figure}[!htb] 
\centering
\includegraphics[trim={1.5cm 1.5cm 1.5cm 1.5cm},clip,height=6cm]{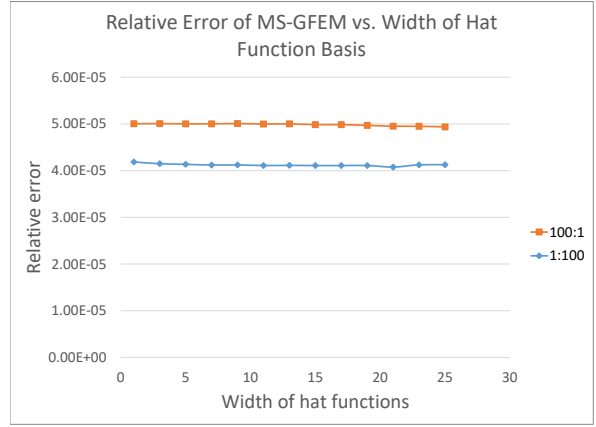}
\caption{Plot of MS-GFEM relative global error versus the number of nodes in the support of the boundary hat functions used  to generate the local spaces $S_{\omega_i^*}^{n_i} \subset H_A(\omega_i^*)$. Here the two cases of matrix to inclusion conductivities of 100 to 1 and 1 to 100 are shown.}
\label{fig:whf1}
\end{figure}

The relative error of the final solutions $u_k$ as compared to the ``overkill'' FEM solution $u$ is computed with respect to the energy norm as
\begin{equation}\label{eq:energy_error_01}
E_{k} = \frac{\norm{u_{k}-u}_{\mathcal{E}(\Omega)}^2}{\norm{u}_{\mathcal{E}(\Omega)}^2}.
\end{equation}
The relative global errors versus the support of the hat functions used to generate the local approximation spaces are reported in Table~\ref{tab:whf1} and are graphed in Figure~\ref{fig:whf1}. The result is that the global error is insensitive to the support of the local hat functions up to a width  of twenty five nodes. Here  the ``Error 1'' column lists the relative error for a matrix of conductivity $1$ and inclusion of conductivity $500$. The ``Error 2'' column lists the relative errors for a matrix of conductivity $500$ and inclusion of conductivity $1$.   The table also shows the size of the spectral matrices $\mathrm{P}$ and $\mathrm{Q}$ and the number of eigenfunctions used in the computation of $u_k$.  
\begin{table}[!ht]
\centering
\begin{tabular}{| r | r | r | r | r | r | r |}
\cline{4-7}
\multicolumn{3}{ c }{\rule{0pt}{3ex}}  &  \multicolumn{2}{|p{15ex}|}{\centering Eigenvalue System Size} & \multicolumn{2}{|p{15ex}|}{\centering Number of Eigenfunctions}  \\ \hline
\rule{0pt}{3ex}\# & Error 1    & Error 2     & \multicolumn{1}{p{1.1cm}|}{\centering$\omega_1$} &  \multicolumn{1}{p{1.1cm}|}{\centering$\omega_2$}  & \multicolumn{1}{p{1.1cm}|}{\centering$\omega_1$} &  \multicolumn{1}{p{1.1cm}|}{\centering$\omega_2$}  \\ \hline
\rule[-.3\baselineskip]{-3pt}{3ex}
1    & 4.185E-05    & 5.004E-05      & 3840  & 961    & 124   & 49\\
3    & 4.147E-05    & 5.007E-05      & 1920   & 481    & 124   & 49\\
5    & 4.135E-05    & 5.002E-05      & 1280   & 321    & 124   & 49\\
7    & 4.120E-05    & 5.003E-05      & 960   & 241     & 124   & 49\\
9    & 4.122E-05    & 5.008E-05      & 768   & 193     & 124    & 49\\
11   & 4.109E-05    & 4.998E-05      & 640   & 161     & 124    & 49 \\
13   & 4.113E-05    & 5.001E-05      & 548   & 138     & 124   & 49\\
15   & 4.108E-05    & 4.984E-05      & 480   & 121     & 124   & 49\\
17   & 4.108E-05    & 4.984E-05      & 426   & 107     & 124   & 49\\
19   & 4.109E-05    & 4.969E-05      & 384   & 97     & 124   & 49\\
21   & 4.077E-05    & 4.951E-05      & 349   & 88     & 123     & 49\\
23   & 4.124E-05    & 4.950E-05      & 320   & 81     & 123    & 49\\
25   & 4.126E-05    & 4.936E-05      & 295   & 74     & 123    & 47\\ \hline
\end{tabular}
\vspace*{0.3cm}
\caption{Relative error for global MS-GFEM solutions using hat functions with different numbers of boundary nodes for support.  The first column indicates the number of nodes in the support of the hat functions.  The first error is measured in the energy norm forwhen the matrix material has conductivity 1 and the inclusions 100, and the second errors are when the conductivities are switched.  Also the size of the spectral matrices $\mathrm{P}$ and $\mathrm{Q}$  (number of hat functions used), and the number of local spectral basis functions used in the approximation (chosen such that $10^{-12} \leq \lambda \leq 1$).}
\label{tab:whf1}
\end{table}
In this example the error \eqref{eq:energy_error_01} is seen to be constant, see Figure \ref{fig:whf1}, while the size of the spectral matrices decays rapidly for the first few sets of wider hat functions.  The number of eigenfunctions selected for the computation begins decreasing after widening the  boundary hat functions beyond 9 nodes. Collecting these observations, we conclude that this method of generating the local A-harmonic space $S^{n_i}_{\omega_i^*}$ is an effective method to reduce computational costs while maintaining acceptable error tolerance. 

\section{Span of Hat Functions Used Directly as Local Approximation Space}
\label{sec:oversampledgfem}
\begin{figure}[htb]
    \centering
    \includegraphics[trim={1cm 1.5cm 1cm 1cm},clip,height=6cm]{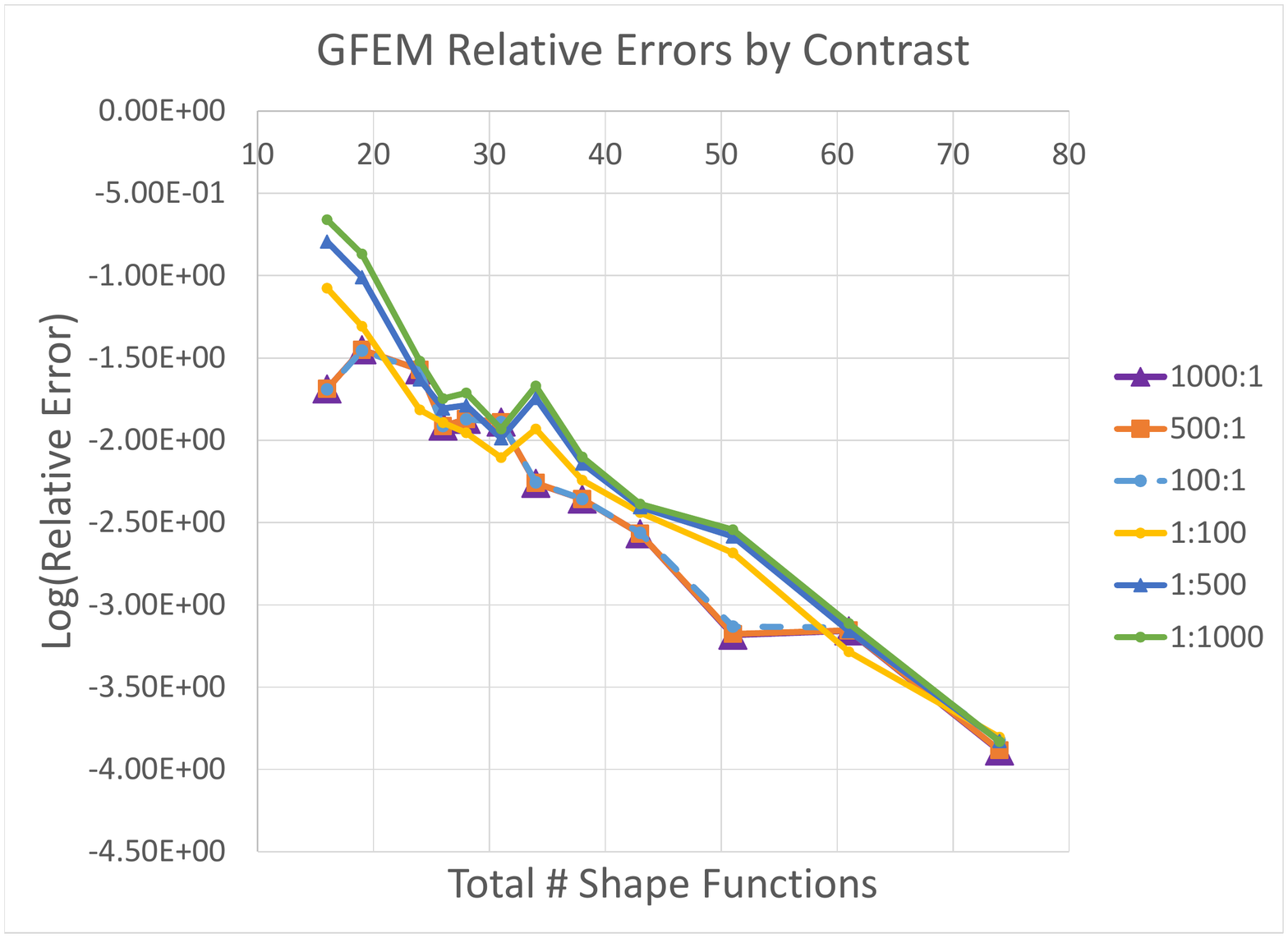}
    \includegraphics[trim={1cm 1.5cm 1cm 1cm},clip,height=6cm]{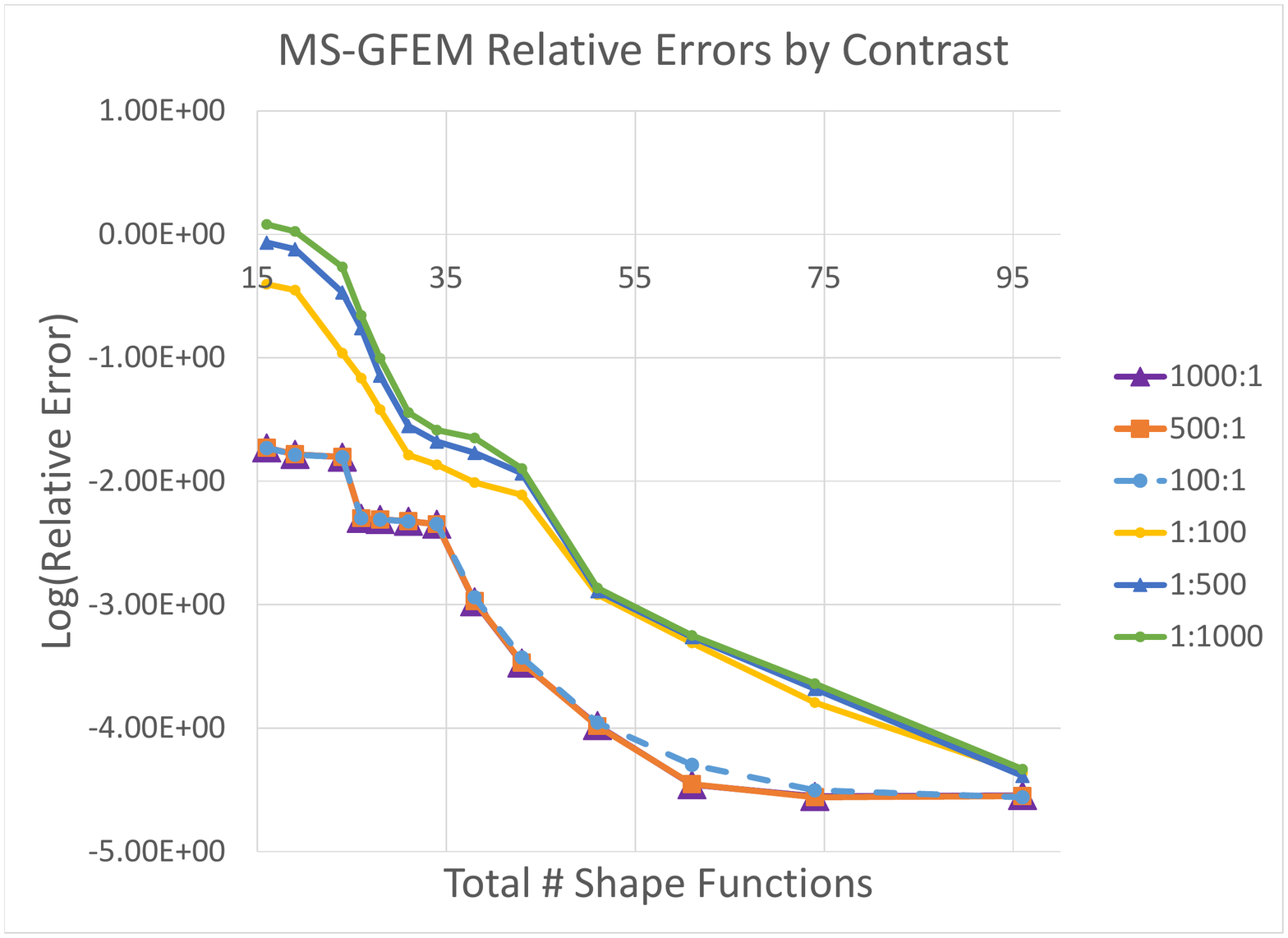}
    \\
    \includegraphics[trim={1cm 1.5cm 1cm 1cm},clip,height=6cm]{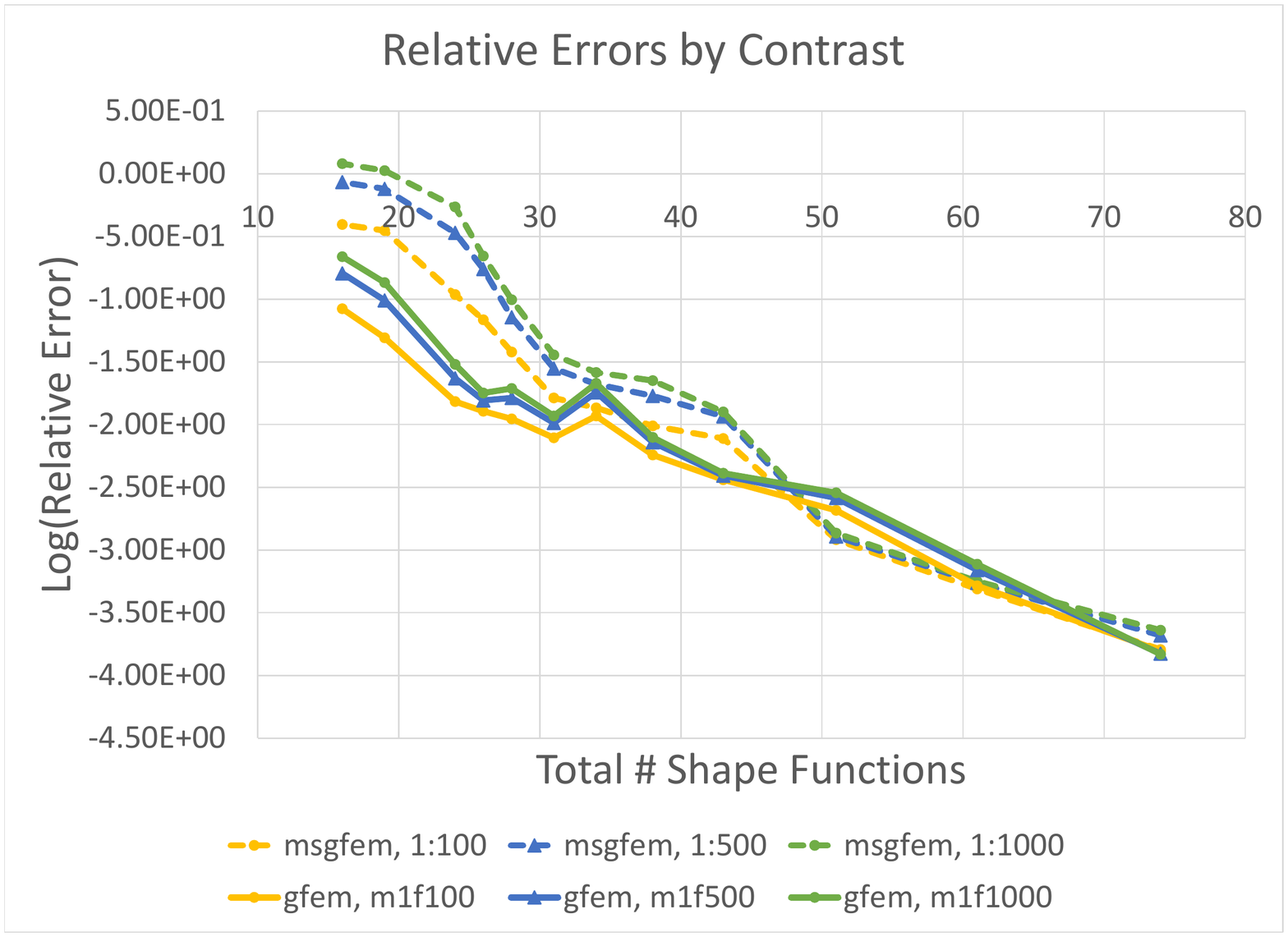}
    \includegraphics[trim={1cm 1.5cm 1cm 1cm},clip,height=6cm]{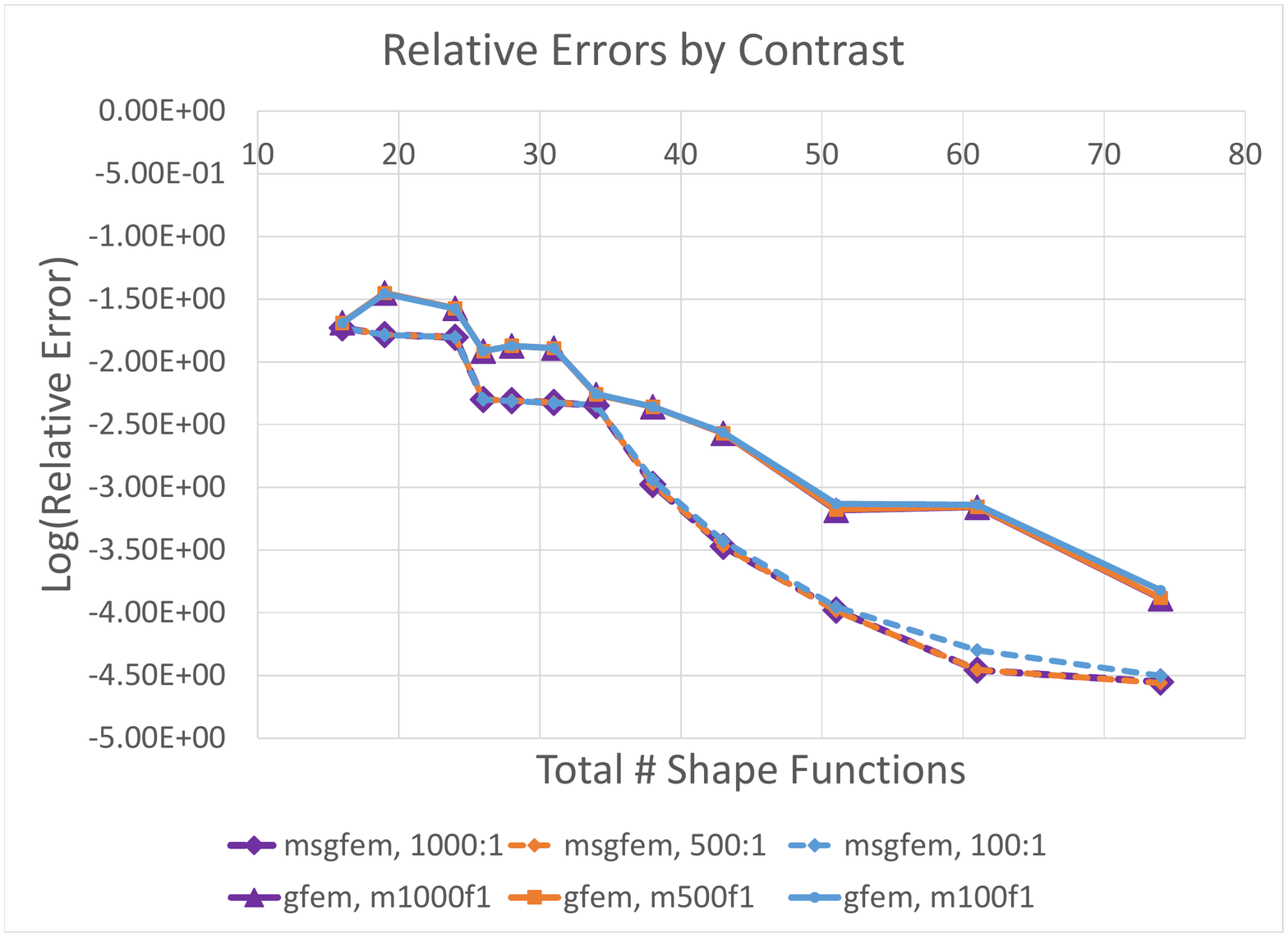}
    \caption{Relative errors of oversampled-GFEM and MS-GFEM solutions for different contrasts versus total number of shape functions used. Errors computed relative to the FEM overkill solution. }
    \label{fig:hats_gfem 100:1}
\end{figure}

In the last section we exploited the penetration of fields inside $\omega^*_i$ obtained by A-harmonically extending hat functions of large support on the boundary. Motivated by this we investigate the direct use of such functions as local approximations on $\omega_i\subset \omega^*_i$.  We call the GFEM method using these ``oversampled'' local basis the {\em oversampled-GFEM}. We investigate the decay rate of the global error using {\em oversampled-GFEM}. Here we prescribe hat functions on the boundary $\partial\omega^*_i$ and compute their A-harmonic extensions into $\omega^*_{i}$ as in \eqref{eq:solSpaces1} and \eqref{eq:solSpaces2}. Then the basis for the local approximation space is given by the restriction of these functions to $\omega_i\subset\omega_i^*$. As in the last section we take the hat functions on the boundary as elements of a partition of unity on the boundary of $\omega^*_i$.  We try coarser to finer boundary partitions of unity to generate our local subspace for the {\em oversampled-GFEM}. We find as we increase the number of basis elements by refining the boundary partition of unity we generate global approximate solutions with accuracy that are nearly exponentially decreasing with the dimension of the local basis. Tables \ref{tab:hats_gfem}, \ref{tab:hats_gfem2}, and \ref{tab:hats_gfem3} provides a comparison of the relative global error obtained using the same dimension of local approximation space given by $m_i$ for the oversampled-GFEM and $V_{\omega_i}^{m_i}$ for the MS-GFEM.  Each table addresses a different material contrast. Here we are comparing with MS-GFEM that use $n$-width eigenfunctions  generated from a local high dimensional approximation space $S^{n_i}_{\omega^*_i}$ with $n_i>>m_i$ that is spanned by A-harmonic extensions of hat functions with support over one boundary element. The tables show that for all contrasts the MS-GFEM does better than oversampled-GFEM at reducing the relative error in the global approximation as the number of bases elements are decreased. However as the dimension of the local approximating space for oversampled-GFEM and MS-GFEM is reduced we see that the oversampled-GFEM  delivers global solutions with relative error comparable to MS-GFEM. This numerical study is carried out for the matrix-inclusion composite with contrasts $1:100$ and $100:1$, $1:500$ and $500:1$, and $1:1000$ and $1000:1$. Here the coarseness of the hat function A-harmonic local basis used in the oversampled-GFEM is measured by the support of the boundary hat functions, see Tables \ref{tab:hats_gfem}, \ref{tab:hats_gfem2}, and \ref{tab:hats_gfem3}. The relative error versus total degrees of freedom given in Tables \ref{tab:hats_gfem}, \ref{tab:hats_gfem2}, and \ref{tab:hats_gfem3} are plotted in Figure \ref{fig:hats_gfem 100:1}. Here the top left chart in Figure \ref{fig:hats_gfem 100:1} shows the exponential decay of oversampled-GFEM with respect to global degrees of freedom for different material contrasts and the top right chart displays the exponential decay of  MS-GFEM with respect to global degrees of freedon for different contrasts. The bottom left chart in Figure \ref{fig:hats_gfem 100:1} plots oversampled-GFEM together with MS-GFEM and their decay with global degrees of freedom with matrix conductivity equal to unity and different particle conductivity's and  the lower right hand chart plots the decay of the error of both oversampled-GFEM and  MS-GFEM with respect to global degrees of freedom for different matrix conductivity's and particle conductivity unity. Note that both MS-GFEM and oversampled-GFEM converge exponentially with respect to the global degrees of freedom for all contrasts but MS-GFEM gives the most rapid decrease in approximation error with respect to global degrees of freedom.

\begin{table}[!htb]
\centering
\begin{tabular}{| c | c | c | c | c | c | c | c | c | c |}
\cline{4-7}
\multicolumn{3}{ c }{\rule{0pt}{3ex}}  &  \multicolumn{2}{|p{15ex}|}{\centering 100:1} & \multicolumn{2}{|p{15ex}|}{\centering 1:100} \\
\hline
\rule[-.3\baselineskip]{-3pt}{3ex}
Hat Widths & $\omega_1$ & $\omega_2$ & GFEM & MS-GFEM & GFEM & MS-GFEM	\\ \hline
\rule[-.3\baselineskip]{-3pt}{3ex}
69  & 109 & 28  & 2.89e-05 & 2.91e-05                    & 4.63e-05 & 2.25e-05                    \\
99  & 76  & 20  & 1.11e-04 & 2.75e-05                    & 1.50e-04 & 4.34e-05                    \\
129 & 59  & 15  & 1.51e-04 & 3.14e-05                    & 1.57e-04 & 1.61e-04                    \\
159 & 48  & 13  & 7.28e-04 & 5.03e-05                    & 5.19e-04 & 4.90e-04                    \\
189 & 40  & 11  & 7.38e-04 & 1.11e-04                    & 2.07e-03 & 1.21e-03                    \\
219 & 34  & 9   & 2.74e-03 & 3.73e-04                    & 3.64e-03 & 7.74e-03                    \\
249 & 30  & 8   & 4.40e-03 & 1.15e-03                    & 5.73e-03 & 9.77e-03                    \\
279 & 27  & 7   & 5.55e-03 & 4.52e-03                    & 1.17e-02 & 1.36e-02                    \\
309 & 24  & 7   & 1.29e-02 & 4.72e-03                    & 7.85e-03 & 1.63e-02                    \\
339 & 22  & 6   & 1.34e-02 & 4.88e-03                    & 1.11e-02 & 3.80e-02                    \\
369 & 20  & 6   & 1.22e-02 & 5.02e-03                    & 1.28e-02 & 6.86e-02                    \\
399 & 19  & 5   & 2.66e-02 & 1.57e-02                    & 1.53e-02 & 1.09e-01                    \\ 
499 & 15  & 4   & 3.50e-02 & 1.64e-02                    & 4.93e-02 & 3.53e-01                    \\
599 & 12  & 4   & 2.04e-02 & 1.86e-02                    & 8.39e-02 & 3.94e-01                    \\
 \hline
\end{tabular}
\vspace*{0.3cm}
\caption{Relative global errors for oversampled-GFEM and MS-GFEM, both computed and compared to the FEM overkill solution (relative energy norm error).  Contrasts 100:1 and 1:100. The columns $\omega_1,\,\omega_2$ are the number of shape functions used in the corresponding subdomain. The local spaces $S_{\omega_i}^{n_i}$ are composed of elements obtained by A-harmonically extending boundary hat functions containing a single-node.}
\label{tab:hats_gfem}
\end{table}

\begin{table}[!htb]
\centering
\begin{tabular}{| c | c | c | c | c | c | c |}
\cline{4-7}
\multicolumn{3}{ c }{\rule{0pt}{3ex}}  &  \multicolumn{2}{|p{15ex}|}{\centering 500:1} & \multicolumn{2}{|p{15ex}|}{\centering 1:500}  \\
\hline
\rule[-.3\baselineskip]{-3pt}{3ex}
Hat Widths & $\omega_1$ & $\omega_2$ & GFEM & MS-GFEM & GFEM & MS-GFEM	\\ \hline
\rule[-.3\baselineskip]{-3pt}{3ex}
69  & 109 & 28  & 6.81E-03 & 3.06E-05  & 3.41E-05  & 4.61E-05 \\
99  & 76  & 20  & 1.04E-04 & 2.81E-05  & 1.33E-04  & 4.10E-05 \\
129 & 59  & 15  & 1.31E-04 & 2.75E-05  & 1.49E-04  & 2.08E-04 \\
159 & 48  & 13  & 6.99E-04 & 3.51E-05  & 6.92E-04  & 5.46E-04 \\
189 & 40  & 11  & 6.65E-04 & 1.04E-04  & 2.60E-03  & 1.28E-03 \\
219 & 34  & 9   & 2.70E-03 & 3.40E-04  & 3.92E-03  & 1.16E-02 \\
249 & 30  & 8   & 4.39E-03 & 1.07E-03  & 7.23E-03  & 1.70E-02 \\
279 & 27  & 7   & 5.49E-03 & 4.50E-03  & 1.81E-02  & 2.10E-02 \\
309 & 24  & 7   & 1.29E-02 & 4.74E-03  & 1.03E-02  & 2.80E-02 \\
339 & 22  & 6   & 1.34E-02 & 4.91E-03  & 1.63E-02  & 7.18E-02 \\
369 & 20  & 6   & 1.22E-02 & 5.00E-03  & 1.56E-02  & 1.74E-01 \\
399 & 19  & 5   & 2.68E-02 & 1.57E-02  & 2.34E-02  & 3.39E-01 \\ 
499 & 15  & 4   & 3.54E-02 & 1.65E-02  & 9.79E-02  & 7.59E-01 \\
599 & 12  & 4   & 2.05E-02 & 1.87E-02  & 1.61E-01  & 8.58E-01 \\
 \hline
\end{tabular}
\vspace*{0.3cm}
\caption{Relative global errors for oversampled-GFEM and MS-GFEM, both computed and compared to the FEM overkill solution (relative energy norm error).  Contrasts 500:1 and 1:500. The columns $\omega_1,\,\omega_2$ are the number of shape functions used in the corresponding subdomain.  The local spaces $S_{\omega_i}^{n_i}$ are composed of elements obtained by A-harmonically extending boundary hat functions containing a single-node.}
\label{tab:hats_gfem2}
\end{table}

\begin{table}[!htbp]
\centering
\begin{tabular}{| c | c | c | c | c | c | c |}
\cline{4-7}
\multicolumn{3}{ c }{\rule{0pt}{3ex}}  &  \multicolumn{2}{|p{15ex}|}{\centering 1000:1} & \multicolumn{2}{|p{15ex}|}{\centering 1:1000}  \\
\hline
\rule[-.3\baselineskip]{-3pt}{3ex}
Hat Widths & $\omega_1$ & $\omega_2$ & GFEM & MS-GFEM & GFEM & MS-GFEM	\\ \hline
\rule[-.3\baselineskip]{-3pt}{3ex}
69  & 109 & 28  & 2.99E-05 & 3.09E-05 & 9.71E-04 & 5.62E-05 \\
99  & 76  & 20  & 1.17E-04 & 2.85E-05 & 5.87E-03 & 4.65E-05 \\
129 & 59  & 15  & 1.30E-04 & 2.80E-05 & 1.48E-04 & 2.29E-04 \\
159 & 48  & 13  & 6.95E-04 & 3.49E-05 & 7.71E-04 & 5.64E-04 \\
189 & 40  & 11  & 6.55E-04 & 1.05E-04 & 2.87E-03 & 1.36E-03 \\
219 & 34  & 9   & 2.70E-03 & 3.38E-04 & 4.09E-03 & 1.26E-02 \\
249 & 30  & 8   & 4.39E-03 & 1.06E-03 & 7.91E-03 & 2.24E-02 \\
279 & 27  & 7   & 5.48E-03 & 4.49E-03 & 2.14E-02 & 2.60E-02 \\
309 & 24  & 7   & 1.29E-02 & 4.75E-03 & 1.17E-02 & 3.60E-02 \\
339 & 22  & 6   & 1.34E-02 & 4.91E-03 & 1.94E-02 & 9.93E-02 \\
369 & 20  & 6   & 1.22E-02 & 5.00E-03 & 1.79E-02 & 2.21E-01 \\
399 & 19  & 5   & 2.68E-02 & 1.57E-02 & 3.03E-02 & 5.44E-01 \\ 
499 & 15  & 4   & 3.55E-02 & 1.65E-02 & 1.36E-01 & 1.06E+00 \\
599 & 12  & 4   & 2.05E-02 & 1.87E-02 & 2.19E-01 & 1.20E+00 \\
 \hline
\end{tabular}
\vspace*{0.3cm}
\caption{Relative global errors for oversampled-GFEM and MS-GFEM, both computed and compared to the FEM overkill solution (relative energy norm error).  Contrasts 1000:1 and 1:1000. The columns $\omega_1,\,\omega_2$ are the number of shape functions used in the corresponding subdomain. The local spaces $S_{\omega_i}^{n_i}$ are composed of elements obtained by A-harmonically extending boundary hat functions containing a single-node.}
\label{tab:hats_gfem3}
\end{table}

\section{Iterative solution implementation}
\label{iterative}
The prime motivation behind the MS-GFEM and its contrast independent exponential convergence rate is its potential use in large parallel implementations. As expressed in the introduction the method has two essential parts. A. Parallel construction of local shape functions leading to the exponential approximation of the actual solution in the energy norm. Here each local basis is defined over separate patches and is computed independently on separate processors using local memory. B. Efficient parallel solution of the global system of linear algebraic equations. Recall that we have the partition of unity $\{\phi_i\}_{i=1}^N$ subordinate to the overlapping patches $\omega_i$, with $\Omega=\cup_{i=1}^N\omega_i$ and $\sum_{i=1}^N\phi_i=1$.
We consider the span $V^{m_i}_{\omega_i}$ of $n$-width functions $\xi_j^i,\,j=1,\ldots,m_i$  defined on $\omega_i^\ast$ restricted to $\omega_i$. Set  ${V^{m_i}_{\omega_i}}^\ast=\phi_i V^{m_i}_{\omega_i}$ and the global approximation space is
\begin{equation}\label{eq:globalapprox2}
   V^N =  \sum_{i=1}^N {V^{m_i}_{\omega_i}}^{\ast} .
\end{equation}
We use the partition of unity structure of the global approximation $V^N$ space to iteratively solve the global problem applying the Schwarz method for overlapping domains.  Because each patch overlaps only with its neighbours we can develop an iterative scheme that can be implemented in parallel. Proceeding in this way we arrive at the preconditioned Richardson iteration 
\begin{equation}\label{matrixform}
   \mathrm{B}_N(\mathbf{x}^{n+1}-\mathbf{x}^n)+\mathrm{G}\mathbf{x}^n=\mathbf{r}_N, \, n=0,1,\ldots, x^0\in\mathbb{R}^N,
\end{equation}
where $B_N$ is taken to be an additive two level Schwartz preconditioner \cite{TarekM}, \cite{TandWidland}.  This scheme is implemented in parallel and the iteration process converges exponentially. Here we will assume a fixed number of cores per patch. For this implementation the time for convergence is anticipated to be constant for a fixed ratio between the size of the problem and the number of cores, so it weakly scales with the number of cores used in the computation \cite{DJN}. This iterative method  will be addressed in detail in a forthcoming paper.


\section{Conclusions}
\label{concludsions}
 The MS-GFEM numerical method is shown to be a viable method for the computation of local fields in heterogeneous materials with any contrast between material properties.  A full implementation of the method along with the theoretical foundation of the method were discussed.  A method for dramatically reducing computational costs in computing local trial fields given by the $n$-width functions was presented along with numerical results supporting the method.  We remark that A-harmonic extensions with long wavelength boundary data when restricted to subdomains have energies on these subdomains that decay slower than short wavelength data as the subdomain boundary moves away from where the data is prescribed, see \cite{BLS08} and \cite{LS16}. Motivated by this fact we have proposed to use  A-harmonic extensions of boundary partitions of unity where the boundary hat functions have large support as local bases. Our numerical investigations show that local basis generated this way provide nearly exponential decay in the relative global error as the dimension of the approximation space is increased. Use of these bases provides a further reduction in computational cost. In future work the implementation of MS-GFEM discussed here will be further parallelized for use on large scale computers.  For much larger problems with many patches greater savings from computational costs are anticipated by using coarser hat functions on the boundary to generate local approximation spaces.  We finish by noting that in \cite{SW17} related work has been recently carried out showing the contrast independence of the local spectral basis functions on their domains and that  the global solutions converge exponentially in the $L^2$ norm. There it is found that the local $n$-width basis functions outperform all other contemporary choices of local trial fields constructed by oversampling.
 
\section{Appendix}
\label{App}
For completeness we show how to estimate the local error for patches that share a common boundary with the computational domain.
For patches with Neumann boundary data \eqref{eq:globalneumann} we introduce the space $\mathcal{U}=\{v\in H^1(\omega^*_i):\,v=0 \hbox{ on }\partial{\omega^*_i}\cap\Omega\}$. The local particular solution $\chi_i$ of \eqref{eq:localparticular} with \eqref{eq:globalneumann} solves the equivalent variational problem: $\chi_i\in\mathcal{U}$ and
\begin{equation}\label{eq:neumanncorrect1}
\begin{aligned}
\int_{\omega_i^*}A(x)\nabla\chi_i\cdot\nabla v\,dx=\int_{\omega^*_i}vf\,dx+\int_{\partial\omega^*_i\cap\partial\Omega}v g\,ds,
\end{aligned}
\end{equation}
for all $v$ in $\mathcal{U}$.
Next observe that the actual solution $u$ also solves
\begin{equation}\label{eq:neumanncorrectu}
\begin{aligned}
\int_{\omega_i^*}A(x)\nabla u\cdot\nabla v\,dx=\int_{\omega^*_i}vf\,dx+\int_{\partial\omega^*_i\cap\partial\Omega}v g\,ds,
\end{aligned}
\end{equation}
for all $v$ in $\mathcal{U}$ so
\begin{equation}\label{eq:neumanncorrectequate}
\begin{aligned}
&\int_{\omega_i^*}A(x)\nabla\chi_i\cdot\nabla v\,dx=\int_{\omega_i^*}A(x)\nabla u\cdot\nabla v\,dx\\
\end{aligned}
\end{equation}
for all $v$ in $\mathcal{U}$. On choosing $v=\chi_i$ in \eqref{eq:neumanncorrectequate} and applying Cauchy's inequality we get
\begin{equation}\label{eq:patchlocintestneumann}
\begin{aligned}
&\Vert\chi_i\Vert_{\mathcal{E}(\omega_i^*)}\leq\Vert u\Vert_{\mathcal{E}(\omega_i^*)},
\end{aligned}
\end{equation}
and as before a simple calculation shows
\begin{equation}\label{eq:patchlocintesquadtneumann}
\begin{aligned}
&\Vert u-\chi_i\Vert_{\mathcal{E}(\omega_i^*)}\leq 2\Vert u\Vert_{\mathcal{E}(\omega_i^*)}.
\end{aligned}
\end{equation}
and \eqref{est3} follows.

For Dirichlet boundary data on $\partial\omega_i^*\cap\partial\Omega$ we use linearity and write
$\chi_i=\chi_i^R+\chi_i^D$, where $\chi_i^R\in H^1_0(\omega_i^*)$ solves
\begin{equation}\label{eq:neumanncorrect}
\begin{aligned}
\int_{\omega_i^*}A(x)\nabla\chi_i^R\cdot\nabla v\,dx=\int_{\omega^*_i}vf\,dx,
\end{aligned}
\end{equation}
for all $v$ in $H_0^1(\omega_i^*)$, and $\chi_i^D=u$  on $\partial\omega_i^*\cap\partial\Omega$, where $n\cdot A(x)\nabla \chi_i^D=0$ on $\partial\omega_i^*\cap\Omega$, and solves
\begin{equation}\label{eq:neumanncorrectsplit}
\begin{aligned}
\div(A(x)\nabla\chi_i^D)=0,
\end{aligned}
\end{equation}
in $\omega_i^*$. 
To conclude we will show that
\begin{equation}\label{eq:patchlocintdirichletestdirichletw}
\begin{aligned}
&\Vert \chi_i^D\Vert_{\mathcal{E}(\omega_i^*)}\leq \Vert u-\chi_i^R\Vert_{\mathcal{E}(\omega^*_i)}\leq2\Vert u\Vert_{\mathcal{E}(\omega^*_i)}.
\end{aligned}
\end{equation}
This follows from
\begin{equation}\label{eq:vartrace}
    \begin{aligned}
    &\Vert \chi_i^D\Vert_{\mathcal{E}(\omega_i^*)}^2=\inf\left\{\Vert w\Vert_{\mathcal{E}(\omega_i^*)}^2:\,w\in \mathcal{A}_{ad}\right\}\\
    &\hbox{where   }\mathcal{A}_{ad}=\left\{ w\in H^1(\omega_i^*),\,\div(A(x)\nabla w)=0,\hbox{ in }\omega_i^*,\,\,w=u\hbox{ on }\partial\omega_i^*\cap\partial\Omega\right\}
    \end{aligned}
\end{equation}
To see \eqref{eq:vartrace} write any $w\in \mathcal{A}_{ad}$ as $\chi_i^D+\delta$, $\delta=w-\chi_i^D$ to observe $(\chi_i^D,\delta)_{\mathcal{E}(\omega_i^*)}=0$ so
\begin{equation}
    \begin{aligned}
    \Vert \chi_i^D\Vert^2_{\mathcal{E}(\omega_i^*)}\leq\Vert w\Vert^2_{\mathcal{E}(\omega_i^*)}=\Vert \chi_i^D\Vert^2_{\mathcal{E}(\omega_i^*)}+\Vert \delta\Vert^2_{\mathcal{E}(\omega_i^*)}
    \end{aligned}
\end{equation}
and since $u-\chi_i^R\in\mathcal{A}_{ad}$ the first inequality of \eqref{eq:patchlocintdirichletestdirichletw} follows. Note as before $\Vert\chi_i^R\Vert_{\mathcal{E}(\omega_i^*)}\leq\Vert u\Vert_{\mathcal{E}(\omega_i^*)}$ and together with the triangle inequality this gives the second inequality of \eqref{eq:patchlocintdirichletestdirichletw}.
On collecting results we conclude using the triangle inequality that
\begin{equation}\label{eq:patchlocintdirichletestdirichletchi}
\begin{aligned}
&\Vert \chi_i\Vert_{\mathcal{E}(\omega_i^*)}=\Vert \chi_i^D+\chi_i^R\Vert_{\mathcal{E}(\omega_i^*)}\leq 3\Vert u\Vert_{\mathcal{E}_(\omega^*_i)},
\end{aligned}
\end{equation}
and as before a simple calculation shows
\begin{equation}\label{eq:patchlocintesquadfinal}
\begin{aligned}
&\Vert u-\chi_i\Vert_{\mathcal{E}(\omega_i^*)}\leq 2\sqrt{5}\Vert u\Vert_{\mathcal{E}(\omega_i^*)}.
\end{aligned}
\end{equation}
and \eqref{est4} follows.

For completeness we finish by showing that \eqref{locerror1} implies \eqref{globalerror1}. Our calculation follows Theorem 3.3 and  remarks 3.4 and 3.5 of \cite{BB04}. We expand
\begin{equation}\label{expand1}
\begin{aligned}
\Vert u-u_T\Vert^2_{{\mathcal E}(\Omega)}& =\int_\Omega\,A\nabla(u-u^T)\cdot\nabla(u-u^T)\,dx\\
& =\int_\Omega\,A\nabla(\sum_{i=1}^N \phi_i(u-\xi_i-\chi_i)\cdot\nabla(\sum_{i=1}^N \phi_i(u-\xi_i-\chi_i)\,dx\\
& = \int_\Omega\,A(\sum_{i=1}^N k_i)\cdot(\sum_{i=1}^N k_i) \,dx
\end{aligned}
\end{equation}
with
\begin{equation}\label{expand2}
\begin{aligned}
k_i =[\nabla\phi_i(u-\xi_i-\chi_i)+\phi_i\nabla(u-\xi_i-\chi_i)].
\end{aligned}
\end{equation}
It follows that
\begin{equation}\label{expand3}
\begin{aligned}
\Vert u-u_T\Vert^2_{{\mathcal E}(\Omega)}
& \leq2\int_\Omega\,A(\sum_{i=1}^N \nabla\phi_i(u-\xi_i-\chi_i)\cdot(\sum_{i=1}^N \nabla\phi_i(u-\xi_i-\chi_i))\,dx\\
& +2\int_\Omega\,A(\sum_{i=1}^N \phi_i\nabla(u-\xi_i-\chi_i)\cdot(\sum_{i=1}^N \phi_i\nabla(u-\xi_i-\chi_i))\,dx\\
&\leq 2\kappa\sum_{i=1}^N\int_{\omega_i}\,A\nabla\phi_i(u-\xi_i-\chi_i)\cdot \nabla\phi_i(u-\xi_i-\chi_i))\,dx\\
& +2\kappa\sum_{i=1}^N\int_{\omega_i}\,A \phi_i\nabla(u-\xi_i-\chi_i)\cdot \phi_i\nabla(u-\xi_i-\chi_i))\,dx\\
&\leq C (\kappa\sum_{i=1}^N\Vert u-\xi_i-\chi_i\Vert_{L_A^2({\omega_i})}+\kappa\sum_{i=1}^N\Vert u-\xi_i-\chi_i \Vert_{\mathcal{E}(\omega_i)}),
\end{aligned}
\end{equation}
where for isotropic material properties $A(x)=Ia(x)$ we have $\Vert w\Vert_{L^2_A(\omega_i)}=(\int_{\omega_i}\,a w^2\,dx)^{1/2}$. For interior patches and patches on the boundary with Neumann data the local spaces $
V^{m_i}_{\omega_i}$ contain constant functions and we can choose them so that the weighted Poincar\'e inequality holds
\begin{equation}\label{expand4}
\begin{aligned}
\Vert u-\xi_i-\chi_i\Vert^2_{{L^2_A}(\omega_i)}\leq C\Vert u-\xi_i-\chi_i\Vert^2_{{\mathcal E}(\omega_i)}
\end{aligned}
\end{equation}
for a constant C independent of $u-\xi_i-\chi_i$ and $\omega_i$.
We also have a similar inequality for boundary patches associated with  Dirichlet data as $\xi_i=0$ and $u-\chi_i=0$ for points on $\partial\omega_i\cap\partial\Omega_D$. This observation together with \eqref{expand3}, \eqref{est2}, \eqref{est3} and \eqref{est4} shows that \eqref{locerror1} implies \eqref{globalerror1}.

\end{document}